\documentclass[11pt]{article}
\usepackage{amssymb,amsthm,amsmath}
\usepackage[dvips]{graphicx}
\usepackage{psfrag}
\theoremstyle{plain}
\newtheorem{theorem}{Theorem}
\newtheorem{corollary}{Corollary}
\newtheorem{proposition}{Proposition}[section]
\newtheorem{lemma}[proposition]{Lemma}

\theoremstyle{definition}
\newtheorem{definition}[proposition]{Definition}
\newtheorem{problem}{Problem}
\newtheorem{example}{Example}

\theoremstyle{remark}
\newtheorem{remark}[proposition]{Remark}

\numberwithin{equation}{section}

\def\C{\mathbb{C}}
\def\R{\mathbb{R}}
\def\Z{\mathbb{Z}}
\def\N{\mathbb{N}}
\def\id{\operatorname{id}}
\def\re{\operatorname{Re}}
\def\im{\operatorname{Im}}
\def\hR{\widehat{R}}
\def\hL{\widehat{L}}
\def\hQ{\widehat{Q}}
\def\hla{\hat{\lambda}}

\def\Diff{\mathrm{Diff}_{\mu}^{1}(M)}
\def\Sympl{\mathrm{Sympl}_{\omega}^{1}(M)}
\def\gld{\mathrm{GL}(d)}
\def\gldr{\mathrm{GL}(d,\mathbb{R})}
\def\gldc{\mathrm{GL}(d,\mathbb{C})}

\def\sldr{\mathrm{SL}(d,\mathbb{R})}
\def\sldc{\mathrm{SL}(d,\mathbb{C})}
\def\spdr{\mathrm{Sp}(2q,\mathbb{R})}
\def\uni{\mathrm{U}(q)}

\def\rp{\mathbb{R}\mathrm{P}^{d-1}}
\def\LE{\mathrm{LE}}
\DeclareMathOperator{\vol}{vol}
\DeclareMathOperator{\diam}{diam}
\def\AA{\mathcal{A}}
\def\BB{\mathcal{B}}
\def\CC{\mathcal{C}}
\def\DD{\mathcal{D}}
\def\EE{\mathcal{E}}
\def\LL{\mathcal{L}}

\def\RR{\mathcal{R}}
\def\UU{\mathcal{U}}

\def\mm{\mathbf{m}}
\def\eps{\varepsilon}
\def\ang{\sphericalangle}
\def\minus{\smallsetminus}
\def\meio{{\textstyle \frac 12}}
\def\hatE{\hat{E}}
\def\wp{\mathord{\wedge}^p}  

\begin{document}


\title{The Lyapunov exponents of generic \\
       volume~preserving and symplectic systems}
\author{Jairo Bochi and Marcelo Viana
\footnote{Partially supported by CNPq and Faperj, Brazil. M.V. is
grateful to the hospitality of Coll\`ege de France, Universit\'e
de Paris-Orsay, and Institut de Math\'ematiques de Jussieu.}}
\date{June 2002}
\maketitle

\hfill To Jacob Palis, with friendship and admiration.

\begin{abstract}
We show that the integrated Lyapunov exponents of $C^1$ volume
preserving diffeomorphisms are simultaneously continuous at a
given diffeomorphism only if the corresponding Oseledets splitting
is trivial (all Lyapunov exponents equal to zero) or else
dominated (uniform hyperbolicity in the projective bundle) almost
everywhere.

We deduce a sharp dichotomy for generic volume preserving
diffeomorphisms on any compact manifold: almost every orbit either
is projectively hyperbolic or has all Lyapunov exponents equal to
zero.

Similarly, for a residual subset of all $C^1$ symplectic
diffeomorphisms on any compact manifold, either the diffeomorphism
is Anosov or almost every point has zero as a Lyapunov exponent,
with multiplicity at least $2$.

Finally, given any closed group $G \subset \gld$ that acts
transitively on the projective space, for a residual subset of all
continuous $G$-valued cocycles over any measure preserving
homeomorphism of a compact space, the Oseledets splitting is
either dominated or trivial.
\end{abstract}




\section{Introduction}

Lyapunov exponents describe the asymptotic evolution of a linear
cocycle over a transformation: positive or negative exponents
correspond to exponential growth or decay of the norm,
respectively, whereas vanishing exponents mean lack of exponential
behavior.

In this work we address two basic, a priori unrelated problems.
One is to understand how frequently do Lyapunov exponents vanish
on typical orbits. The other, to analyze the dependence of
Lyapunov exponents as functions of the system. We are especially
interested in dynamical cocycles, i.e. given by the derivatives of
conservative diffeomorphisms, but we discuss the general situation
as well.

Several approaches have been proposed for proving existence of
non-zero Lyapunov exponents. Let us mention
Furstenberg~\cite{Furstenberg}, Herman~\cite{Herman2},
Kotani~\cite{Kotani}, among others. In contrast, we show
here that vanishing Lyapunov exponents are actually very
frequent: {\em for a residual (dense $G_\delta$) subset
of all volume-preserving $C^1$ diffeomorphisms, and for
almost every orbit, all Lyapunov exponents are equal to
zero or else the Oseledets splittings is dominated.\/}
This extends to generic continuous $G$-valued cocycles
over any transformation, for any matrix group $G$ that
acts transitively on the projective space.

Domination, or uniform hyperbolicity in the projective bundle,
means that each Oseledets subspace is more expanded/less
contracted than the next, by a definite \emph{uniform} factor.
This is a very strong property.
In particular, {\em domination implies that the angles between
the Oseledets subspaces are bounded from zero,
and the Oseledets splitting extends to a continuous splitting
on the closure.\/}
For this reason, it can often be excluded a priori:

\begin{example}
Let $f:S^1\to S^1$ be a homeomorphism and $\mu$ be any invariant
ergodic measure with supp $\mu=S^1$.
Let $\mathcal{N}$ be the set of all continuous
$A:S^1\to\operatorname{SL}(2,\R)$ non-homotopic to a constant. For
a residual subset of $\mathcal{N}$, the Lyapunov exponents of the
corresponding cocycle over $(f,\mu)$ are zero. That is because the
cocycle has no invariant continuous subbundle if $A$ is
non-homotopic to a constant.
\end{example}

These results generalize to arbitrary dimension the work of Bochi~\cite{Bochi},
where it was shown that generic area preserving $C^1$ diffeomorphisms
on any compact surface either are uniformly hyperbolic (Anosov) or have
no hyperbolicity at all: both Lyapunov exponents equal to zero almost
everywhere. This fact had been announced by Ma\~n\'e~\cite{Mane,Mane2} in
the early eighties.

The high dimensional setting requires a conceptually different approach.
That is partly because of the difficulty involved in handling several
subbundles, with variable dimensions, and partly because one has to deal
with projectively hyperbolic, instead of uniformly hyperbolic, sets.
The properties of projectively hyperbolic sets are much weaker
(e.g. they need not be robust) and not yet understood.

Our strategy is to analyze the dependence of Lyapunov exponents on
the dynamics. We obtain the following characterization of the
continuity points of Lyapunov exponents in the space of volume
preserving $C^1$ diffeomorphisms on any compact manifold: {\em
they must have all exponents equal to zero or else the Oseledets
splitting must be dominated, over almost every orbit.\/} Similarly
for continuous linear cocycles over any transformation, and in
this setting the necessary condition is known to be sufficient.

The issue of continuous or differentiable dependence of Lyapunov exponents
on the underlying system is subtle, and not well understood.
See Ruelle~\cite{Ruelle} and also Bourgain, Jitomirskaya~\cite{BJ} for a
discussion and further references. We also mention the following simple
application of the result we just stated, in the context of quasi-periodic
Schr\"odinger cocycles:

\begin{example}
Let $f:S^1\to S^1$ be an irrational rotation, $\mu$ be Lebesgue
measure, and $A:S^1\to\operatorname{SL}(2,\R)$ be given by
$$
A=\left(\begin{array}{cc}E-V(\theta) & -1 \\ 1 & 0
\end{array}\right)
$$
for some $E\in\R$ and $V:S^1\to\R$ continuous. Then $A$ is a point
of discontinuity for the Lyapunov exponents, among all continuous
cocycles over $(f,\mu)$, if and only if the exponents are non-zero
and $E$ is in the spectrum of the associated Schr\"odinger
operator. Compare~\cite{BJ}. This is because $E$ is in the
complement of the spectrum if and only if the cocycle is uniformly
hyperbolic, which for $\operatorname{SL}(2,\R)$ cocycles is
equivalent to domination.
\end{example}

We extend the two-dimensional result of Ma\~n\'e--Bochi also in a
different direction, namely to symplectic diffeomorphisms on any
compact symplectic manifold.
Firstly, we prove that continuity points for the Lyapunov exponents
either are uniformly hyperbolic or have at least $2$ Lyapunov
exponents equal to zero at almost every point. Consequently, {\em
generic symplectic $C^1$ diffeomorphisms either are Anosov or have
vanishing Lyapunov exponents with multiplicity at least $2$ at
almost every point.\/}


Topological results in the vein of our present theorems were
obtained by Millionshchikov~\cite{Million}, in the early eighties,
and by Bonatti, D\'\i az, Pujals, Ures~\cite{BDP,DPU}, in their
recent characterization of robust transitivity for
diffeomorphisms. A counterpart of the latter for symplectic maps
had been obtained by Newhouse~\cite{Newhouse} in the seventies,
recently extended by Arnaud~\cite{Arnaud}.

\subsection{Dominated splittings}

Let $M$ be a compact manifold of dimension $d \geq 2$. Let $f:M
\to M$ be a diffeomorphism and $\Gamma \subset M$ be an
$f$-invariant set. Suppose for each $x \in \Gamma$ one is given
non-zero subspaces $E^1_x$ and $E^2_x$ such that $T_x M = E^1_x
\oplus E^2_x$\,, the dimensions of $E_x^1$ and $E_x^2$ are
constant, and the subspaces are $Df$-invariant: $Df_x(E^i_x) =
E^i_{f(x)}$ for all $x\in \Gamma$ and $i=1,2$.
\begin{definition}
Given $m \in \N$, we say that $T_\Gamma\ M=E^1 \oplus E^2$ is an $m$-\emph{dominated}
splitting if for every $x \in \Gamma$ we have
\begin{equation} \label{e.normdomination}
\| Df^m_x |_{E^2_x} \| \,\cdot\,  \| (Df^m_x |_{E^1_x})^ {-1} \|
\leq \meio\,.
\end{equation}
We call $T_\Gamma M=E^1 \oplus E^2$ a \emph{dominated splitting} if it is
$m$-dominated for some $m \in \N$. Then we write $E^1 \succ E^2$.
\end{definition}

Condition~\eqref{e.normdomination} means that, for typical tangent
vectors, their forward iterates converge to $E^1$ and their
backward iterates converge to $E^2$, at uniform exponential rates.
Thus, $E^1$ acts as a global hyperbolic attractor, and $E^2$ acts
as a global hyperbolic repeller, for the dynamics induced by $Df$
on the projective bundle.


More generally, we say that a splitting $T_\Gamma M =E^1\oplus
\cdots \oplus E^k$, into any number of sub-bundles, is dominated
if
$$
E^1\oplus \cdots \oplus E^j \succ E^{j+1} \oplus \cdots \oplus E^k
\quad\text{for every $1 \leq j < k$.}
$$
We say that a splitting $T_\Gamma M = E^1 \oplus \cdots \oplus E^k$,
is \emph{dominated at $x$}, for some point $x\in\Gamma$, if it is
dominated when restricted to the orbit $\{f^n(x);\; n\in \Z \}$ of
$x$.

\subsection{Dichotomy for volume preserving diffeomorphisms}\label{ss.intro volume}

Let $f\in\Diff$. By the theorem of Oseledets~\cite{Oseledets}, for
$\mu$-almost every point $x\in M$, there exists $k(x) \in \N$,
real numbers $\hla_1(f,x) > \cdots > \hla_{k(x)}(f,x)$, and a
splitting $T_x M = E^1_x \oplus \cdots \oplus E^{k(x)}_x$ of the
tangent space at $x$, all depending measurably on the point $x$,
such that
$$
\lim_{n \to \pm \infty} \frac 1n \log \| Df^n_x (v) \| = \hla_j(f,x)
\quad \text{for all $v\in E^j_x \minus \{ 0 \}$}.
$$
Let $\lambda_1(f,x) \geq \lambda_2(f,x) \geq \cdots \geq
\lambda_d(f,x)$ be the numbers $\hla_j(x)$, in non-increasing
order and each repeated with multiplicity $\dim E^j_x$. They are
called the \emph{Lyapunov exponents} of $f$ at $x$. Note that
$\lambda_1(f,x) + \cdots + \lambda_d(f,x) =0$, because $f$
preserves volume. We say that the Oseledets splitting is
\emph{trivial} at $x$ when $k(x)=1$, that is, when all Lyapunov
exponents vanish.

\begin{theorem} \label{t.vol}
There exists a residual set $\RR \subset \Diff$ such that, for
each $f \in \RR$ and $\mu$-almost every $x \in M$, the Oseledets
splitting of $f$ is either trivial or dominated at $x$.
\end{theorem}

For $f \in \RR$ the ambient manifold $M$ splits, up to zero
measure, into disjoint invariant sets $Z$ and $D$
corresponding to trivial splitting and dominated splitting,
respectively. Moreover, $D$ may be written as an increasing union
$D = \cup_{m \in \N} D_m$ of compact $f$-invariant sets, each
admitting a dominated splitting of the tangent bundle.

If $f \in \RR$ is ergodic then either $\mu(Z)=1$ or there is
$m\in\N$ such that $\mu(D_m)=1$. The first case means that all the
Lyapunov exponents vanish almost everywhere. In the second case,
the Oseledets splitting extends continuously to a dominated
splitting of the tangent bundle over the whole ambient manifold $M$.

\begin{example}
Let $f_t : N \to N$, $t\in S^1$, be a smooth family of volume
preserving diffeomorphisms on some compact manifold $N$, such that
$f_t = \id$ for $t$ in some interval $I \subset S^1$, and $f_t$ is
partially hyperbolic for $t$ in another interval $J \subset S^1$.
Such families may be obtained, for instance, using the
construction of partially hyperbolic diffeomorphisms isotopic to
the identity in~\cite{BD}. Then $f: S^1 \times N\to S^1 \times
N$, $f(t,x) = (t, f_t(x))$ is a volume preserving diffeomorphism
for which $D \supset S^1 \times J$ and $Z \supset S^1 \times I$.
\end{example}

Thus, in general we may have $0 < \mu(Z) < 1$.
However, we ignore
whether such examples can be made generic (see also
section~\ref{ss.intro symplectic}):
\begin{problem}
Is there a residual subset of $\Diff$ for which invariant sets
with a dominated splitting have either zero or full measure ?
\end{problem}

Theorem~\ref{t.vol} is a consequence of the following result about
continuity of Lyapunov exponents as functions of the dynamics. For
$j=1,\ldots,d-1$, define
$$
\LE_j (f) = \int_M \left[\lambda_1(f,x) + \cdots + \lambda_j(f,x)
\right] d\mu(x).
$$
It is well-known that the functions $f\in \Diff \mapsto \LE_j(f)$
are upper semi-continuous.
Our next main theorem shows that \emph{lower} semi-continuity is
much more delicate:

\begin{theorem} \label{t.vol.continuity}
Let $f_0 \in \Diff$ be such that the map
$$
f \in \Diff \mapsto \big( \LE_1(f), \ldots, \LE_{d-1} (f) \big) \in \R^{d-1}
$$
is continuous at $f=f_0$. Then for $\mu$-almost every $x\in M$,
the Oseledets splitting of $f_0$ is either dominated or
trivial at $x$.
\end{theorem}

The set of continuity points of a semi-continuous
function on a Baire space is always a residual subset of the space
(see e.g. \cite[\S 31.X]{Kur});
therefore theorem~\ref{t.vol} is an immediate corollary of theorem~\ref{t.vol.continuity}.

\begin{problem}
Is the necessary condition in theorem~\ref{t.vol.continuity} also
sufficient for continuity ?
\end{problem}

Diffeomorphisms with all Lyapunov exponents equal to zero almost
everywhere, or else whose Oseledets splitting extends to a
dominated splitting over the whole manifold, are always continuity
points. Moreover, the answer is affirmative in the context of
linear cocycles, as we shall see.

\subsection{Dichotomy for symplectic diffeomorphisms}\label{ss.intro symplectic}

Now we turn ourselves to symplectic systems. Let $(M^{2q},
\omega)$ be a compact symplectic manifold without boundary. We
denote by $\mu$ the volume measure associated to the volume form
$\omega^q = \omega \wedge \cdots \wedge \omega$. The space
$\Sympl$ of all $C^1$ symplectic diffeomorphisms is a subspace of
$\Diff$. We also fix a Riemannian metric on $M$, the particular
choice being irrelevant for all purposes.

The Lyapunov exponents of symplectic diffeomorphisms have a
symmetry property: $\lambda_j(f,x) = - \lambda_{2q-j}(f,x)$ for
all $1\le j\le q$. In particular, $\lambda_q(x) \geq 0$ and
$\LE_q(f)$ is the integral of the sum of all non-negative
exponents. Consider the splitting
$$
T_x M = E^+_x \oplus E^0_x \oplus E^-_x,
$$
where $E^+_x$, $E^0_x$, and $E^-_x$ are the sums of all Oseledets
spaces associated to positive, zero, and negative Lyapunov
exponents, respectively. Then $\dim E^+_x = \dim E^-_x$ and $\dim
E^0_x$ is even.

\begin{theorem} \label{t.sympl.continuity}
Let $f_0 \in \Sympl$ be such that the map
$$
f \in \Sympl \mapsto \LE_q(f) \in \R
$$
is continuous at $f=f_0$. Then for $\mu$-almost every $x\in M$,
either $\dim E^0_x \geq 2$ or the splitting $T_x M = E^+_x \oplus
E^-_x$ is hyperbolic along the orbit of $x$.
\end{theorem}

In the second alternative, what we actually prove is that the
splitting is dominated at $x$.
This is enough because, for symplectic diffeomorphisms,
dominated splittings into two subspaces of the same dimension
are uniformly hyperbolic.

As in the volume preserving case, the function $f \mapsto
\LE_q(f)$ is continuous on a residual subset $\RR_1$ of $\Sympl$.
Also, we show that there is a residual subset
$\RR_2 \subset \Sympl$ such that for every $f\in\RR_2$ either $f$
is an Anosov diffeomorphism or all its hyperbolic sets have zero
measure. Taking $\RR = \RR_1 \cap \RR_2$, we obtain:

\begin{theorem} \label{t.sympl}
There exists a residual set $\RR \subset \Sympl$ such that every
$f \in \RR$ either is Anosov or has at least two zero Lyapunov
exponents at almost every point.
\end{theorem}

For $d=2$ one recovers the two-dimensional result of Ma\~n\'e-Bochi.

\subsection{Linear cocycles} \label{ss.intro cocycles}

Now we comment on corresponding statements for linear cocycles.
Let $M$ be a compact Hausdorff space, $\mu$ a Borel regular
probability measure, and $f:M\to M$ a homeomorphisms that
preserves $\mu$.

Let $G \subset \gldr$ be a closed group and $C(M,G)$ represent the
space of all continuous maps $M \to G$, endowed with the
$C^0$-topology. To each $A \in C(M,G)$ one associates the linear
cocycle
\begin{equation}\label{eq.cocycletrivialbundle}
F_A : M \times \R^d \to M \times \R^d \ , \quad F(x,v) = (f(x),
A(x) v).
\end{equation}
Oseledets theorem extends to this setting, and so does the concept
of dominated splitting; see sections~\ref{ss.oseledets} and
\ref{ss.basicdominated}.

\begin{theorem}\label{t.cocycle}
Let $G$ be a closed subgroup of $\gldr$ acting transitively on
$\rp$. Then $A_0 \in C(M,G)$ is a point of continuity of
$$
C(M,G)\ni A \mapsto (\LE_1(A), \ldots, \LE_{d-1}(A))\in\R^{d-1}
$$
if and only if the Oseledets splitting of the cocycle $F_A$ at $x$
is either dominated or trivial at $\mu$-almost every $x\in M$.

Consequently, there exists a residual subset $\RR \subset C(M,G)$
such that for every $A\in\RR$ and almost every $x\in X$, the
Oseledets splitting of $F_A$ at $x$ is either trivial or
dominated.
\end{theorem}

The most common matrix groups satisfy the hypothesis of the
theorem, e.g., $\gldr$, $\sldr$, $\spdr$,
as well as $\sldc$, $\gldc$ (which are isomorphic to subgroups of
$\mathrm{GL}(2d,\R)$). Notice that compact groups are not of
interest in this context, because all Lyapunov exponents vanish
identically.

\begin{corollary}
Assume $(f,\mu)$ is ergodic. For any $G$ as in
Theorem~\ref{t.cocycle}, there exists a residual subset $\RR
\subset C(M,G)$ such that every $A\in\RR$ either has all exponents
equal at almost every point, or there exists a dominated splitting
of $M \times \R^d$ which coincides with the Oseledets splitting
almost everywhere.
\end{corollary}



\subsection{Extensions and related problems}

\begin{problem}
For generic smooth families $\R^p\to \Diff$, $\Sympl$, $C(M,G)$,
what can be said of the Lebesgue measure of the subset of
parameters corresponding to zero Lyapunov exponents ?
\end{problem}

\begin{problem}
What are the continuity points of Lyapunov exponents in
$\operatorname{Diff}^{1+r}_\mu(M)$ or $C^r(M,G)$ for $r>0$ ?
\end{problem}

Most of the results stated above were announced in~\cite{IHP}.
Actually, our theorems~\ref{t.sympl.continuity} and~\ref{t.sympl}
do not quite give the full strength of theorem~4 in~\cite{IHP}.
The difficulty is that the symplectic analogue of our construction
of realizable sequences is less satisfactory, unless the subspaces
involved have the same dimension; see remark~\ref{r.solereason}.

\begin{problem}
The Oseledets splitting of generic symplectic $C^1$ diffeomorphisms
is either trivial or partially hyperbolic at almost every point,
\end{problem}

Theorem~\ref{t.cocycle} and the corollary remain true if one
replaces $C(M,G)$ by $L^\infty(M,G)$. We only need $f$ to be an
invertible measure preserving transformation.

\section{Preliminaries} \label{s.preliminaries}

%
%

\subsection{Lyapunov exponents, Oseledets splittings} \label{ss.oseledets}

Let $M$ be a compact Hausdorff space and $\pi:\EE\to M$ be a
continuous finite-dimensional vector bundle endowed with a
continuous Riemann structure. A {\em cocycle\/} over a
homeomorphism $f:M\to M$ is a continuous transformation
$F:\EE\to\EE$ such that $\pi\circ F= f\circ \pi$ and
$F_x:\EE_x\to\EE_{f(x)}$ is a linear isomorphism on each fiber
$\EE_x=\pi^{-1}(x)$. Notice that \eqref{eq.cocycletrivialbundle}
corresponds to the case when the vector bundle is trivial.

\subsubsection{Oseledets theorem}

Let $\mu$ be any $f$-invariant Borel probability measure in $M$. The
theorem of Oseledets~\cite{Oseledets} states that for $\mu$-almost
every point $x$ there exists a splitting
\begin{equation}
\label{eq.oseledets}
\EE_x = E^1_x \oplus \cdots \oplus E^{k(x)}_x\,,
\end{equation}
and real numbers $\hla_1(x) > \cdots > \hla_{k(x)}(x)$ such that
$F_x(E^j_x) = E^j_{f(x)}$ and
$$
\lim_{n\to\pm\infty} \frac 1n \log \|F_x^n (v) \| = \hla_j(x)
$$
for $v \in E_x^j \minus \{0\}$ and $j=1, \dots, k(x)$.
Moreover, if $J_1$ and $J_2$ are any disjoint subsets of
the set of indices $\{1,\dots,k(x)\}$, then
\begin{equation}\label{eq.oseledetsangle}
\lim_{n\to\pm\infty} \frac 1n \log \ang\Big({\bigoplus}_{j \in J_1}
E^j_{f^n(x)}, {\bigoplus}_{j \in J_2} E^j_{f^n(x)}\Big) = 0.
\end{equation}

Let $ \lambda_1(x) \geq \lambda_2(x) \geq \cdots \geq\lambda_d(x)$
be the numbers $\hla_j(x)$, each repeated with multiplicity $\dim
E^j_x$ and written in non-increasing order.
When the dependence on $F$ matters, we write $\lambda_i(F,x)=\lambda_i(x)$.
In the case when $F=Df$, we write $\lambda_i(f,x)=\lambda_i(F,x)=\lambda_i(x)$.

\subsubsection{Exterior products}

Given a vector space $V$ and a positive integer $p$,
let $\wp(V)$ be the $p$:th exterior power of $V$.
This is a vector space of dimension $\binom{d}{p}$,
whose elements are called \emph{$p$-vectors}.
It is generated by the $p$-vectors of the form
$v_1 \wedge \dots \wedge v_p$ with $v_j \in V$,
called the \emph{decomposable $p$-vectors}.
A linear map $L:V \to W$ induces a linear map
$\wp(L): \wp(V) \to \wp(W)$ such that
$$
\wp(L) (v_1 \wedge \dots \wedge v_p) = L(v_1) \wedge \dots \wedge L(v_p).
$$
If $V$ has an inner product, then we always endow
$\wp(V)$ with the inner product such that
$\| v_1 \wedge \dots \wedge v_p \|$ equals the $p$-dimensional volume of the
parallelepiped spanned by $v_1$, \ldots, $v_p$.
See~\cite[section~3.2.3]{LArnold}.

More generally, there is a vector bundle $\wp(\EE)$,
with fibers $\wp(\EE_x)$, associated to $\EE$,
and there is a vector bundle automorphism $\wp(F)$,
associated to $F$.
If the vector bundle $\EE$ is endowed with a continuous inner product,
then $\wp(\EE)$ also is.
The Oseledets data of $\wp(F)$ can be obtained
from that of $F$, as shown by the proposition below.
For a proof, see~\cite[theorem~5.3.1]{LArnold}.

\begin{proposition} \label{p.oseledets exterior}
The Lyapunov exponents (with multiplicity) of the automorphism $\wp(F)$ at
a point $x$ are the numbers
$$
\lambda_{i_1}(x) + \dots + \lambda_{i_p}(x), \quad
\text{where $1 \leq i_1 < \dots < i_p \leq d$.}
$$
Let $\{e_1(x), \dots, e_d(x)\}$ be a basis of $\EE_x$ such that
$$
e_i(x) \in E^\ell_x
\quad \text{for $\dim E^1_x + \dots + \dim E^{\ell-1}_x
< i \leq \dim E^1_x + \dots + \dim E^\ell_x$.}
$$
Then the Oseledets space $E_x^{j,\wedge p}$ of $\wp(F)$ corresponding to the
Lyapunov exponent $\hla_j(x)$ is the sub-space of $\wp(\EE_x)$ generated by
$$
e_{i_1} \wedge \dots \wedge e_{i_p}, \quad
\text{with $1 \leq i_1 < \dots < i_p \leq d$ and
$\lambda_{i_1}(x) + \dots + \lambda_{i_p}(x) = \hla_j(x)$.}
$$
\end{proposition}

\subsubsection{Semi-continuity of integrated exponents} \label{sss.semicont}

Let us indicate $\Lambda_p(F,x) = \lambda_1(F,x) + \cdots + \lambda_p(F,x)$,
for $p = 1,\ldots, d-1$.
We define the \emph{integrated Lyapunov exponent}
$$
\LE_p (F) = \int_M \Lambda_p(F, x) \, d\mu(x) .
$$
More generally, if $\Gamma \subset M$ is a measurable $f$-invariant subset,
we define
$$
\LE_p (F, \Gamma) = \int_\Gamma \Lambda_p(F, x) \, d\mu(x) .
$$

By proposition~\ref{p.oseledets exterior},
$\Lambda_p(F,x) = \lambda_1(\wp F, x)$ and so
$\LE_p(F,\Gamma) = \LE_1(\wp(F), \Gamma)$.
When $F=Df$, we write $\Lambda_p(f,x)=\Lambda_i(F,x)$ and
$\LE_p (f, \Gamma)=\LE_p (F,\Gamma)$.

\begin{proposition} \label{p.formula}
If $\Gamma \subset M$ is a measurable $f$-invariant subset then
$$
\LE_p(F, \Gamma) = \inf_{n\geq 1} \frac 1n \int_\Gamma \log
\|\wp(F_x^n)\| \, d\mu(x).
$$
\end{proposition}

\begin{proof}
The sequence $a_n = \int_\Gamma \log \|\wp(F_x^n)\| \, d\mu$ is subadditive
($a_{n+m} \leq a_n + a_m$),
therefore $\lim \frac{a_n}{n} = \inf \frac{a_n}{n}$.
\end{proof}

As a consequence of proposition~\ref{p.formula}, the map
$f\in \Diff \mapsto \LE_p(f)$ is upper semi-continuous,
as mentioned in the introduction.

\subsection{Dominated splittings} \label{ss.basicdominated}

Let $\Gamma\subset M$ be an $f$-invariant set. A splitting
$\EE_\Gamma=E^1\oplus E^2$ is \emph{dominated} for $F$ if it is
$F$-invariant, the dimensions of $E^i_x$ are constant on $\Gamma$,
and there exists $m\in\N$ such that, for every $x\in\Gamma$,
\begin{equation} \label{e.normdominationagain}
\frac{\|F^m_x |_{E^2_x} \|}{ \mm(F^m_x |_{E^1_x})} \leq \frac{1}{2} \,.
\end{equation}
We denote $\mm(L) = \| L^{-1} \|^{-1}$ the \emph{co-norm} of a
linear isomorphism $L$.
The dimension of the space $E^1$ is called the \emph{index} of the splitting.

A few elementary properties of dominated decompositions follow.
The proofs are left to the reader.

\paragraph{Transversality:} If $\EE_\Gamma = E^1 \oplus E^2$ is a
dominated splitting then the angle $\ang(E^1_x,E^2_x)$ is bounded
away from zero, over all $x\in\Gamma$.

\paragraph{Uniqueness:} If $\EE_\Gamma = E^1 \oplus E^2$ and
$\EE_\Gamma = \hatE^1 \oplus \hatE^2$ are dominated decompositions with
$\dim E^i=\dim \hatE^i$ then $E^i=\hatE^i$ for $i=1, 2$.

\paragraph{Continuity:} A dominated splitting $\EE_\Gamma=E^1\oplus
E^2$ is continuous, and extends continuously to a dominated
splitting over the closure of $\Gamma$.

\subsection{Dominance and hyperbolicity for symplectic maps}

Let $(V,\omega)$ be a symplectic vector space of dimension $2q$.
Given a subspace $W \subset V$,
its \emph{symplectic orthogonal} is the space (of dimension $2q - \dim W$)
$$
W^\omega = \{w \in W;\; \omega(v,w) = 0 \text{ for all } v \in V \}.
$$
The subspace $W$ is called \emph{symplectic} if $W^\omega \cap W = \{ 0\}$,
that is, $\omega|_{W \times W}$ is a non-degenerate form.
$W$ is called \emph{isotropic} if $W \subset W^\omega$, that is,
$\omega|_{W \times W} \equiv 0$.
The subspace $W$ is called \emph{Lagrangian} if $W = W^\omega$, that is,
it is isotropic and $\dim W = q$.

\medskip

Now let $(M, \omega)$ be a symplectic manifold of dimension $d=2q$.
We also fix in $M$ a Riemannian structure.
For each $x\in M$, let $J_x: T_x M \to T_x M$ be the anti-symmetric isomorphism
defined by $\omega(v,w)=\langle J_x v,w \rangle$ for all
$v,w \in T_x M$.
Denote
\begin{equation} \label{e.def C omega}
C_\omega = \sup_{x \in M} \| J_x^{\pm 1} \|.
\end{equation}
In particular, we have
\begin{equation} \label{e.C omega}
|\omega (v,w)|  \leq C_\omega \|v\| \, \|w\| \quad
\text{for all $v,w \in T_x M$.}
\end{equation}

\begin{lemma} \label{l.sympl 1}
If $E$, $F \subset T_x M$ are two Lagrangian subspaces
with $E \cap F = \{0\}$
and $\alpha = \ang(E,F)$ then:
\begin{enumerate}
\item For every $v \in E\minus \{0\}$ there exists $w \in F\minus \{0\}$
such that
$$
|\omega (v,w)|  \geq C_\omega^{-1} \sin \alpha \, \|v\| \, \|w\|.
$$

\item If $S: T_x M \to T_y M$ is any symplectic linear map
and $\beta = \ang(S(E),S(F))$ then
$$
C_\omega^{-2} \sin \alpha
\leq \mm(S|_E) \, \| S|_F \| \leq
C_\omega^2 (\sin \beta)^{-1}.
$$
\end{enumerate}
\end{lemma}

\begin{proof}
To prove part 1, let $p: T_x M \to F$ be the projection parallel to $E$.
Given a non-zero $v \in E$, take $w = p(J_x v)$.
Since $E$ is isotropic,
$\omega (v,w) = \omega (v, J_x v) = \|J_x v\|^2
\geq C_\omega^{-1} \|v\| \, \|J_x v \|$.
Also
$\|w\| \leq \|p\| \, \|J_x v\|$ and $\|p\| = 1/\sin\alpha$,
so the claim follows.

To prove part 2, take a non-zero $v\in E$ such that
$\|Sv\|/\|v\| = \mm(S|_E)$ and let $w$ be given by part 1.
Then
$$
C_\omega^{-1} \sin \alpha \, \|v\| \, \|w\| \leq
|\omega (v,w)|  =
|\omega (Sv,Sw)|  \leq
C_\omega \|Sv\| \, \|Sw\|.
$$
Thus $\mm(S|_E)\, \|Sw \| / \|w\| \geq C_\omega^{-2} \sin \alpha$,
proving the lower inequality in part 2.
The upper inequality follows from the lower one applied to $S(F)$,
$S(E)$ and $S^{-1}$ in the place of $E$, $F$, and $S$, respectively.
\end{proof}

\begin{lemma} \label{l.sympl 2}
Let $f \in \Sympl$, and let $x$ be a regular point.
Assume that $\lambda_q(f,x)>0$, that is, there are no zero exponents.
Let $E^+_x$ and $E^-_x$ be the sum of all Oseledets subspaces associated
to positive and to negative Lyapunov exponents, respectively. Then
\begin{enumerate}
\item The subspaces $E^+_x$ and $E^-_x$ are Lagrangian.
\item If the splitting $E^+ \oplus E^-$ is dominated at $x$ then
$E^+$ is uniformly expanding and $E^-$ is uniformly contracting
along the orbit of $x$.
\end{enumerate}
\end{lemma}

\begin{proof}
To prove part 1, we only have to show that the spaces $E^+_x$ and $E^-_x$ are isotropic.
Take vectors  $v_1,v_2\in E^-_x$.
Take $\eps>0$ with $\eps < \lambda_q(f,x)$.
For every large $n$ and $i=1,2$, we have
$\| Df_x^n v_i\| \leq e^{-n \eps} \|v_i\|$.
Hence, by \eqref{e.C omega},
$$
|\omega (v_1, v_2)| = |\omega (Df^n_x v_1, Df^n_x  v_2)| \leq
C_\omega e^{-2n \eps} \|v_1\| \, \|v_2\|,
$$
that is, $\omega(v_1,v_2)=0$.
A similar argument, iterating backwards, gives that $E^+_x$ is isotropic.

Now assume that $E^+ \succ E^-$ at $x$.
Let $\alpha>0$ be a lower bound for $\ang(E^+,E^-)$ along the orbit of $x$,
and let $C = C_\omega^2 (\sin \alpha)^{-1}$.
By domination, there exists $m \in \N$ be such that
$$
\frac{\|Df^m_{f^n(x)}|_{E^-} \|}{\mm(Df^m_{f^n(x)}|_{E^+} )} < \frac{1}{4C},
\quad \text{for all $n \in \Z$.}
$$
By part 2 of lemma~\ref{l.sympl 2}, we have
$C^{-1} \leq \mm( Df^m_{f^n(x)}|_{E^+}) \, \|Df^m_{f^n(x)}|_{E^-}\|  \leq C$.
Therefore
$$
\mm( Df^m_{f^n(x)}|_{E^+}) > 2 \quad\text{and}\quad
\|Df^k_{f^n(x)}|_{E^-}\| < \meio \quad \text{for all $n \in \Z$.}
$$
This proves part 2.
\end{proof}

\begin{remark} \label{r.partial hyperbolicity}
More generally, existence of a dominated splitting implies partial
hyperbolicity:
\emph{If $E \oplus \widehat{F}$ is a dominated splitting, with
$\dim E \leq \dim \widehat{F}$, then $\widehat{F}$ splits invariantly as
$\widehat{F} = C \oplus F$, with $\dim F = \dim E$. Moreover,
$E$ is uniformly expanding and $F$ is uniformly contracting.}
This fact was pointed out by Ma\~n\'e in \cite{Mane2}. A proof in
dimension $4$ was given recently by Arnaud~\cite{Arnaud}.
Since the present paper does not use this result, we omit the proof.
\end{remark}

\subsection{Angle estimation tools}\label{ss.angle tools}

Here we collect a few useful facts from elementary linear algebra.
We begin by noting that, given any one-dimensional subspaces $A$,
$B$, and $C$ of $\R^d$, then
\begin{multline*}
\sin \ang(A, B) \, \sin \ang(A + B, C) =
\sin \ang(C, A) \, \sin \ang(C + A, B) \\
= \sin \ang(B, C) \, \sin \ang(B + C, A).
\end{multline*}
Indeed, this quantity is the $3$-dimensional volume of the
parallepipid with unit edges in the directions $A$, $B$ and $C$.
As a corollary, we get:

\begin{lemma} \label{l.abc}
Let $A$, $B$ and $C$ be subspaces (of any dimension) of $\R^d$. Then
$$
\sin \ang (A, B + C) \geq \sin \ang (A,B) \, \sin \ang (A + B, C).
$$
\end{lemma}

Let $v$, $w$ be non-zero vectors.
For any $\alpha \in \R$,
$\| v + \alpha w \| \geq \| v \| \sin \ang(v,w)$, with
equality when $\alpha = \langle v,w \rangle / \| w \|^2$.
Given $L \in \gldr$,
let $\beta  = \langle Lv,Lw \rangle / \| Lw \|^2$ and
$z = v + \beta w$.
By the previous remark, $\| z \| \geq \| v \| \sin \ang(v,w)$ and
$\| Lz \| = \| Lv \| \sin \ang(Lv,Lw)$.
Therefore
\begin{equation} \label{e.useful}
\sin \ang(Lv,Lw) = \frac{\| Lz \|}{\| Lv \|}
\geq \frac{\mm(L) \| v \|} {\| Lv \|} \, \sin \ang(v,w).
\end{equation}
As a consequence of~\eqref{e.useful}, we have:

\begin{lemma}\label{l.quase conforme}
Let $L: \R^d \to \R^d$ be a linear map and
let $v$, $w$ be non-zero vectors. Then
$$
\frac{\mm(L)}{\| L \|} \leq
\frac{\sin \ang(Lv,Lw)}{\sin \ang(v,w)} \leq
\frac{\| L \|}{\mm(L)}.
$$
\end{lemma}


Thus $\|L\| / \mm(L)$ measures
how much angles can be distorted by $L$.
At last, we give a bound for this quantity when $d=2$.

\begin{lemma} \label{l.4}
Let $L:\R^2 \to \R^2$ be an invertible linear map and let $v,w \in \R^2$ be
linearly independent unit vectors.
Then
$$
\frac{\|L\|}{\mm(L)} \leq
4 \max \left\{\frac{\|Lv\|}{\|Lw\|} , \frac{\|Lw\|}{\|Lv\|}\right\}
\, \frac{1}{\sin\ang(v,w)}\,\frac{1}{\sin\ang(Lv,Lw)}.
$$
\end{lemma}

\begin{proof}
We may assume that $L$ is not conformal, for in the conformal case
the left hand side is $1$ and the inequality is obvious. Let $\R
s$ be the direction most contracted by $L$, and let $\theta$,
$\phi \in [0, \pi]$ be the angles that the directions $\R v$ and
$\R w$, respectively, make with $\R s$. Suppose that $\|Lv\| \geq
\| Lw \|$. Then $\phi \leq \theta$ and so $\ang(v,w) \leq
2\theta$. Hence
$$
\| Lv \| \geq \| L \| \sin \theta \geq \meio \| L \| \sin 2\theta \geq
\meio \| L \| \sin \ang(v,w).
$$
Moreover, $\lvert \det L \rvert = \mm (L) \|L \|$ and
$$
\| Lv \| \| Lw \| \sin \ang(Lv,Lw) = \lvert \det L \rvert \sin
\ang(v,w).
$$
The claim is an easy consequence of these relations.
\end{proof}

\subsection{Coordinates, metrics, neighborhoods} \label{ss.basic}


Let $(M,\omega)$ be a symplectic manifold of dimension $d=2q\ge 2$.
According to Darboux's theorem, there exists an atlas
$\AA^*=\{\varphi_i:V_i^*\to\R^d\}$ of {\em canonical\/} local
coordinates, that is, such that
$$
(\varphi_i)_*\omega = dx_1 \wedge dx_2 + \cdots + dx_{2q-1} \wedge dx_{2q}
$$
for all $i$. Similarly, cf. \cite[Lemma~2]{Moser}, given any
volume structure $\beta$ on a $d$-dimensional manifold $M$, one
can find an atlas $\AA^*=\{\varphi_i:V_i^*\to\R^d\}$ consisting of
charts $\varphi_i$ such that
$$
(\varphi_i)_*\beta = dx_1 \wedge \cdots \wedge dx_d \,.
$$

In either case, assuming $M$ is compact one may choose $\AA^*$ finite.
Moreover, we may always choose $\AA^*$ so that every $V_i^*$
contains the closure of an open set $V_i$, such that the restrictions
$\varphi_i:V_i\to\R^d$ still form an atlas of $M$.
The latter will be denoted $\AA$. Let $\AA^*$ and $\AA$ be fixed once
and for all.

\smallskip

By compactness, there exists $r_0>0$ such that for each $x \in M$,
there exists $i(x)$ such that the Riemannian ball of radius $r_0$
around $x$ is contained in $V_{i(x)}$.
For definiteness, we choose $i(x)$ smallest with this property.
For technical convenience, when dealing with the point $x$ we express
our estimates in terms of the Riemannian metric $\|\cdot\|=\|\cdot\|_x$
defined on that ball of radius $r_0$ by $\|v\|=\|D\varphi_{i(x)}v\|$.
Observe that these Riemannian metrics are (uniformly) equivalent to
the original one on $M$, and so there is no inconvenience in replacing
one by the other.

We may also view any linear map $A: T_{x_1} M \to T_{x_2} M$ as
acting on $\R^d$, using local charts $\varphi_{i(x_1)}$ and
$\varphi_{i(x_2)}$. This permits us to speak of the distance $\|A
- B \|$ between $A$ and another linear map $B : T_{x_3}M \to
T_{x_4}M$ whose base points are different:
$$
\|A - B \| = \| D_2 A D_1^{-1}
- D_4 B D_3^{-1} \|, \quad \text{where $D_j =
(D\varphi_{i(x_j)})_{x_j}$.}
$$

For $x\in M$ and $r>0$ small (relative to $r_0$), $B_r(x)$ will denote
the ball of radius $r$ around $x$ relative to the new metric.
In other words,
$B_r (x) = \varphi_{i(x)}^{-1} \big( B(\varphi_{i(x)}(x),r) \big)$.
We assume that $r$ is small enough so that the closure of
$B_r(x)$ is contained in $V_{i(x)}^{\ast}$.

\begin{definition} \label{d.basic}
Let $\eps_0 > 0$. The \emph{$\eps_0$-basic neighborhood} $\UU
(\mathrm{id}, \eps_0)$ of the identity in $\Diff$, or in $\Sympl$,
is the set $\UU(\id,\eps_0)$ of all $h\in\Diff$, or $h\in\Sympl$,
such that $h^{\pm 1} (\overline{V}_i)\subset V_i^*$  for each $i$
and
$$
h(x) \in B(x,\eps_0) \quad \text{and} \quad \| Dh_{x} - I  \| < \eps_0
\quad \text{for every $x \in M$.}
$$
For a general $f\in\Diff$, or $f\in\Sympl$, the $\eps_0$-basic
neighborhood $\UU(f,\eps_0)$ is defined by: $g\in\UU(f,\eps_0)$ if
and only if $f^{-1} \circ g \in \UU(\id,\eps_0)$ or $g
\circ f^{-1} \in \UU(\id,\eps_0)$.
\end{definition}

\subsection{Realizable sequences} \label{ss.realizable}

The following notion, introduced in~\cite{Bochi}, is crucial to
the proofs of theorems~\ref{t.vol} through~\ref{t.sympl}. It
captures the idea of sequence of linear transformations that can
be (almost) realized {\em on subsets with large relative
measure\/} as tangent maps of diffeomorphisms close to the
original one.

\begin{definition} \label{d.sr}
Given $f\in \Diff$ or $f\in\Sympl$, constants $\eps_0 >0$,
and $0<\kappa<1$, and a non-periodic point $x \in M$,
we call a sequence of linear maps (volume preserving or symplectic)
$$
T_x M \xrightarrow{L_0} T_{fx} M \xrightarrow{L_1} \cdots
\xrightarrow{L_{n-1}} T_{f^n x} M
$$
an $(\eps_0, \kappa)$\emph{-realizable sequence of length $n$ at
$x$\/} if the following holds:

\smallskip

For every $\gamma >0$ there is $r>0$ such that the iterates
$f^j(\overline{B}_r(x))$ are two-by-two disjoint for $0 \leq j
\leq n$, and given any non-empty open set $U\subset B_{r}(x)$,
there are $g\in \UU(f, \eps_0)$ and a measurable set $K\subset U$
such that
\begin{itemize}
\item[(i)] $g$ equals $f$ outside the disjoint union
$\bigsqcup_{j=0}^{n-1}f^{j}(\overline{U})$;

\item[(ii)] $\mu (K) > (1-\kappa) \mu (U)$;

\item[(iii)] if $y\in K$ then $\left\| Dg_{g^{j}y}-L_{j}\right\| <\gamma$
for every $0\leq j \leq n-1$.
\end{itemize}
\end{definition}


Some basic properties of realizable sequences are collected in the
following

\begin{lemma}\label{l.basicproperties}
Let $f\in \Diff$ or $f \in \Sympl$, $x\in M$ not periodic and $n \in \N$.
\begin{enumerate}
\item The sequence $\{Df_x,\ldots ,Df_{f^{n-1}(x)}\}$ is
$(\eps_0, \kappa)$-realizable for every $\eps_0$ and $\kappa$ (we
call this a \emph{trivial} realizable sequence).

\item Let $\kappa_1 ,\kappa_2 \in (0,1)$ be such that $\kappa=\kappa_1+\kappa_2<1$.
If\/ $\{L_{0},\ldots ,L_{n-1}\}$ is $(\eps_0,
\kappa_1)$-realizable at $x$, and $\{L_{n},\ldots,L_{n+m-1}\}$ is
$(\eps_0, \kappa_2)$-realizable at $f^n(x)$, then $\{L_{0},\ldots
,L_{n+m-1}\}$ is $(\eps_0, \kappa)$-realizable at $x$.

\item If $\{L_{0},\ldots ,L_{n-1}\}$ is $(\eps_0, \kappa)$-realizable
at $x$, then $\{L_{n-1}^{-1},\ldots ,L_{0}^{-1}\}$ is an $(\eps_0,
\kappa)$-realizable sequence at $f^n(x)$ for the diffeomorphism
$f^{-1}$.

\end{enumerate}
\end{lemma}

\begin{proof}
The first claim is obvious. For the second one, fix $\gamma>0$.
Let $r_1$ be the radius associated to the $(\eps_0,
\kappa_1)$-realizable sequence, and $r_2$ be the radius associated
to the $(\eps_0,\kappa_2)$-realizable sequence. Fix $0<r<r_1$ such
that $f^n(B_r(x)) \subset B(f^n(x),r_2)$. Then the
$f^j(\overline{B}_r(x))$ are two-by-two disjoint for $0\le j \le
n+m$. Given an open set $U \subset B_r(x)$, the realizability of
the first sequence gives us a diffeomorphism $g_1 \in
\UU(f,\eps_0)$ and a measurable set $K_1 \subset U$. Analogously,
for the open set $f^n(U) \subset B(f^n(x),r_2)$ we find $g_2 \in
\UU(f,\eps_0)$ and a measurable set $K_2 \subset f^n(U)$. Then
define a diffeomorphism $g$ as $g=g_1$ inside $U\cup\cdots\cup
f^{n-1}(U)$ and $g=g_2$ inside $f^n(U)\cup\cdots\cup
f^{n+m-1}(U)$, with $g=f$ elsewhere. Consider also $K=K_1 \cap
g^{-n}(K_2)$. Using that $g$ preserves volume, one checks that $g$
and $K$ satisfy the conditions in definition~\ref{d.sr}. For claim
3, notice that $\UU(f, \eps_0) = \UU(f^{-1}, \eps_0)$.
\end{proof}

The next lemma makes it simpler to verify that a sequence is realizable:
we only have to check the conditions for certain open sets $U\subset B_r(x)$.

\begin{definition}
A family of open sets $\{W_\alpha\}$ in $\R^d$ is a {\em Vitali
covering\/} of $W=\cup_\alpha W_\alpha$ if there is $C>1$ and for
every $y \in W$, there are sequences of sets $W_{\alpha_n}\ni y$
and positive numbers $s_n \to 0$ such that
$$
B_{s_n}(y)\subset W_{\alpha_n} \subset B_{Cs_n}(y)
\quad \text{for all $n\in\N$.}
$$
A family of subsets $\{U_\alpha\}$ of $M$ is a {\em Vitali
covering\/} of $U=\cup_\alpha U_\alpha$ if each $U_\alpha$ is
contained in the domain of some chart $\varphi_{i(\alpha)}$ in the
atlas $\AA$, and the images $\{\varphi_{i(\alpha)}(U_\alpha)\}$
form a Vitali covering of $W=\varphi(U)$, in the previous sense.
\end{definition}

\begin{lemma} \label{l.simplifying}
Let $f\in \Diff$ or $f\in\Sympl$, and let $\eps_0>0$ and $\kappa>0$.
Consider any sequence
$L_j:T_{f^j(x)}M \to T_{f^{j+1}(x)}$, $0 \le j \le n-1$ of linear
maps at a non-periodic point $x$, and let $\varphi:V\to\R^d$ be a
chart in the atlas $\AA$, with $V \ni x$. Assume the conditions in
definition~\ref{d.sr} are valid for every element of some Vitali
covering $\{U_\alpha\}$ of $B_r(x)$. Then the sequence $L_j$ is
$(\eps_0,\kappa)$-realizable.
\end{lemma}

\begin{proof}
Let $U$ be an arbitrary open subset of $B_r(x)$. By Vitali's
covering lemma (see~\cite{McShane}), there is a countable family
of two-by-two disjoint sets $U_\alpha$ covering $U$ up to a zero
Lebesgue measure subset. Thus we can find a finite family of
$U_\alpha$ with disjoint closures and such that $\mu
\left(U-\bigsqcup_{\alpha} U_\alpha\right) $ is as small as we
please. For each $U_\alpha$ there are, by
hypothesis, a perturbation $g_\alpha\in \UU(f, \eps_0)$ and a
measurable set $K_\alpha \subset U_\alpha$ with the properties
(i)-(iii) of definition~\ref{d.sr}. Let $K=\bigcup K_\alpha$ and
define $g$ as being equal to $g_\alpha$ on each $f^j(U_\alpha)$
with $0\le j\le n-1$. Then $g\in \UU(f, \eps_0)$ and the pair
$(g,K)$ have the properties required by definition~\ref{d.sr}.
\end{proof}

\section{Geometric consequences of non-dominance} \label{s.p.geom}

The aim of this section is to prove the following key result, from
which we shall deduce theorem~\ref{t.vol.continuity} in
section~\ref{s.t.vol}:

\begin{proposition} \label{p.geom}
Given $f \in \Diff$, $\eps_0>0$ and $0< \kappa <1$,
if $m \in \N$ is sufficiently large then the following holds:
Let $y \in M$ be a
non-periodic point and suppose one is given a non-trivial
splitting $T_y M = E \oplus F$ such that
$$
\frac{\| Df^m_y|_F\|}{\mm(Df^m_y|_E)} \geq \frac 12\,.
$$
Then there exists an $(\eps_0, \kappa)$-realizable sequence
$\{L_{0}, \ldots, L_{m-1}\}$ at $y$ of length $m$ and there are
non-zero vectors $v \in E$ and $w\in Df_y^m(F)$ such that
$$
L_{m-1} \cdots L_{0} (v) = w.
$$
\end{proposition}

\subsection{Nested rotations}

Here we present some tools for the construction of realizable
sequences. The first one yields sequences of length $1$:

\begin{lemma} \label{l.length 1}
Given $f\in \Diff$, $\eps_0>0$, $\kappa>0$, there exists $\eps>0$
with the following properties:

Suppose we are given a non-periodic point $x \in M$, a splitting
$\R^d = X \oplus Y$ with $X \perp Y$ and $\dim Y = 2$, and a
elliptic linear map $\hR: Y \to Y$ with $\| \hR - I
\|<\eps$. Consider the linear map $R: T_x M \to T_x M$ given by
$R(u + v ) = u + \hR(v)$, for $u\in X$, $v\in Y$. Then $\{Df_x
R\}$ is an $(\eps_0,\kappa)$-realizable sequence of length $1$ at~$x$
and $\{R \, Df_{f^{-1}(x)} \}$ is an $(\eps_0,\kappa)$-realizable
sequence of length $1$ at the point~$f^{-1}(x)$.
\end{lemma}

We also need to construct long realizable sequences. Part 2 of
lemma~\ref{l.basicproperties} provides a way to do this, by
concatenation of shorter sequences. However, simple concatenation
is far too crude for our purposes because it worsens $\kappa$:
the relative measure of the set where the sequence can be (almost)
realized decreases when the sequence increases. This problem is
overcome by lemma~\ref{l.nested} below, which allows us to obtain
certain non-trivial realizable sequences with arbitrary length
while keeping $\kappa$ controlled.

In short terms, we do concatenate several length $1$ sequences, of
the type given by lemma~\ref{l.length 1}, but we also impose that
the supports of successive perturbations be mapped one to the
other. More precisely, there is a domain $\CC_0$ invariant under
the sequence, in the sense that $L_{j-1} \cdots L_0 (\CC_0) =
Df^j_x (\CC_0)$ for all $j$. Following~\cite{Bochi}, where a
similar notion was introduced for the $2$-dimensional setting, we
call such $L_j$ \emph{nested rotations}. When $d>2$ the domain
$\CC_0$ is not compact, indeed it is the product $\CC_0 = X_0
\oplus \BB_0$ of a codimension~$2$ subspace $X_0$ by an ellipse
$\BB_0 \subset X_0^{\perp}$.

\smallskip

Let us fix some terminology to be used in the sequel. If $E$ is a
vector space with an inner product and $F$ is a subspace of $E$,
we endow the quotient space $E/F$ with the inner product that
makes $v \in F^\perp \mapsto (v+F) \in E/F$ an isometry. If $E'$
is another vector space, any linear map $L: E \to E'$ induces a
linear map $L/F : E/F \to E'/F'$, where $F'=L(F)$. If $E'$ has an
inner product, then we indicate by $\| L/F \|$ the usual operator
norm.

\begin{lemma} \label{l.nested}
Given $f\in \Diff$, $\eps_0>0$, $\kappa>0$, there exists $\eps>0$
with the following properties: Suppose we are given a non-periodic
point $x \in M$ and, for $j=0, 1, \ldots, n-1$,
\begin{itemize}
\item codimension $2$ spaces $X_j \subset T_{f^j(x)} M$
      such that $X_j= Df^j_x (X_0)$;
\item ellipses $\BB_j \subset (T_{f^j(x)} M) / X_j$ centered at
zero with $\BB_j = (Df^j_x / X_0) (\BB_0)$.
\item linear maps $\hR_j: (T_{f^j(x)} M) / X_j \to (T_{f^j(x)} M) / X_j$
      such that $\hR_j (\BB_j) \subset \BB_j$ and
      $\| \hR_j - I \| < \eps$.
\end{itemize}
Consider the linear maps $R_j: T_{f^j(x)} M \to T_{f^j(x)} M$ such
that $R_j$ restricted to $X_j$ is the identity, $R_j(X_j^\perp) =
X_j^\perp$ and $R_j / X_j = \hR_j$. Define
$$
L_j = Df_{f^j(x)} R_j:T_{f^j(x)}M\to T_{f^{j+1}(x)}M
\quad\text{for $0 \leq j \leq n-1$.}
$$
Then $\{L_0,\ldots ,L_{n-1}\}$ is an $(\eps_0,\kappa)$-realizable
sequence of length $n$ at~$x$.
\end{lemma}

We shall prove lemma~\ref{l.nested} in
section~\ref{sss.nestedproof}. Notice that lemma~\ref{l.length 1}
is contained in lemma~\ref{l.nested}: take $n=1$ and use also part
3 of lemma~\ref{l.basicproperties}. Actually, lemma~\ref{l.length
1} also follows from the forthcoming lemma~\ref{l.rotcil}.

\subsubsection{Cylinders and rotations}

We call a \emph{cylinder} any affine image $\CC$ in $\R^d$ of a
product $B^{d-i}\times B^{i}$, where $B^j$ denotes a ball in
$\R^j$. If $\psi$ is the affine map,  the \emph{axis}
$\AA=\psi(B^{d-i}\times\{0\})$ and the \emph{base} $\BB=\psi(\{0\}
\times B^{i})$ are ellipsoids. We also write $\CC=\AA\oplus\BB$.
The cylinder is called \emph{right} if $\AA$ and $\BB$ are
perpendicular. The case we are most interested in is when $i=2$.

The present section contains three preliminary lemmas that we use
in the proof of lemma~\ref{l.nested}. The first one explains how
to rotate a right cylinder, while keeping the complement fixed.
The assumption $a> \tau b$ means that the cylinder $\CC$ is thin
enough, and it is necessary for the $C^1$ estimate in part (ii) of
the conclusion.

\begin{lemma} \label{l.rotcil}
Given $\eps_0>0$ and $0<\sigma<1$, there is $\eps>0$ with the
following properties: Suppose we are given a splitting $\R^d = X
\oplus Y$ with $X \perp Y$ and $\dim Y = 2$, a right cylinder
$\AA\oplus\BB$ centered at the origin with $\AA \subset X$ and
$\BB \subset Y$, and a linear map $\hR:Y\to Y$ such that
$\hR(\BB)=\BB$ and $\|\hR-I\|<\eps$. Then there exists $\tau>1$
such that the following holds:

Let $R:\R^d\to\R^d$ be the linear map defined by $R(u+v)=u + \hR
v$, for $u \in X$, $v\in Y$. For $a$, $b>0$ consider the cylinder
$\CC=a\AA\oplus b\BB$. If $a > \tau b$ and $\diam \CC < \eps_0$
then there is a $C^1$ volume preserving diffeomorphism $h:\R^d \to
\R^d$ satisfying
\begin{itemize}
\item[(i)] $h(z)=z$ for every $z\notin \CC$ and
           $h(z)=R(z)$ for every $z\in \sigma\CC$;
\item[(ii)] $\| h(z) - z \| < \eps_0$ and\ \
             $\| Dh_z - I \| < \eps_0$ for all $z\in\R^d$.
\end{itemize}
\end{lemma}

\begin{proof}
We choose $\eps>0$ small enough so that
\begin{equation} \label{e.1}
\frac{18 \eps}{1-\sigma} < \eps_0.
\end{equation}
Let $\AA$, $\BB$, $X$, $Y$, $\hR$, $R$ be as in the statement of
lemma. Let $\{e_1, \dots, e_d\}$ be an orthonormal basis of $\R^d$
such that $e_1, e_2 \in Y$ are in the directions of the axes of
the ellipse $\BB$ and $e_j\in X$ for $j=3, \dots, d$. We shall
identify vectors $v = x e_1 + y e_2 \in Y$ with the coordinates
$(x,y)$. Then there are constants $\lambda\ge 1$ and $\rho >0$
such that $\BB=\{(x,y);\; \lambda^{-2}x^2 + \lambda^2 y^2 \le
\rho^2\}$. Relative to the basis $\{e_1, e_2\}$, let
$$
H_\lambda=
\begin{pmatrix}
\lambda & 0 \\
0 & \lambda^{-1}
\end{pmatrix}
\quad\text{and}\quad
R_{\alpha }=
\begin{pmatrix}
\cos \alpha & -\sin \alpha \\
\sin \alpha & \cos \alpha
\end{pmatrix}.
$$
The assumption $\hR(\BB)=\BB$ implies that $\hR = H_\lambda
R_\alpha H_\lambda^{-1}$ for some $\alpha$. Besides, the condition
$\|\hR-I\|<\eps$ implies
\begin{equation} \label{e.2}
\lambda^2 \lvert \sin\alpha \rvert \leq \|(\hR-I)(0,1)\| < \eps.
\end{equation}

Let $\varphi:\R \to \R$ be a $C^\infty$ function such that
$\varphi(t)=1$ for $t\leq \sigma$, $\varphi(t)=0$ for $t\geq 1$,
and $0 \leq -\varphi'(t)\leq 2/(1-\sigma)$ for all $t$. Define
smooth maps $\psi:Y\to\R$ and $\tilde{g}_t:Y\to Y$ by
$$\psi(x,y)=\alpha \varphi(\sqrt{x^2+y^2})
 \quad\text{and}\quad
 \tilde{g}_t(x,y) = R_{\varphi(t) \psi(x,y)}(x,y).
$$
On the one hand, $\tilde{g}_t(x,y)=(x,y)$ if either $t\ge 1$ or
$x^2+y^2\ge 1$. On the other hand,
$\tilde{g}_t(x,y)=R_\alpha(x,y)$ if $t\le\sigma$ and
$x^2+y^2\ge\sigma^2$. We are going to check that the derivative of
$\tilde{g}_t$ is close to the identity if $\eps$ is close to zero;
note that $|\sin\alpha|$ is also close to zero, by \eqref{e.2}. We
have
\begin{align*}
D(\tilde{g}_t)_{(x,y)} & =
\begin{pmatrix}
\cos (t \psi) & -\sin (t \psi) \\
\sin (t \psi) &  \cos (t \psi)
\end{pmatrix}
+
\begin{pmatrix}
-x\sin (t \psi) - y\cos (t \psi) \\
 x\cos (t \psi) - y\sin (t \psi)
\end{pmatrix}
\cdot
\begin{pmatrix} t\partial_x \psi & t\partial_y\psi
\end{pmatrix} \\
& = R_{t \psi(x,y)} + t \left[ R_{\pi/2+ t\psi(x,y)}(x,y) \right]
\cdot D\psi_{(x,y)}
\end{align*}
Consider $0 \leq t \leq 1$ and $x^2+y^2 \leq 1$. Then
\begin{align*}
\| D(\tilde{g}_t)_{(x,y)} -I \|
&= \|R_{t \psi(x,y)} -I \| + \| R_{\pi/2+ t\psi(x,y)}(x,y)\| \cdot \|D\psi_{(x,y)} \|   \\
&\leq \big| \sin \big( t\psi(x,y) \big)\big| +
      \big\|\big(2\alpha x \varphi'(x^2+y^2) \ ,\
                 2\alpha y \varphi'(x^2+y^2)\big)\big\|
\end{align*}
Taking $\eps$ small enough, we may suppose that $\alpha\le 2
|\sin\alpha|$. In view of the choice of $\varphi$ and $\psi$, this
implies
\begin{equation}
\label{e.est1}
\| D(\tilde{g}_t)_{(x,y)} - I \|
\leq \lvert \sin \alpha \rvert + 4 \lvert \alpha \rvert/(1-\sigma)
\leq 9 \lvert \sin \alpha \rvert/(1-\sigma).
\end{equation}
We also need to estimate the derivative with respect to $t$:
\begin{equation}
\label{e.est2}
\left\|\partial_t \tilde{g}(x,y)\right\| \leq
\left\|\varphi'(t) \psi(x,y) R_{\pi/2+ t\psi(x,y)}(x,y) \right\|
\leq 4 \lvert \sin \alpha \rvert/(1-\sigma).
\end{equation}
Now define $g_t:Y \to Y$ by $g_t = H_\lambda \circ \tilde{g}_t
\circ H_\lambda^{-1}$. Each $g_t$ is an area preserving
diffeomorphism equal to the identity outside $\BB$. Thus
\begin{equation}
\label{e.2b}
\|g_t(x,y) - (x,y) \| <\diam \BB,
\end{equation}
for every $(x,y)\in\BB$. Moreover, $g_t = \hR = H_\lambda R_\alpha
H_\lambda^{-1}$ on $\sigma\BB$ for all $t\le \sigma$.
By~\eqref{e.est1},
$$
\|D(g_t)_{(x,y)} -I \|
= \big\|
H_\lambda \big( D(\tilde{g}_t)_{(\lambda^{-1} x,\lambda y)}-I \big) H_\lambda^{-1}
\big\|
\leq \lambda^2 \left( \frac{9 \lvert \sin\alpha \rvert}{1-\sigma} \right),
$$
and, applying~\eqref{e.2} and~\eqref{e.1}, we deduce that
\begin{equation}
\label{e.3} \|D(g_t)_{(x,y)} -I \| < \frac{9 \eps}{1-\sigma} <
\frac{\eps_0}{2}
\end{equation}
for all $(x,y)\in\BB$. Similarly, by \eqref{e.est2},
\begin{equation}
\label{e.4} \|\partial_t g_t(x,y)\|
 \leq \lambda^2 \|\partial_t \tilde{g}_t (\lambda^{-1} x,\lambda y)\|
 \leq \lambda^2 \left( \frac{4 \lvert \sin\alpha \rvert}{1-\sigma} \right)
 < \frac{\eps_0}{2}\,.
\end{equation}

Now let $Q: X \to \R$ be a quadratic form such that $\AA=\{u\in X
; \; Q(u) \leq 1\}$, and let $q:\R^d\to X$ and $p:\R^d\to Y$ be
the orthogonal projections. Given $a, b >0$, define $h:\R^d \to
\R^d$ by
$$
h(z) = z' + b g_{a^{-2}Q(z')}(b^{-1} z''), \quad\text{where $z' =
q(z)$ and $z'' = p(z)$.}
$$
It is clear that $h$ is a volume preserving diffeomorphism. The
subscript $t=a^{-2} Q(z')$ is designed so that $t \leq 1$ if and
only if $z'\in a\AA$. Then $h(z)=z$ if either $z'\notin a\AA$ or
$z''\notin b\BB$. Moreover, $h(z)=z'+ \hR(z'')=R(z)$ if
$z'\in\sigma a \AA$ and $z'' \in \sigma b \BB$. This proves
property (i) in the statement. The hypothesis $\diam\CC<\eps_0$
and \eqref{e.2b} give
\begin{align*}
\| h(z)-z \|
& =  b \|g_{a^{-2}Q(z')} (b^{-1} z'') -  b^{-1} z''\|
\\
& < b \diam \BB \le \diam (a\AA\oplus b\BB) < \eps_0
\end{align*}
which is the first half of (ii). Finally, fix $\tau>1$ such that
$\|DQ_u\| \leq \tau \|u\|$ for all $u\in\R^d$, and assume that $a>
\tau b$. Clearly,
$$
Dh = q + \frac{b}{a^2} (\partial_t g) (DQ) q + (Dg) p.
$$
Using~\eqref{e.3}, \eqref{e.4}, and the fact that $\|q \| = \|p \|
= 1$ (these are orthogonal projections),
\begin{align*}
\|Dh-I\|
& \leq \|\frac{b}{a^2} (\partial_t g)(DQ) q \| + \|(Dg -I) p \| \\
& \leq \frac{b}{a^2} \, \| \partial_t g\| \, \tau a \, \|q \| + \|
Dg -I \| \, \| p \| <    \eps_0 \,. &
\end{align*}
This completes the proof of property (ii) and the lemma.
\end{proof}


The second of our auxiliary lemmas says that the image of a small
cylinder by a $C^1$ diffeomorphism $h$ contains the image by $Dh$
of a slightly shrunk cylinder. Denote $\CC(y,\rho) = \rho  \CC +
y$, for each $y\in\R^d$ and $\rho >0$.

\begin{lemma} \label{l.nlinear}
Let $h:\R^d \to \R^d$ be a $C^1$ diffeomorphism with $h(0)=0$,
$\CC\subset \R^d$ be a cylinder centered at $0$, and $0 < \lambda
< 1$. Then there exists $r>0$ such that for any
$\CC(y,\rho)\subset B_r(0)$,
$$
h(\CC (y,\rho)) \supset Dh_0(\CC(0,\lambda\rho)) + h(y) .
$$
\end{lemma}

\begin{proof}
Fix a norm $\| \cdot \|_0$ in $\R^d$ for which $\CC=\{ z
\in\R^d;\; \|z\|_0 < 1 \}$. Such a norm exists because $\CC$ is
convex and $\CC = - \CC$.
Let $H=Dh_0$ and $g:\R^d \to \R^d$ be such that $h=H \circ
g$. Since $g$ is $C^1$ and $Dg_0 = I$, we have
$$
g(z)-g(y)=z-y+\xi(z,y)
\quad\text{with}\quad
\lim_{(z,y)\to(0,0)}\frac{\xi(z,y)}{\|z-y\|_0}=0.
$$
Choose $r>0$ such that
$\|z\| ,\|y\| \leq r \Rightarrow \| \xi (z,y)\|_0 < (1-\lambda) \| z-y \|_0$
(where $\|\cdot \|$ denotes the Euclidean norm in $\R^d$).
Now suppose $\CC(y,\rho) \subset B_r(0)$, and let $z\in \partial \CC (y,\rho )$.
Then $\|z - y \|_0 = \rho$ and
$$
\| g(z)-g(y) \|_0 \geq \| z-y \|_0 - \|\xi(z,y)\|_0 > \lambda \rho.
$$
This proves that the sets
$g(\partial \CC(y,\rho))-g(y)$ and $\lambda \CC$ are disjoint.
Applying the linear map $H$, we find that
$h(\partial \CC(y,\rho ))-h(y)$ and $\lambda H\CC$
are disjoint.
From topological arguments,
$h(\CC(y,\rho))-h(y) \supset \lambda H\CC$.
\end{proof}

\begin{figure}[ht]
\begin{center}
\psfrag{A}{$a\AA$}
\psfrag{B}{$b\BB$}
\psfrag{A1}{$\lambda a\AA_1$}
\psfrag{B1}{$b\BB_1$}
\psfrag{LB}{$bL(\BB)$}
\psfrag{L}{$L$}
\includegraphics[width=4in]{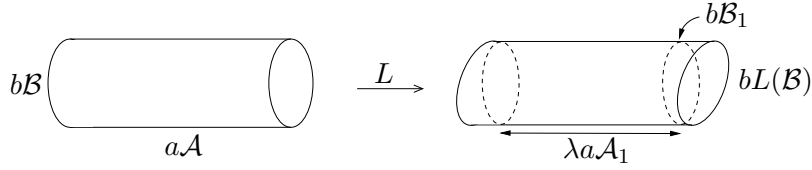}
\caption{\label{f.truncate} Truncating a thin cylinder to make it
right}
\end{center}
\end{figure}

The third lemma says that a linear image of a sufficiently thin
cylinder contains some right cylinder with almost the same volume.
The idea is contained in figure~\ref{f.truncate}. The proof of the
lemma is left to the reader.

\begin{lemma} \label{l.trunca}
Let $\AA\oplus\BB$ be a cylinder centered at the origin,
$L:\R^d\to\R^d$ be a linear isomorphism, $\AA_1 = L(\AA)$ and
$\BB_1 = p(L(\BB))$, where $p$ is the orthogonal projection onto
the orthogonal complement of $\AA_1$. Then, given any
$0<\lambda<1$, there exists $\tau > 1$ such that if $a > \tau b$,
$$
L\big( a\AA \oplus b\BB \big)
\supset \lambda a \AA_1  \oplus b \BB_1.
$$
\end{lemma}

\subsubsection{Proof of the nested rotations lemma~\ref{l.nested}}\label{sss.nestedproof}

\begin{proof}
Let $f$, $\eps_0$, and $\kappa$ be given. Define $\sigma =
(1-\kappa)^{1/2d}$ and then take $\eps>0$ as given by
lemma~\ref{l.rotcil}. Now let $x$, $n$, $X_j$, $\BB_j$, $\hR_j$,
$R_j$, $L_j$ be as in the statement. We want to prove that $\{L_0,
\ldots, L_n\}$ is an $(\eps_0,\kappa)$-realizable sequence of
length at $x$, cf. definition~\ref{d.sr}.

In short terms, we use lemma~\ref{l.rotcil} to construct the
realization $g$ at each iterate. The subset $U\minus K$, where
we have no control on the approximation, has two sources:
lemma~\ref{l.rotcil} gives $h=R$ only on a slightly smaller
cylinder $\sigma\CC$; and we need to straighten out
(lemma~\ref{l.nlinear}) and to ``rightify'' (lemma~\ref{l.trunca})
our cylinders at each stage. These effects are made small by
considering cylinders that are small and very thin. That is how we
get $U\minus K$ with relative volume less than $\kappa$,
independently of~$n$.

For clearness we split the proof into three main steps:

\paragraph{Step 1:} Fix any $\gamma>0$. We explain how to find
$r>0$ as in definition~\ref{d.sr}.

\smallskip

We consider local charts $\varphi_j : V_j \to \R^d$ with
$\varphi_j=\varphi_{i(f^j x)}$ and $V_j = V_{i(f^j x)}$, as
introduced in section~\ref{ss.basic}.  Let $r'>0$ be small enough
so that
\begin{itemize}
\item $f^j( \overline{B}_{r'}(x)) \subset V_j^\ast$ for every $j=0,1\ldots,n$;

\item the sets $f^j(\overline{B}_{r'}(x))$ are two-by-two disjoint;

\item $\| Df_z - Df_{f^j(x)} \| \, \| R_j \| <\gamma $
for every $z\in f^j (B_{r'}(x))$ and  $j=0,1\ldots,n$.

\end{itemize}

We use local charts to translate the situation to $\R^d$. Let $f_j
= \varphi_{j+1} \circ f\circ \varphi_j^{-1}$ be the expression of
$f$ in local coordinates near $f^j(x)$ and $f^{j+1}(x)$. To
simplify the notations, we suppose that each $\varphi_j$ has been
composed with a translation to ensure $\varphi_j(f^j(x))=0$ for
all $j$. Up to identification of tangent spaces via the charts
$\varphi_j$ and $\varphi_{j+1}$, we have $L_j = (Df_j)_0 R_j$.

Let $\AA_0 \subset X_0$ be any ellipsoid centered at the origin (a
ball, for example), and let $\AA_j = Df_x^j(\AA_0)$ for $j \geq
1$. We identify $(T_{f^j(x)} M)/X_j$ with $X_j^\perp$, so that we
may consider $\BB_j \subset X_j^\perp$.
In these terms, the assumption $\BB_j = (Df^j_x /X_0)(\BB_0)$ means that
$\BB_j$ is the orthogonal projection of $Df^j_x (\BB_0)$ onto $X_j^\perp$.

Fix $0<\lambda<1$ close enough to $1$ so that $\lambda^{4n(d-1)}
> 1-\kappa$. Let $\tau_j>1$ be associated to the data
$(\AA_j\oplus\BB_j, (Df_j)_0,\lambda)$ by lemma~\ref{l.trunca}: if
$a > \tau_j b$ then
\begin{equation} \label{e.trunca}
(Df_j)_0 (a\AA_j\oplus b\BB_j) \supset \lambda a \AA_{j+1} \oplus
b \BB_{j+1}
\end{equation}
and let $\tau'_j>1$ be associated to the data $(\eps_0, \sigma,
X_j \oplus X_j^\perp, \AA_j \oplus \BB_j, \hR_j)$ by
lemma~\ref{l.rotcil}. Fix $a_0>0$ and $b_0>0$ such that
\begin{equation} \label{e.fininho}
a_0 > b_0 \lambda^{-n} \max\{\tau_j,\tau'_j ; \; 0 \leq j \leq n-1\}.
\end{equation}
For $0 \leq j \leq n$, define $\CC_j = \lambda^{2j} a_0 \AA_j
\oplus \lambda^j b_0 \BB_j$. For $z\in\R^d$ and $\rho>0$, denote
$\CC_j (z, \rho) = \rho \CC_j + z$. Applying lemma~\ref{l.nlinear}
to the data $(f_j, \CC_j, \lambda)$ we get $r_j>0$ such that
\begin{equation} \label{e.def rj}
\CC(z,\rho) \subset B_{r_j} (0) \quad\Rightarrow\quad f_j(\CC_j(z,
\rho)) \supset (Df_j)_0 (\CC_j(0, \lambda \rho)) + f_j(z).
\end{equation}
Now take $r>0$ such that $r<r'$ and, for each $j=1,\ldots, m-1$,
\begin{equation}\label{e.def r}
f_{j-1} \cdots f_0 (B_r(0)) \subset B_{r_j} (0).
\end{equation}

\paragraph{Step 2:} Let $U$ be fixed. We find $g \in
\UU(f,\eps_0)$ and $K \subset U$ as in definition~\ref{d.sr}.

\smallskip

For this we take advantage of lemma~\ref{l.simplifying}: it suffices to
consider open sets of the form $U=\varphi_0^{-1}(\CC_0(y_0,\rho))$,
because the cylinders $\CC_0(y_0,\rho)$ contained in $B_r(0)$
constitute a Vitali covering.

We claim that, for each $j=0,1,\ldots,m-1$, and every $t\in [0, \rho]$,
\begin{equation} \label{e.red}
\CC_j (y_j, t)
\subset f_{j-1} \cdots f_0 (B_r(0))
\end{equation}
and
\begin{equation} \label{e.black}
f_j (\CC_j (y_j,t)) \supset \CC_{j+1} (y_{j+1}, t)
\end{equation}
For $j=0$, relation~\eqref{e.red} means $\CC_0 (y_0,t) \subset
B_r(0)$, which is true by assumption. We proceed by induction.
Assume~\eqref{e.red} holds for some $j\ge 0$. Then,
by~\eqref{e.def r} and~\eqref{e.def rj},
\begin{align*}
f_j(\CC_j(y_j, t))
&\supset (Df_j)_0 (\CC_j(0, \lambda t)) + y_{j+1} \\
& = (Df_0)_0 \left[(\lambda^{2j+1} t a_0 \AA_j) \oplus
(\lambda^{j+1} t b_0 \BB_j)\right] + y_{j+1}.
\end{align*}
Relation~\eqref{e.fininho} implies that $\lambda^{2j+1} t a_0 >
\tau_j (\lambda^{j+1} t b_0)$. So, we may use~\eqref{e.trunca} to
conclude that
$$
f_j(\CC_j(y_j, t)) \supset (\lambda^{2j+2} t a_0  \AA_j) \oplus
(\lambda^{j+1} t b_0 \BB_0) + y_{j+1} = \CC_{j+1} (y_{j+1}, t)
$$
This proves that \eqref{e.black} holds for the same value of $j$.
Moreover, it is clear that if \eqref{e.black} holds for all $0\le
i \le j$ then \eqref{e.red} is true with $j+1$ in the place of
$j$. This completes the proof of \eqref{e.red} and
\eqref{e.black}.

\smallskip

Condition~\eqref{e.fininho} also implies $\lambda^{2j} a_0 >
\tau_j'\, (\lambda^{j} b_0)$. So, we may use lemma~\ref{l.rotcil}
(centered at $y_j$) to find a volume preserving diffeomorphism
$h_j:\R^d \to \R^d$ such that
\begin{enumerate}
\item $h_j(z) = z$ for all $z\notin \CC_j(y_j, \rho)$ and
      $h_j(z) = y_j + R_j(z-y_j)$ for all $z \in \CC_j(y_j, \sigma
      \rho)$ and, consequently,
\begin{equation}\label{eq.consequently}
h_j(\CC_j(y_j, \sigma\rho)) = \CC_j(y_j, \sigma\rho)
 \quad\text{and}\quad
 h_j(\CC_j(y_j, \rho)) = \CC_j(y_j, \rho).
\end{equation}
\item $\| h_j(z)-z\| < \eps_0$ and $\|(Dh_j)_z- I \|<\eps_0$
      for all $z\in\R^d$.
\end{enumerate}

$R_j$ is the linear map $T_{f^j(x)}\to T_{f^{j+1}(x)}$ in the
statement of the theorem or, more precisely, its expression in
local coordinates $\varphi_j$. Let $S_j = \varphi_j^{-1}(\{ z;\;
h(z) \neq z \}) \subset M$. By property 1 above and the
inclusion~\eqref{e.red},
$$
S_j \subset \varphi_j^{-1}(f_{j-1} \cdots f_0 (B_r(0))) = f^j(B_r(x)).
$$
In particular, the sets $S_j$ have pairwise disjoint closures.
This permits us to define a diffeomorphism $g\in \Diff$ by
$$
g= \left\{\begin{array}{ll}
 \varphi_{j+1}^{-1} \circ (f_j \circ h_j) \circ \varphi_j &
 \text{on $S_j$ for each } 0 \le j \le n-1 \\
 f & \text{outside $S_0 \sqcup\cdots \sqcup S_{n-1}$}
 \end{array}\right.
$$
Property 2 above gives that $f^{-1} \circ g \in \UU(\id,\eps_0)$,
and so $g \in \UU(f, \eps_0)$.

\paragraph{Step 3:} Now we define $K\subset U$ and check
the conditions (i)--(iii) in definition~\ref{d.sr}.

\smallskip

By construction, $h_j = \id$ outside $C_j(y_j,\rho)$, and so
$$
\varphi_{j+1}^{-1} \circ (f_j \circ h_j) \circ \varphi_j = f
\quad\text{outside $\varphi^{-j}(C_j(y_j,\rho))$.}
$$
Using \eqref{e.black} and \eqref{eq.consequently}, we have
$\varphi^{-j}(C_j(y_j,\rho))\subset f^j(U)$ for all $0\le j \le
n-1$. Recall that $U = \varphi_0^{-1}(\CC_0(y_0,\rho))$. Hence,
$g=j$ outside the disjoint union $\sqcup_{j=0}^{n-1} f^j(U)$. This
proves condition (i).

\smallskip

Define $K=g^{-n}(\varphi_n^{-1}(\CC_n(y_n,\sigma \rho)))$.
Using~\eqref{e.black} and \eqref{eq.consequently} in the same way
as before, we see that $K\subset U$.
Also, since all the maps $f$, $g$, $h_j$, $\varphi_j$ are volume
preserving, and all the cylinders $\CC_j(y_j,\rho)$,
$\CC_j(y_j,\sigma\rho)$, are right
$$
\frac{\vol K}{\vol U} = \frac{\vol (\sigma\rho\,\lambda^{2n}  a
\AA_n \oplus \sigma\rho\,\lambda^n b \BB_n)}
       {\vol (\rho a \AA_0 \oplus \rho b \BB_0)}
= \frac{ (\lambda^{2n} \sigma)^{d-2} \vol\AA_n
           \, (\lambda^n \sigma)^2 \, \vol \BB_n}
        {\vol \AA_0 \, \vol \BB_0}\,.
$$
Notice also that $\vol\AA_n \, \vol \BB_n = \vol\AA_0 \, \vol \BB_0$,
since the cylinders $Df^n_x(\AA_0 \oplus \BB_0)$ and
$\AA_n \oplus \BB_n$ differ by a sheer.
So, the right hand side is equal to
$\lambda^{2n(d-1)} \sigma^d$. Now, this expression is larger than
$1-\kappa$, because we have chosen $\sigma=(1-\kappa)^{1/2}$ and
$\lambda>(1-\kappa)^{1/4n(d-1)}$. This gives condition (ii).

\smallskip

Finally, let $z\in K$. Recall that $L_j=Df_{f^j(x)} R_j$.
Moreover, $(Dh_j)_{\varphi_jg^j(z)} = R_j$ (we continue to
identify $R_j$ with its expression in the local chart
$\varphi_j$), because
$$
g^j(z) \in g^{-n+j}(\varphi_n^{-1}(\CC_n(y_n, \sigma\rho)))
\subset \varphi_j^{-1}(\CC_j(y_j,\sigma\rho)).
$$
Therefore, writing $z_j=h_j(\varphi_j(g^j(z)))$ for simplicity,
$$
\big\| Dg_{g^j(z)}-L_j \big\| = \big\| D(f_j)_{z_j} R_j - D(f_j)_0
R_j\big\| \leq \big\|D(f_j)_{z_j}-D(f_j)_0\big\| \big\| R_j \big\|
< \gamma.
$$
The last inequality follows from our choice of $r'$. This gives
condition (iii) in definition~\ref{d.sr}. The proof of
lemma~\ref{l.nested} is complete.
\end{proof}

\begin{remark}
This last step explains why it is technically more convenient to
require $\|Dg_{g^j(z)}-L_j\|<\gamma$, rather than
$Dg_{g^j(z)}=L_j$, when defining realizable sequence.
\end{remark}

\subsection{Proof of the directions interchange proposition~\ref{p.geom}}

\begin{proof}

First, we define some auxiliary constants. Fix $0 < \kappa' <
\meio \kappa$. Let $\eps_1>0$, depending on $f$, $\eps_0$ and
$\kappa'$, be given by lemma~\ref{l.length 1}. Let $\eps_2>0$,
depending on $f$, $\eps_0$ and $\kappa$, be given by
lemma~\ref{l.nested}. Take $\eps = \min \{\eps_1,\eps_2\}$. Fix
$\alpha>0$ such that $\sqrt{2} \, \sin \alpha < \eps$. Take
\begin{equation}\label{eq.def K}
K \ge (\sin \alpha)^{-2} \quad\text{and}\quad K \ge
\max\big\{{\|Df_x\|}/{\mm(Df_x)};\; x\in M\big\}.
\end{equation}
Let $\beta>0$ be such that
\begin{equation} \label{e.def beta}
{8 \sqrt{2} \; K \sin \beta} < \eps \, {\sin^6 \alpha}.
\end{equation}
Finally, assume $m\in \N$ satisfies $m \geq 2\pi / \beta$.

\smallskip

Let $y\in M$ be a non-periodic point and $T_y M = E \oplus F$ be a
splitting as in the hypothesis:
\begin{equation}\label{e.not dominated}
\frac{\|Df^m_y|_F \|}{\mm(Df^m_y|_E)} \geq \frac 12\,.
\end{equation}
We write $E_j = Df^j_y (E)$ and $F_j = Df^j_y (F)$ for
$j=0,1,\ldots,m$. The proof is divided in three cases.
Lemma~\ref{l.length 1} suffices for the first two, in the third
step we use the full strength of lemma~\ref{l.nested}.

\paragraph{First case:}
Suppose there exists $\ell \in \{0,1,\ldots,m\}$ such that
\begin{equation} \label{e.I}
\ang(E_\ell, F_\ell) < \alpha.
\end{equation}

Fix $\ell$ as above. Take unit vectors $\xi \in E_\ell$ and $\eta
\in F_\ell$ such that $\ang(\xi,\eta) < \alpha$. Let $Y = \R \xi
\oplus \R \eta$ and $X = Y^\perp$. Let $\hR : Y \to Y$ be a
rotation such that $\hR (\xi) = \eta$. Then  $\| \hR - I \|
=\sqrt{2} \; \sin \ang(\xi,\eta)< \eps$. Let $R : T_{f^\ell(y)}M
\to T_{f^\ell(y)}M$ be such that $R$ preserves both $X$ and $Y$,
$R|_X =I$ and $R|_Y = \hR$.

Consider first $\ell < m$. By lemma~\ref{l.length 1}, the length
$1$ sequence $\{Df_{f^\ell(y)} R\}$ is $(\kappa',
\eps_0)$-realizable at $f^\ell(y)$. Using part 2 of
lemma~\ref{l.basicproperties} we conclude that
$$
\{ L_0, \ldots L_{m-1} \} = \{ Df_y, \ldots, Df_{f^{\ell-1}(y)},
Df_{f^\ell(y)} R, Df_{f^{\ell+1}(y)}, \ldots, Df_{f^{m-1}(y)} \}
$$
is a $(\kappa, \eps_0)$-realizable sequence of length $m$ at $y$.
The case $\ell = m$  is similar. By lemma~\ref{l.length 1}, the
length $1$ sequence $\{ R \, Df_{f^{m-1}(y)} \}$ is $(\kappa',
\eps_0)$-realizable at $f^{m-1}(y)$. Then, by part 2 of
lemma~\ref{l.basicproperties},
$$
\{ L_0, \ldots L_{m-1} \} = \{ Df_y, \ldots, Df_{f^{m-2}(y)}, R \,
Df_{f^{m-1}(y)} \}.
$$
is a $(\kappa, \eps_0)$-realizable sequence of length $m$ at $y$.
In either case, $L_{m-1} \cdots L_0$ sends the vector $v =
Df^{-\ell} (\xi) \in E_0$ to a vector $w$ collinear to
$Df^{m-\ell} (\eta) \in F_m$.

\paragraph{Second case:}
Assume there exist $k, \ell \in \{0,\ldots,m\}$, with
$k<\ell$, such that
\begin{equation} \label{e.II}
\frac{\|Df^{\ell-k}_{f^k(y)} |_{F_k} \|}{\mm (Df^{\ell-k}_{f^k(y)}
|_{E_k})} > K.
\end{equation}

Fix $k$ and $\ell$ as above.
Let $\xi \in E_k$, $\eta \in F_k$ be unit vectors
such that
$$
\| Df^{\ell-k}(\xi)\| = \mm (Df^{\ell-k} |_{E_k}) \quad \text{and}
\quad \| Df^{\ell-k}(\eta) \| = \| Df^{\ell-k} |_{F_k} \|
$$
($Df^{\ell-k}$ is always meant at the point $f^k(y)$). Define also
unit vectors
$$
\xi ' = \frac{Df^{\ell-k}(\xi)}{\| Df^{\ell-k}(\xi) \|} \in E_\ell
\quad\text{and}\quad \eta' = \frac{Df^{\ell-k}(\eta)}{\|
Df^{\ell-k}(\eta) \|} \in F_\ell \,.
$$
Let $\xi_1 = \xi + (\sin \alpha) \eta$. Then $\theta = \ang(\xi,
\xi_1) \leq \alpha$, simply because $\|\xi\|=\|\eta\|=1$. In
particular, if $\hR : \R \xi \oplus \R \eta \to \R \xi \oplus \R
\eta$ is a rotation of angle $\pm \theta$, sending $\R \xi$ to $\R
\xi_1$ then
$$
\|\hR - I  \| = \sqrt{2}\; \sin \theta < \eps.
$$
Let $Y = \R \xi \oplus \R \eta$ and $X=Y^\perp$. Let $R :
T_{f^k(y)}M \to T_{f^k(y)}M$ be such that $R$ preserves both $X$
and $Y$, with $R|_X =I$ and $R|_Y = \hR$. By lemma~\ref{l.length
1}, the length 1 sequence $\{Df_{f^k(y)} R\}$ is $(\kappa',
\eps_0)$-realizable at $f^k(y)$. Let $\eta'_1 = s \xi' + \eta'$,
where
$$
s
 = \frac{1}{\sin \alpha}\,
   \frac{\|Df^{\ell-k}(\xi)\|}{\| Df^{\ell-k}(\eta)\|}
 = \frac{1}{\sin \alpha}\,
   \frac{\mm (Df^{\ell-k}|_{E_k})}{\|Df^{\ell-k}|_{F_k} \|}.
$$
Then the vectors $Df^{\ell-k}\xi_1$ and $\eta_1$ are collinear.
Besides, $s < 1/(K\sin \alpha) < \sin \alpha$, because of
\eqref{eq.def K} and \eqref{e.II}. Hence, $\theta' = \ang(\eta'_1,
\eta) < \alpha$. Then, as before, there exists $R' : T_{f^\ell(y)}M
\to T_{f^\ell(y)}M$ such that $R'(\R \eta'_1) = \R \eta$ and $\{R'
\, Df_{f^{\ell-1}(y)} \}$ is a $(\kappa', \eps_0)$-realizable
sequence of length $1$ at $f^{\ell-1}(y)$.

Notice that~\eqref{eq.def K} and~\eqref{e.II} imply $\ell - 1
> k$. Then we may define a sequence $\{L_0$, \ldots, $L_{m-1} \}$
of linear maps as follows:
$$
L_j = \left\{\begin{array}{ll} Df_{f^k(y)} \, R & \text{for } j=k \\
                               R' \, Df_{f^{\ell-1} (y)} & \text{for } j=\ell-1 \\
                               Df_{f^j(y)} & \text{for all other $j$.}
      \end{array}\right.
$$
By parts 1 and 2 of lemma~\ref{l.basicproperties}, this is a
$(\kappa,\eps_0)$-realizable sequence of length $m$ at $y$. By
construction, $L_{m-1} \cdots L_0$ sends $v=Df^{-k} (\xi) \in E_0$
to a vector $w$ collinear to $Df^{m-\ell} (\eta') \in F_m$.

\paragraph{Third case:}
We suppose that we are not in the previous
cases, that is, we assume
\begin{equation} \label{e.not I}
\text{for every $j \in \{0,1,\ldots,m\}$,}\quad
\ang(E_j, F_j) \geq \alpha.
\end{equation}
and
\begin{equation} \label{e.not II}
\text{for every $i,j \in \{0,\ldots,m\}$ with $i<j$,} \quad
\frac{\|Df^{j-i}_{f^i(y)} |_{F_i} \|}{\mm (Df^{j-i}_{f^i(y)}
|_{E_i})} \leq K.
\end{equation}
We now use the assumption~\eqref{e.not dominated}, and the choice
of $m$ in \eqref{e.def beta}. Take unit vectors $\xi  \in E_0$ and
$\eta \in F_0$ such that
 $\| Df^m \xi \| = \mm (Df^m |_{E_0}) $ and
 $\| Df^m \eta\| = \| Df^m |_{F_0} \|$ ($Df^m$ is always computed at $y$). 
Let also $\eta' = Df^m(\eta)/ \| Df^m(\eta) \| \in F_m$.

\smallskip

Define $G_0=E_0 \cap \xi^\perp$ and $G_j = Df^j_y (G_0)\subset
E_j$ for $0<j\le m$. Dually, define $H_m= F_m \cap {\eta'}^\perp$
and $H_j = Df^{j-m}(H_m) \subset F_j$ for $0\le j <m$. In
addition, consider unit vectors $v_j \in E_j \cap G_j^\perp$ and
$w_j \in F_j \cap H_j^\perp$ for $0\le j \le m$. These vectors are
uniquely defined up to a choice of sign, and $v_0 = \pm \xi$ and
$w_m = \pm \eta'$. See Figure~\ref{f.spaces}. For $j=0,\ldots,m$,
define
$$
X_j = G_j \oplus H_j
 \quad\text{and}\quad
 Y_j = \R v_j \oplus \R w_j\,.
$$
The spaces $X_j$ are invariant: $Df_{f^j(y)}(X_j) = X_{j+1}$ (the
$Y_j$ are not). We shall prove, using \eqref{e.not I} and
\eqref{e.not II}, that the maps $Df^j_y / X_0 : T_yM / X_0 \to
T_{f^j(y)}M / X_j$ do not distort angles too much:

\begin{lemma}\label{l.eccentricity}
For every $j=0,1,\ldots,m$,
$$
\frac{\| Df^j_y / X_0 \|}{\mm( Df^j_y / X_0 )} \leq
\frac{8K}{\sin^6 \alpha}.
$$
\end{lemma}

\begin{figure}[h]
\begin{center}
\psfrag{E0}{$E_0$}\psfrag{Ej}{$E_j$} \psfrag{Em}{$E_m$}
\psfrag{F0}{$F_0$}\psfrag{Fj}{$F_j$} \psfrag{Fm}{$F_m$}
\psfrag{G0}{$G_0$}\psfrag{Gj}{$G_j$} \psfrag{Gm}{$G_m$}
\psfrag{H0}{$H_0$}\psfrag{Hj}{$H_j$} \psfrag{Hm}{$H_m$}
\psfrag{v0}{$\xi=v_0$}\psfrag{vj}{$v_j$} \psfrag{vm}{$v_m$}
\psfrag{w0}{$w_0$}\psfrag{wj}{$w_j$} \psfrag{wm}{$\eta'=w_m$}
\psfrag{f+}{$Df^j$}\psfrag{f-}{$Df^{j-m}$}
\includegraphics[width=4in]{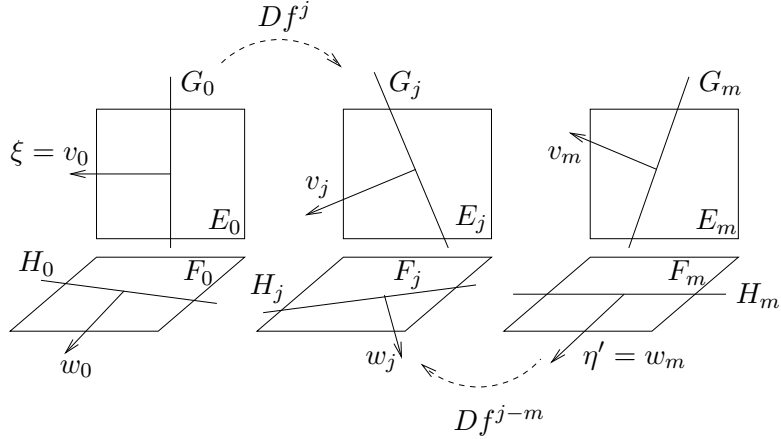}
\caption{\label{f.spaces} Setup for application of the nested
rotations lemma}
\end{center}
\end{figure}

Let us postpone the proof of this fact for a while, and proceed
preparing the application lemma~\ref{l.nested}. Let $\BB_0 \subset
(T_yM)/X_0$ be a ball and $\BB_j = (Df_y^j / X_0)(\BB_0)$ for
$0<j\le m$. Since $m\beta \geq 2\pi$, it is possible to choose
numbers $\theta_0, \ldots, \theta_{m-1}$ such that $|\theta_j|
\leq \beta$ for all $j$ and
\begin{equation}\label{eq.sumtheta}
\sum_{j=0}^{m-1} \theta_j = \ang(v_0 + X_0, w_0 + X_0).
\end{equation}
Let $P_j : (T_yM)/X_0 \to (T_yM)/X_0$ be the rotation of angle
$\theta_j$. Define linear maps $\hR_j : (T_{f^j(y)}M)/X_j \to
(T_{f^j(y)}M)/X_j$ by
$$
\hR_j = \big( Df^j_y / X_0 \big) \, P_j \, \big( Df^j_y / X_0
\big)^{-1}.
$$
Since $P_j$ preserves the ball $\BB_0$, we have $\hR_j
(\BB_j)=\BB_j$ for all $j$. Moreover,
$$
\| \hR_j -I \| \le \frac{\|Df^j_y/X_0\|}{\mm(Df^j_y/X_0)}\,\|P_j-I\|
\le \frac{8K}{\sin^6\alpha} \, \sqrt{2}\sin\beta < \eps,
$$
by lemma~\ref{l.eccentricity}, the relation  $\|P_j - I \| \leq
\sqrt{2} \; \sin\beta$, and our choice~\eqref{e.def beta} of
$\beta$.

\smallskip

Applying lemma~\ref{l.nested} to these data $(\eps_0, \kappa, x=y,
n=m, X_j, \BB_j, \hR_j)$ we obtain an $(\eps_0,
\kappa)$-realizable sequence $\{L_0, \ldots, L_{m-1}\}$ at the
point $y$, with $ L_j |_{X_j} = Df_{f^j(y)} |_{X_j}$ and
$$
L_j / X_j = (Df_{f^j(y)} / X_j) \, \hR_j
 = (Df_y^{j+1} / X_j) \, P_j \, (Df_y^j / X_0)^{-1}.
$$
Let $\LL = L_{m-1} \cdots L_0$. Then $\LL / X_0 = (Df_y^m / X_0)
P_{m-1} \cdots P_0$. In particular, by \eqref{eq.sumtheta},
$$
\LL(v_0 + X_0) = (Df_y^m/X_0)(w_0 + X_0)= Df^m_y(w_0)+X_m\,.
$$
Recall that $X_m = G_m\oplus H_m$ by definition. Then we may write
$$\LL (v_0) = Df_y^m(w_0) + u_m + u_m'$$ with $u_m \in G_m$ and
$u_m' \in H_m$. Let $u_0=(Df_y^m)^{-1}(u_m)\in G_0 \subset X_0\cap
E_0$. Since $\LL$ equals $Df^m_y$ on $X_0$, we have
$\LL(u_0)=u_m$. This means that the vector $v = v_0 - u_0 \in E_0$
is sent by $\LL$ to the vector $Df_y^m(w_0) + u_m' \in F_m$. This
finishes the third and last case of proposition~\ref{p.geom}.

\smallskip

Now we are left to give the

\begin{proof}[Proof of lemma~\ref{l.eccentricity}]
Recall that $X_j=G_j\oplus H_j$\,, $G_j\subset E_j$ and
$H_j\subset F_j$\,, $v_j \in E_j$ and $w_j\in F_j$\,, and
$v_j\perp G_j$ and $w_j\perp H_j$. Hence, using \eqref{e.not I},
\begin{alignat*}{1}
\ang (X_j, v_j) & = \ang (H_j, v_j) \ge \ang(F_j,E_j) \ge \alpha
 \quad\text{and} \\
 \ang (X_j \oplus \R v_j, w_j) & = \ang (\R v_j \oplus G_j, w_j) \ge \ang(E_j,F_j) \geq \alpha
\end{alignat*}
Using lemma~\ref{l.abc} with $A=X_j$, $B=\R v_j$, $C=\R w_j$, we
deduce the following lower bound for the angle between the spaces
$X_j$ and $Y_j=\R v_j\oplus\R w_j$ :
$$
\sin \ang (X_j,Y_j) \geq \sin\ang (X_j, v_j) \sin\ang(\R v_j
\oplus X_j, w_j)  \geq \sin^2 \alpha.
$$
Let $\pi_j : Y_j \to (T_{f^j(y)} M)/X_j$ be the canonical map
$\pi_j(w) = w + X_j$. Then $\pi_j$ is an isomorphism, $\| \pi_j \|
= 1$ and
\begin{equation}\label{eq.angXY}
\| \pi_j^{-1} \| = 1/\sin\ang(Y_j, X_j) \leq 1/\sin^2 \alpha.
\end{equation}
(the quotient space has the norm that makes $X_j^\perp\ni w
\mapsto w + X_j$ an isometry).

\smallskip

Now let $p_j : T_{f_j(y)}M \to Y_j$ be the projection onto $Y_j$
associated to the splitting $T_{f_j(y)}M = X_j \oplus Y_j$. Let
$\DD_j : Y_j \to Y_{j+1}$ be given by $\DD_j = p_{j+1} \circ
(Df_{f^j(y)}|_{Y_j})$. Define
$$
\DD^{(j)}:Y_0 \to Y_j \quad\text{by}\quad \DD^{(j)} = \DD_{j-1}
\circ\cdots\circ \DD_{0} = p_j \circ (Df^j_y |_{Y_0}).
$$
We claim that the following inequalities hold:
\begin{equation}\label{e.supersonic}
\frac{1}{2K} <
\frac{\| \DD^{(j)} (w_0) \|}{\| \DD^{(j)} (v_0) \|}
\leq K
\quad \text{for every $j$ with $0\leq j \leq m$.}
\end{equation}
To prove this, consider the matrix of $\DD_j$ relative to
bases $\{ v_j, w_j \}$ and $\{ v_{j+1}, w_{j+1} \}$:
$$
\DD_j =
\begin{pmatrix}
a_j & 0 \\ 0 & b_j
\end{pmatrix}.
$$
Then $\| \DD^{(j)} (v_0) \| = |a_{j-1} \cdots a_0 |$ and  $\|
\DD^{(j)} (w_0) \| = |b_{j-1} \cdots b_0 |$\,, since $v_j$ and
$w_j$ are unit vectors. Moreover, for $0 \leq i < j \leq m$ we
have
\begin{alignat*}{2}
|a_{j-1} \cdots a_i | &= \| p_{j} \circ Df^{j-i}_{f^i(y)} (v_i) \|
&&= \|p_i \circ Df^{-(j-i)}_{f^j(y)} (v_{j}) \|^{-1},
\\
|b_{j-1} \cdots b_i | &= \| p_{j} \circ Df^{j-i}_{f^i(y)} (w_i) \|
&&= \|p_i \circ Df^{-(j-i)}_{f^j(y)} (w_{j}) \|^{-1}.
\end{alignat*}
Recall that $v_s\in E_s$ and $w_s\in F_s$ for all $s$. When
restricted to $E_s$ (or $F_s$), the map $p_s$ is the orthogonal
projection to the direction of $v_s$ (or $w_s$).
In particular, $\| p_i |_{E_i}\| = \| p_j |_{F_j} \| = 1$ and so
$$
|a_{j-1} \cdots a_i |
 \geq \| Df^{-(j-i)}_{f^j(y)} |_{E_{j}} \|^{-1}
 = \mm(Df^{j-i}_{f^i(y)} |_{E_i})
 \quad\text{and}\quad
 |b_{j-1} \cdots b_i | \leq \| Df^{j-i}_{f^i y} |_{F_i}  \|\,
$$
Using~\eqref{e.not II}, we obtain that
\begin{equation}\label{eq.gettingthere}
\frac{|b_{j-1} \cdots b_i |}{|a_{j-1} \cdots a_i |} \leq K
 \quad\text{for all $0\le i \le k \le m$}.
\end{equation}
Taking $i=0$ gives the upper inequality in \eqref{e.supersonic}.
For the same reasons, and the definitions of $v_0=\xi$ and
$w_m=\eta'=Df^m_y(\eta)/\|Df^m_y(\eta)\|$, we also have
\begin{alignat*}{2}
|a_{m-1} \cdots a_0 | & \leq \| Df^m_y (v_0) \|
 = \| Df^m_y (\xi) \| = \mm(Df^m |_{E_0}),
\\
|b_{m-1} \cdots b_0 | & \geq \| Df^{-m}_{f^m(y)} (w_m) \|^{-1} = \|
Df^m_y (\eta) \| = \|Df^m |_{F_0}\|\,.
\end{alignat*}
Now~\eqref{e.not dominated} translates into
$$
\frac{|b_{m-1} \cdots b_0 |}{|a_{m-1} \cdots a_0 |} > \frac 12.
$$
Combining this inequality and \eqref{eq.gettingthere}, with $i$,
$j$ replaced by $j$, $m$, we find
$$
\frac{|b_{j-1} \cdots b_0 |}{|a_{j-1} \cdots a_0 |} =
\frac{|b_{m-1} \cdots b_0 |}{|a_{m-1} \cdots a_0 |} \; / \;
\frac{|b_{m-1} \cdots b_j |}{|a_{m-1} \cdots a_j |} >
\frac{1/2}{K},
$$
which is the remaining inequality in \eqref{e.supersonic}.

\smallskip

Now, combining lemma~\ref{l.4} with \eqref{e.supersonic} and
$\ang(v_s,w_s)\geq \alpha$, we get
$$
\frac{\| \DD^{(j)} \|}{\mm( \DD^{(j)} )} \leq \frac{8K}{\sin^2
\alpha}\,.
$$
Moreover, $ Df^j_y / X_0 = \pi_j \circ \DD^{(j)} \circ
\pi_0^{-1}$. So, using the relation \eqref{eq.angXY},
$$
\frac{\| Df^j_y / X_0 \|}{\mm( Df^j_y / X_0 )} \leq \frac{\| \pi_j
\|}{\mm( \pi_j )} \cdot \frac{\| \DD^{(j)} \|}{\mm( \DD^{(j)} )}
\cdot \frac{\| \pi_0 \|}{\mm( \pi_0 )} \le \frac{8K}{\sin^6
\alpha}.
$$
This finishes the proof of lemma~\ref{l.eccentricity}.
\end{proof}

The proof of proposition~\ref{p.geom} is now complete.
\end{proof}

\section{Proof of theorems~\ref{t.vol} and~\ref{t.vol.continuity}}\label{s.t.vol}

Let us define some useful invariant sets. Given $f \in \Diff$, let
$\mathcal{O}(f)$ be the set of the regular points, in the sense of
the theorem of Oseledets. Given $p\in\{1,\ldots,d-1\}$ and $m \in
\N$, let $\DD_p(f,m)$ be the set of points $x$ such that there is
an $m$-dominated splitting of index $p$ along the orbit of $x$.
That is, $x \in \DD_p (f, m)$ if and only if there exists a
splitting $T_{f^n x} M = E_n \oplus F_n$ ($n\in \Z$) such that for
all $n\in \Z$, $\dim E_n = p$, $Df_{f^n x}(E_n) = E_{n+1}$,
$Df_{f^n x}(F_n) = F_{n+1}$ and
$$
\frac{\| Df^m_{f^n(x)}|_{F_n} \|}{\mm(Df^m_{f^n(x)}|_{E_n})} \leq
\frac 12.
$$
By section~\ref{ss.basicdominated}, $\DD_p (f,m)$ is a closed set.
Let
\begin{align*}
\Gamma_p (f, m) &= M \minus \DD_p(f,m), \\
\Gamma_p^\sharp (f, m) &= \big\{ x \in \Gamma_p(f,m) \cap \mathcal{O}(f); \;
\lambda_p(f,x) > \lambda_{p+1}(f,x) \, \big\}, \\
\Gamma_p^* (f, m) &= \big\{ x \in \Gamma_p^\sharp(f,m); \;
\text{$x$ is not periodic} \big\}.
\end{align*}
Define also
$$
\Gamma_p (f, \infty) = \bigcap_{m\in \N} \Gamma_p(f, m)
\quad \text{and} \quad
\Gamma_p^\sharp (f, \infty) = \bigcap_{m\in \N} \Gamma_p^\sharp(f, m).
$$
It is clear that all these sets are invariant under $f$.

\begin{lemma} \label{l.aperiodic}
For every $f$ and $p$, the set $\Gamma_p^\sharp (f, \infty)$ contains
no periodic points of $f$.
In other words,
$\bigcap_{m\in \N}
\big( \Gamma_p^\sharp(f, m) \minus \Gamma_p^*(f, m) \big)
= \varnothing$.
\end{lemma}

\begin{proof}
Suppose that $x \in \mathcal{O}(f)$ is periodic,
say, of period $n$, and $\lambda_p(f,x) > \lambda_{p+1}(f,x)$.
The eigenvalues of $Df_x^n$ are $\nu_1,\ldots,\nu_d$,
with $| \nu_i | = e^{n\lambda_i(f,x)}$.
Let $E$ (resp. $F$) be the sum of the eigenspaces
of $Df^n_x$ associated to the eigenvalues $\nu_1,\ldots,\nu_p$
(resp. $\nu_{p+1},\ldots,\nu_d$).
Then the splitting $T_x M = E \oplus F$ is $Df^n_x$-invariant.
Spreading it along the orbit of $x$, we obtain a dominated splitting.
That is, $x \in \DD_p(f,m)$ for some $m \in \N$,
and so $x \not\in \Gamma_p^\sharp (f, \infty)$.
\end{proof}

\subsection{Lowering the norm along an orbit segment}

Recall that we write $\Lambda_p(f,x)=\lambda_1(x)+\cdots+\lambda_p(x)$
for each $p=1, \ldots, d$.

\begin{proposition}  \label{p.lower}
Let $f \in \Diff$, $\eps_0>0$, $\kappa >0$, $\delta>0$ and
$p\in\{1,\ldots,d-1\}$.
Then, for every sufficiently large $m\in\N$, there exists a measurable
function $N: \Gamma_p^*(f,m) \to \N$ with the following properties:
For almost every $x\in \Gamma_p^*(f,m) $ and every $n \geq N(x)$ there
exists an $(\eps_0, \kappa)$-realizable sequence
$\{\hL_{0}^{(x,n)},\ldots ,\hL_{n-1}^{(x,n)}\}$ at $x$ of
length $n$ such that
$$
\frac{1}{n} \log \| \wp ( \hL_{n-1}^{(x,n)} \cdots \hL_{0}^{(x,n)} ) \|  \, \leq \,
\frac{\Lambda_{p-1}(f,x) + \Lambda_{p+1}(f,x)}{2} + \delta.
$$
\end{proposition}

\begin{proof} Fix $f$, $\eps_0$, $\kappa$, $\delta$ and $p$.
For clearness, we divide the proof into two parts:

\paragraph{Part 1:} Definition of $N(\cdot)$ and the sequence $\hL_j^{(x,n)}$.
\smallskip

Fix $f$, $\eps_0$, $\kappa$, $\delta$ and $p$.
Assume $m\in \N$ is sufficiently large so that
the conclusion of proposition~\ref{p.geom} holds for
$f$, $\eps_0$ and $\meio \kappa$ (in the place of $\kappa$).
To simplify the notation, let $\Gamma = \Gamma^*(p,m)$.
We may suppose that $\mu(\Gamma) >0$, otherwise there is nothing to prove.
Consider the splitting $T_\Gamma M = E \oplus F$, where $E$ is the sum of
the Oseledets subspaces corresponding to the first $p$ Lyapunov
exponents $\lambda_1 \geq \cdots \geq \lambda_p$ and $F$ is the
sum of the subspaces corresponding to the other exponents
$\lambda_{p+1} \geq \cdots \geq \lambda_d$. This makes sense since
$\lambda_p > \lambda_{p+1}$ on $\Gamma$. Let $A \subset \Gamma$ be
the set of points $y$ such that the non-domination condition
\eqref{e.not dominated} holds. By definition of $\Gamma=\Gamma_p^*(f,m)$,
\begin{equation} \label{e.A Gamma}
\Gamma = \bigcup_{n \in \Z} f^n(A).
\end{equation}

Let $\lambda_i^{\wedge p}(x)$, $1 \leq  i \leq \binom{d}{p}$
denote the Lyapunov exponents of the cocycle $\wp(Df)$ over $f$,
in non-increasing order.
Let $V_x$ denote the Oseledets subspace associated to
the upper exponent $\lambda_1^{\wedge p}(x)$ and let
$H_x$ be the sum of all other Oseledets subspaces.
This gives us a splitting $\wp(TM) = V \oplus H$.
By proposition~\ref{p.oseledets exterior}, we have
\begin{equation*}\begin{aligned}
\lambda_1^{\wedge p}(x) &= \lambda_1(x) + \cdots +\lambda_{p-1}(x) + \lambda_p(x), \\
\lambda_2^{\wedge p}(x) &= \lambda_1(x) + \cdots +\lambda_{p-1}(x) + \lambda_{p+1}(x).
\end{aligned}\end{equation*}
If $x\in \Gamma$ then $\lambda_p(x) > \lambda_{p+1}(x)$ and so
$\lambda_1^{\wedge p}(x) >\lambda_2^{\wedge p}(x)$.
That is, the subspace $V_x$ is one-dimensional.

For almost every $x\in \Gamma$, Oseledets' theorem gives $Q(x) \in
\N$ such that for all $n \geq Q(x)$, we have:
\begin{itemize}
\item $\frac 1n \log \frac{\| \wp(Df_x^n) v \|}{\| v \|} < \lambda_1^{\wedge p}(x) + \delta$
for every $v \in V_x \minus \{0\}$;
\item $\frac 1n \log \frac{\| \wp(Df_x^n) w \|}{\| w \|} < \lambda_2^{\wedge p}(x) + \delta$
for every $w \in H_x \minus \{0\}$;
\item $\frac 1n \log \sin \ang (V_{f^n x}, H_{f^n x}) > - \delta$.
\end{itemize}
For $q\in \N$, let $B_q = \{x \in \Gamma; \; Q(x) \leq q\}$.
Then $B_q \uparrow \Gamma$, that is,
the $B_q$ form a non-decreasing sequence and
their union is a full measure subset of $\Gamma$.
Define $C_0 = \varnothing$ and
\begin{equation} \label{e.Cq}
C_q = \bigcup_{n\in \Z} f^n(A \cap f^{-m}(B_q)).
\end{equation}
Since $f^{-m}(B_q) \uparrow \Gamma$ and \eqref{e.A Gamma}, we have
$C_q \uparrow \Gamma$.
To prove the proposition we must define the function $N$ on $\Gamma$.
We are going to define it on each of the sets $C_q \minus C_{q-1}$ separately.
From now on, let $q\in\N$ be fixed.

We need the following recurrence result, proved in \cite[lemma~3.12]{Bochi}.

\begin{lemma} \label{l.recurrence}
Let $f\in\Diff$. Let $A \subset M$ be a measurable set with
$\mu(A)>0$, and let $\Gamma = \cup_{n\in\Z} f^n(A)$. Fix any
$\gamma>0$. Then there exists a measurable function $N_0:\Gamma
\to \N$ such that for almost every $x\in\Gamma$, and for all $n\ge
N_0(x)$ and $t\in (0,1)$, there exists $\ell \in \{0,1,\ldots,n\}$
such that $t-\gamma \leq \ell/n \leq t+\gamma$ and $f^\ell(x) \in A$.
\end{lemma}

Let $c$ be a strict upper bound for $ \log \| \wp ( Df ) \|$ and
$\gamma = \min \{ c^{-1} \delta, 1/10 \}$. Using \eqref{e.Cq} and
lemma~\ref{l.recurrence}, we find a measurable function $N^{(q)}_0
: C_q \to \N$ such that for almost every $x\in C_q$, every $n \geq
N^{(q)}_0(x)$ and every $t\in(0,1)$ there is $\ell \in
\{0,1,\ldots,n\}$ with $|\ell/n - t| < \gamma$ and $f^\ell x \in A
\cap f^{-m}(B_q)$.
We define $N(x)$ for $x \in C_q \minus C_{q-1}$ as the least integer such that
$$
N(x)\geq\max\{N^{(q)}_0(x), \ 10Q(x), \ m\gamma^{-1},  \ \delta^{-1} \log[4/\sin \ang (V_x, H_x)]\}.
$$

Now fix a point $x \in C_q \minus C_{q-1}$ and $n \geq N(x)$.
We will now construct the sequence $\{ \hL_j^{(x,n)} \} $.
Since $n\ge N_0^{(q)}(x)$, there exists $\ell\in\N$ such that
$$
\left| \frac{\ell}{n} - \frac{1}{2} \right| < \gamma
\quad\text{and}\quad
y = f^{\ell} (x) \in A \cap f^{-m}(B_q).
$$
Since $y \in A$, where the non-domination condition \eqref{e.not dominated} holds,
proposition~\ref{p.geom} gives a sequence
$\{L_{0}, \ldots, L_{m-1}\}$ which is $(\eps_0, \meio
\kappa)$-realizable , such that there are non-zero vectors $v_0
\in E_y$, $w_0 \in F_{f^m(y)}$ for which
$$
L_{m-1} \dots L_{0} (v_0) = w_0\,.
$$
We form the sequence $\{\hL_{0}^{(x,n)}, \ldots, \hL_{n-1}^{(x,n)}\}$ of length $n$
by concatenating
$$
\{Df_{f^i(x)}; \; 0 \leq i < \ell\},
\quad
\{L_{0}, \ldots, L_{m-1}\},
\quad
\{Df_{f^i(x)}; \; \ell+m \le i < m\}.
$$
According to parts 1 and 2 of lemma~\ref{l.basicproperties}, the concatenation is
an $(\eps_0, \kappa)$-realizable sequence at~$x$.

\paragraph{Part 2:} Estimation of $ \| \wp ( \hL_{n-1}^{(x,n)} \cdots \hL_0^{(x,n)} ) \|$.
\smallskip

Write $\wp ( \hL_{n-1}^{(x,n)} \cdots \hL_0^{(x,n)} ) = D_1 \LL D_0 $, with
$D_0 = \wp (Df^\ell_x)$,\
$D_1 = \wp ( Df^{n- \ell- m}_{f^{\ell + m}(x)} )$, and
$\LL = \wp ( L_{m-1} \cdots L_{0} )$. The key observation is:

\begin{lemma}    \label{l.V goes to H}
The map $\LL:\wp(T_y M) \to \wp(T_{f^m(y)} M)$ satisfies
$\LL (V_y) \subset H_{f^m(y)}$.
\end{lemma}

\begin{proof}
Proposition~\ref{p.oseledets exterior} describes
the spaces $V$ and $H$. Let $z \in \Gamma$ and consider a basis
$\{e_1(z), \ldots e_d(z)\}$ of $T_zM$ such that
$$
e_i(x) \in E^j_x
\quad \text{for $\dim E^1_x + \dots + \dim E^{j-1}_x
< i \leq \dim E^1_x + \dots + \dim E^j_x$.}
$$
Then $V_z$ is the space generated by $e_1(z) \wedge \cdots \wedge e_p(z)$ and
$H_z$ is generated by the vectors $e_{i_1}(z) \wedge \cdots \wedge e_{i_p}(z)$
with $1 \leq i_1 < \cdots <i_p \leq d$, $i_p > p$.
Also notice that $\{e_1(z), \ldots, e_p(z)\}$
and  $\{e_{p+1}(z), \ldots, e_d(z)\}$ are bases for the spaces
$E_z$ and $F_z$, respectively.
Consider the vectors $v_0 \in E_y$ and $w_0 = L (v_0) \in F_{f^m y}$,
where $L = L_{m-1} \dots L_{0}$.
There is $\nu \in \{1,\ldots,p\}$ such that
$$
\{ v_0, e_1(y), \ldots, e_{\nu-1}(y),e_{\nu+1}(y), \ldots, e_p(y) \}
$$
is a basis for $E_y$.
Therefore $V_y$ is generated by the vector
$$
v_0 \wedge e_1(y) \wedge \cdots \wedge e_{\nu-1}(y) \wedge e_{\nu+1}(y) \wedge \cdots \wedge e_p(y),
$$
which is mapped by $\LL =\wp (L)$ to
\begin{equation} \label{e.ble}
w_0 \wedge Le_1(y) \wedge \cdots \wedge Le_{\nu-1}(y) \wedge
Le_{\nu+1}(y) \wedge \cdots \wedge Le_p(y),
\end{equation}
Write $w_0$ as a linear combination of vectors $e_{p+1}(f^m(y))$,
\ldots, $e_d(f^m(y))$ and write each $L e_i(y)$ as a linear
combination of vectors $e_1(f^m(y))$, \ldots, $e_d(f^m(y))$.
Substituting in \eqref{e.ble}, we get a linear
combination of $e_{i_1}(f^m(y)) \wedge \cdots \wedge e_{i_p}(f^m(y))$
where $e_1(f^m(y)) \wedge \cdots \wedge e_p(f^m(y))$ does
\emph{not} appear. This proves that the vector in \eqref{e.ble}
belongs to $H_{f^m(y)}$.
\end{proof}

To carry on the estimates, we introduce a more convenient norm:
For $x_0, x_1 \in \Gamma$ we represent a linear map $T :
\wp(T_{x_0} M) \to \wp(T_{x_1} M)$ by its matrix
$$
T = \begin{pmatrix}
T^{++} & T^{+-} \\
T^{-+} & T^{--}
\end{pmatrix}
$$
with respect to the splittings
$T_{x_0} M = V_{x_0} \oplus H_{x_0} $ and $T_{x_1} M = V_{x_1} \oplus H_{x_1}$.
Then we define
$$
\|T\|_{\max} = \max\Big\{\|T^{++}\|, \|T^{+-}\|, \|T^{-+}\|, \|T^{--}\|\Big\}.
$$
The following elementary lemma relates this norm with the original one $\|T\|$
(that comes from the metric in $\wp(T_\Gamma M)$ ).

\begin{lemma} \label{l.twonorms}
Let $\theta_{x_0}=\ang(V_{x_0},H_{x_0})$ and
    $\theta_{x_1}=\ang(V_{x_1},H_{x_1})$.
Then:
\begin{enumerate}
\item $\|T\|        \leq 4 \, (\sin\theta_{x_0})^{-1} \, \|T\|_{\max}$;
\item $\|T\|_{\max} \leq      (\sin\theta_{x_1})^{-1} \, \|T\| $.
\end{enumerate}
\end{lemma}

\begin{proof}
Let $v=v_+ + v_- \in V_{x_0} \oplus H_{x_0}$. We have $\|v_*\| \le
\|v\| /\sin\theta_{x_0}$ for $*=+$ and $*=-$. So
$$
\|T v\| \le \|T^{++} v_+\| + \|T^{++} v_-\| + \|T^{--} v_+\| +
\|T^{--} v_-\| \le 4 \|T\|_{\max} \|v\|/\sin\theta_{x_0}.
$$
This proves part 1. The proof of part 2 is similar. Let $v_+\in V_{x_0}$.
Its image splits as $Tv_+ = T^{++}v_+ + T^{-+}v_+ \in V_{x_1} \oplus H_{x_1}$.
Hence,
$$
\|T^{*+}v_+\| \le \|T v_+\| (\sin\theta_{x_1})^{-1}
\le \|T\| \, \|v_+\| (\sin\theta_{x_1})^{-1}
$$
for $*=+$ and $*=-$.
Together with a corresponding estimate for $T^{*-}v_-$, $v_-\in H_{x_0}$,
this gives part 2.
\end{proof}

For the linear maps we were considering, the matrices have the form:
$$
D_i =
\begin{pmatrix}
D_i^{++} & 0 \\
0        & D_i^{--}
\end{pmatrix}, \
i=0,1,
\quad \text{and} \quad
\LL =
\begin{pmatrix}
0        & \LL^{+-} \\
\LL^{-+} & \LL^{--}
\end{pmatrix}:
$$
$D_i^{+-}=0=D_i^{-+}$ because $V$ and $H$ are $\wp(Df)$-invariant,
and $\LL^{++}=0$ because of lemma~\ref{l.V goes to H}.
Then
\begin{equation} \label{e.product}
\wp( \hL_{n-1} \cdots \hL_0 )=
\begin{pmatrix}
0                          & D_1^{++} \LL^{+-} D_0^{--} \\
D_1^{--} \LL^{-+} D_0^{++} & D_1^{--} \LL^{--} D_0^{--}
\end{pmatrix}.
\end{equation}

\begin{lemma} \label{l.est1}
For $i=0,1$, $x\in C_q \minus C_{q-1}$ and $n\ge N(x)$,
$$
\log \| D_i^{++} \| < \meio n (\lambda^{\wedge p}_1(x) + 5\delta)
\quad\text{and}\quad
\log \| D_i^{--} \| < \meio n (\lambda^{\wedge p}_2(x) + 5\delta).
$$
\end{lemma}
\begin{proof}
Since $\ell > (\frac 12 - \gamma) n >  \frac{1}{10}n \geq \frac{1}{10} N(x) \geq Q(x)$, we have
\begin{align*}
\log \| D_0^{++} \| = \log \| \wp(Df^\ell_x)|_{V_x} \|
&< \ell (\lambda^{\wedge p}_1(x) + \delta), \\
\log \| D_0^{--} \| = \log \| \wp(Df^\ell_x)|_{H_x} \|
&< \ell (\lambda^{\wedge p}_2(x) + \delta).
\end{align*}
Let $\lambda$ be either $\lambda^{\wedge p}_1(x)$ or $\lambda^{\wedge p}_2(x)$.
Using $\gamma\lambda < \gamma c \leq \delta$ and $\gamma<1$, we find
$$
\ell (\lambda + \delta) < n (\meio + \gamma ) (\lambda + \delta)
                      < n (\meio \lambda + \meio \delta + \delta + \delta)
              = \meio n (\lambda + 5\delta).
$$
This proves the case $i=0$.
We have
$n- \ell - m > n(\frac 12 - \gamma) - n\gamma > \frac{1}{10} n \geq Q(x) \geq q$.
Also $f^\ell(x) \in f^{-m}(B_q)$, and so $Q(f^{\ell + m}(x)) \leq q$.
Therefore
\begin{align*}
\log \| D_1^{++} \| = \log \| \wp(Df^{n-\ell-m}_{f^{\ell + m} x})|_{V_{f^{\ell + m} x}} \|
&< (n-\ell-m) (\lambda^{\wedge p}_1(x) + \delta), \\
\log \| D_1^{--} \| = \log \| \wp(Df^{n-\ell-m}_{f^{\ell + m} x})|_{H_{f^{\ell + m} x}} \|
&< (n-\ell-m) (\lambda^{\wedge p}_2(x) + \delta).
\end{align*}
As before,
$(n-\ell-m)(\lambda+\delta) < n (\meio +\gamma) (\lambda+\delta)\ \meio n(\lambda + 5\delta)$.
This proves case $i=1$.
\end{proof}

\begin{lemma} \label{l.est2}
$\log \| \LL \|_{\max} < 2n\delta$.
\end{lemma}
\begin{proof}
Since the sequence $\{ L_0, \ldots, L_{m-1} \}$ is realizable,
each $L_j$ is close to the value of $Df$ at some point.
Therefore we may assume that $\log \| \wp(L_j) \| < c$.
In particular, $\log \| \LL \| < mc \le nc\gamma \le n\delta$.
We have $\ell + m \ge n(\frac 12 -\gamma) \ge \frac {1}{10} n \ge Q(x)$.
So $\log [1/\sin \ang(V_{f^{\ell+m} x}, H_{f^{\ell+m} x})] < \delta$
and, by part 2 of lemma~\ref{l.twonorms},
$\log \| \LL \|_{\max} < 2n\delta$.
\end{proof}

Using lemmas~\ref{l.est1} and \ref{l.est2}, we bound each of the entries in~\eqref{e.product}:
\begin{align*}
\log \| D_1^{++} \LL^{+-} D_0^{--} \| &<
\meio n (\lambda^{\wedge p}_1(x) + \lambda^{\wedge p}_2(x) + 14\delta)
\\
\log \| D_1^{--} \LL^{-+} D_0^{++} \| &<
\meio n (\lambda^{\wedge p}_1(x) + \lambda^{\wedge p}_2(x) + 14\delta)
\\
\log \| D_1^{--} \LL^{--} D_0^{--} \| &<
\meio n (2 \lambda^{\wedge p}_2(x) + 14\delta)
\end{align*}
The third expression is smaller than either of the first two,
so we get
$$
\log\| \wp( \hL_{n-1} \cdots \hL_0 ) \|_{\max} <
n \Big(\frac{\lambda^{\wedge p}_1(x) + \lambda^{\wedge p}_2(x)}{2} + 7\delta\Big).
$$
Therefore, by part 1 of lemma~\ref{l.twonorms}
and $\log [ 4 / \sin \ang(V_x, H_x) ] < n\delta$,
$$
\log\| \wp( \hL_{n-1} \cdots \hL_0 ) \|
< n \Big(\frac{\lambda^{\wedge p}_1(x) + \lambda^{\wedge p}_2(x)}{2} + 8\delta\Big).
$$
We also have
$\lambda^{\wedge p}_1(x) + \lambda^{\wedge p}_2(x) = \Lambda_{p-1}(f,x) + \Lambda_{p+1}(f,x)$.
This proves proposition~\ref{p.lower}
(replace $\delta$ with $\delta/8$ along the proof).
\end{proof}

\subsection{Globalization}
The following proposition renders global the construction of proposition~\ref{p.lower}:

\begin{proposition} \label{p.global}
Let $f \in \Diff$, $\eps_0>0$, $p\in\{1,\ldots,d-1\}$ and $\delta>0$.
Then there exist $m\in\N$ and a diffeomorphism $g \in \UU(f, \eps_0)$
that equals $f$ outside the open set $\Gamma_p(f,m)$ and such that
$$
\int_{\Gamma_p(f,m)} \Lambda_p(g,x) \, d\mu(x) <
\delta +
\int_{\Gamma_p(f,m)} \frac{\Lambda_{p-1}(f,x) + \Lambda_{p+1}(f,x)}{2} \, d\mu(x).
$$
\end{proposition}

We need some preparatory terminology:

\begin{definition}
Let $f \in \Diff$.
An \emph{$f$-tower} (or simply \emph{tower})
is a pair of measurable sets $(T, T_{\mathrm{b}})$
such that there is a positive integer $n$,
called the \emph{height} of the tower, such that the sets
$T_{\mathrm{b}}, f(T_{\mathrm{b}}), \ldots,f^{n-1}(T_{\mathrm{b}})$
are pairwise disjoint and their union is $T$.
$T_{\mathrm{b}}$ is called the \emph{base} of the tower.

An \emph{$f$-castle} (or simply \emph{castle})
is a pair of measurable sets $(Q, Q_{\textrm{b}})$
such that there exists a finite or countable family of pairwise disjoint
towers $(T_i, T_{i\mathrm{b}})$ such that
$Q = \bigcup T_i$ and $Q_{\textrm{b}} = \bigcup T_{i\mathrm{b}}$.
$Q_{\mathrm{b}}$ is called the \emph{base} of the castle.

A castle $(Q, Q_{\textrm{b}})$ is a \emph{sub-castle} of a castle $(Q', Q'_{\textrm{b}})$
if $Q_{\textrm{b}} \subset Q'_{\textrm{b}}$ and
for every point $x \in Q_{\textrm{b}}$, if $n$ (respectively $n'$) denotes the height
of a tower of $(Q, Q_{\textrm{b}})$ (respectively $(Q', Q'_{\textrm{b}})$)
that contains $x$, then $n=n'$. In particular, $Q \subset Q'$.
\end{definition}

We shall frequently omit reference to the base of a castle $Q$ in our notations.

\begin{definition}
Given $f\in \Diff$ and a positive measure set $A\subset M$,
consider the return time $\tau:A\to \N$
defined by $\tau (x)=\inf \{n\geq 1;\;f^n(x)\in A\}$. If we denote
$A_n=\tau ^{-1}(n)$ then $T_n=A_n\cup f(A_n)\cup \cdots \cup
f^{n-1}(A_n)$ is a tower. Consider the castle $Q$, with base $A$,
given by the union of the towers $T_n$.
$Q$ is called the \emph{Kakutani castle} with base $A$.
\end{definition}

Note that $Q = \bigcup_{n\in \Z}f^n(A)$ mod 0, in particular the
set $Q$ is invariant.

\begin{proof}[Proof of proposition~\ref{p.global}]
Let $f$, $\eps_0$, $p$ and $\delta$ be given. For simplicity, we write
$$
\phi(x) = \frac{\Lambda_{p-1}(f,x) + \Lambda_{p+1}(f,x)}{2}\,.
$$

\paragraph{Step 1:} Construction of families of castles $\hQ_i\supset Q_i$\,.
\smallskip

Let $\kappa = \delta^2$.
Take $m\in \N$ large enough so that the conclusion of
proposition~\ref{p.lower} holds:
there exists a measurable function $N: \Gamma_p^*(f,m) \to \N$
such that for a.e. $x\in \Gamma_p^*(f,m)$ and
every $n \geq N(x)$ there exists a $(\eps_0, \kappa)$-realizable sequence
$\{\hL_{0}^{(x,n)},\ldots ,\hL_{n-1}^{(x,n)}\}$
at $x$ of length $n$ such that
\begin{equation}\label{e.seq}
\frac{1}{n} \log \| \wp ( \hL_{n-1}^{(x,n)} \cdots \hL_{0}^{(x,n)} ) \|  \leq
\phi(x) + \delta.
\end{equation}
We shall also (see lemma~\ref{l.aperiodic}) assume that $m$ is large enough so that
\begin{equation} \label{e.aperiodic}
\mu \big( \Gamma_p^\sharp(f, m) \minus \Gamma_p^*(f, m) \big) < \delta.
\end{equation}
Let $C > \sup_{g\in\UU(f,\eps_0)} \sup_{y\in M} \log \| \wp(Dg_y)\|$ and
$\ell = \lceil C/\delta \rceil$. For $i=1,2,\ldots,\ell$, let
$$
Z^i = \{ x \in \Gamma_p^*(f,m);\; (i-1)\delta \leq \phi(x) < i\delta \}.
$$
Each $Z_i$ is an $f$-invariant set.
Since $\phi < C$, we have $\Gamma_p^*(f,m) = \bigsqcup_{i=1}^{\ell} Z^i$.
Define the sets
$Z_n^i = \{ x \in Z^i; \; N(x)\leq n \} $ for $n\in \N$ and $1\leq i\leq\ell$.
Obviously, $Z_n^i \uparrow Z^i$ when $n\to\infty$.
Fix $H\in \N$ such that, for all $i=1,2,\ldots,\ell$,
\begin{equation}\label{e.Z minus ZH}
\mu (Z^i \minus Z_H^i) < \delta^2 \mu(Z^i).
\end{equation}
Using the fact that $\Lambda_p(f)$ equals $\phi$ in the $f$-invariant set
$\Gamma_p(f,m) \minus \Gamma_p^\sharp(f,m)$, and proposition~\ref{p.formula},
we may also assume that $H$ is large enough so that
\begin{equation}\label{e.banal}
\int_{\Gamma_p(f,m) \minus \Gamma_p^\sharp(f,m)}
\frac 1n \log \| \wp (Df^n) \|
< \delta +
\int_{\Gamma_p(f,m) \minus \Gamma_p^\sharp(f,m)} \phi
\end{equation}
for all $n\ge H$.

A measure preserving transformation is \emph{aperiodic} if the set of periodic
points has zero measure.The following result was proved in \cite[Lemma~4.1]{Bochi}:

\begin{lemma} \label{l.maximal}
For every aperiodic invertible measure preserving transformation $f$ on a
probability space $X$, every subset $U\subset X$ of positive measure, and
every $n\in \N$, there exists a positive measure set $V\subset U$
such that the sets $V$, $f(V)$,\ldots, $f^{n}(V)$ are two-by-two disjoint.
Besides, $V$ can be chosen maximal \emph{in the measure-theoretical sense}:
no set that includes $V$ and has larger measure than $V$ has the stated properties.
\end{lemma}

By definition of the set $\Gamma_p^*(f,m)$, the map $f:\Gamma_p^*(f,m)\to\Gamma_p^*(f,m)$ is
aperiodic. So, by lemma~\ref{l.maximal}, for each $i$ there is
$B^i \subset Z_H^i$ such that $B^i, f(B^i),\ldots, f^{H-1}(B^i)$
are two-by-two
disjoint and such that $B^i$ is maximal for these properties
(in the measure-theoretical sense).
Consider the following $f$-invariant set:
$$
\hQ^i = \bigcup_{n\in \Z}f^n(B^i).
$$
$\hQ^i$ is the Kakutani castle with base $B^i$.
It is contained in $Z^i$ and, by the maximality of $B^i$,
it contains $Z_H^i$ up to a zero measure subset.
Thus, by \eqref{e.Z minus ZH},
\begin{equation} \label{e.Z minus hQ}
\mu(Z^i \minus \hQ^i) < \delta^2 \mu(Z^i).
\end{equation}
Let $Q^i \subset \hQ^i$ be the sub-castle
consisting of all the towers of $\hQ^i$ with heights at most $3H$
floors.
The following is a key property of the construction:

\begin{lemma} \label{l.hQ minus Q}
For each $i=1,2,\ldots,\ell$, we have
$\mu (\hQ^i \minus Q^i) \leq 3 \mu(Z^i \minus Z_H^i)$.
\end{lemma}

\begin{proof}
We split the castle $\hQ^i$ into towers as
$\hQ^i = \bigsqcup_{k=H}^{\infty} T_k^i$ where
$B^i = \bigsqcup_{k=H}^{\infty}B_k^i$ is the base $\hQ^i_{\mathrm{b}}$ and
$T_k^i = \bigsqcup_{j=0}^{k-1}f^j(B_k^i)$ is the tower with base $B_k^i$ and of
height $k$ floors.
Take $k\geq 2H$ and $H\leq j\leq k-H$.
The sets $f^j(B_k^i),\ldots f^{j+H-1}(B_k^i)$ are disjoint and do not intersect
$B^i\sqcup \cdots \sqcup f^{H-1}(B^i)$.
Since $B^i$ is maximal, we conclude that
$$
k\geq 2H \text{ and }\ H\leq j\leq k-H \quad\Rightarrow\quad \mu(f^j(B_k^i)\cap Z_H^i)=0
$$
(otherwise we could replace $B^i$ with $B^i \sqcup (f^j(B_k^i) \cap Z_H^i)$,
contradicting the maximality of $B^i$.)
Thus
$$
k\geq 2H \quad\Rightarrow\quad \mu (T_k^i \minus Z_H^i)\geq
\sum_{j=H}^{k-H} \mu(f^j(B_k^i)) = \frac{k-2H+1}{k}\mu(T_k^i).
$$
In particular,
$$
k\geq 3H+1 \Rightarrow \mu (T_k^i \minus Z_H^i)>\frac{1}{3}\mu (T_k^i)
$$
and so
\begin{multline*}
\mu (\hQ^i \minus Q^i)=
\sum_{k=3H+1}^{\infty} \mu (T_k^i)\leq
\sum_{k=3H+1}^{\infty} 3\mu (T_k \minus Z_H^i) \\=
3\mu \left(\bigsqcup_{k=3H+1}^{\infty}T_k^i \minus Z_H^i \right) \leq
3\mu (Z^i \minus Z_H^i),
\end{multline*}
as claimed.
\end{proof}

\smallskip
\noindent{\bf Step 2:} Construction of the diffeomorphism $g$.
\smallskip

\begin{lemma}\label{l.ufa}
For almost every $x\in \Gamma_p^*(f,m)$ and every $n \geq N(x)$,
there exists $r>0$ such that for every ball $U = B_{r'}(x)$ with
$0<r'<r$ there exist $h \in \UU(f,\eps_0)$ and a measurable set
$K\subset B_{r'}(x)$ such that
\begin{itemize}
\item[(i)]   $h$ equals $f$ outside $\sqcup_{j=0}^{n-1} f^j(B_{r'}(x))$;
\item[(ii)]  $\mu(K) > (1-\kappa) \mu(B_{r'}(x))$;
\item[(iii)] if $y\in K$ then
$\frac 1n \log \| \wp (Dh^n_y) \| <  \phi(x) + 2\delta$.
\end{itemize}
\end{lemma}

\begin{proof}
Fix $x$ and $n\geq N(x)$. Recall the point $x$ is not periodic.
Let $\gamma>0$ be very small.
Since the sequence $\{\hL_{j}^{(x,n)}\}$ given by proposition~\ref{p.lower}
is $(\kappa,\eps_0)$-realizable, there exists $r>0$ such that
for every ball $U= B_{r'}(x)$ with $0<r'<r$ there exists $h \in \UU(f,\eps_0)$
satisfying condition (i) above and there exists
$K\subset B_{r'}(x)$ satisfying condition (ii) and
$$
y \in K \text{ and } \ 0\leq j \leq n-1 \quad  \Rightarrow \quad
\| Dh_{h^j y} - \hL_{j}^{(x,n)} \| < \gamma.
$$
Taking $\gamma$ small enough, this inequality and \eqref{e.seq} imply
$$
y \in K  \  \Rightarrow \
\frac 1n \log \| \wp (Dh^n_y) \|
 < \frac 1n \log \| \wp ( \hL_{n-1}^{(x,n)} \cdots \hL_{0}^{(x,n)} )  \| + \delta
\le  \phi(x) + 2\delta,
$$
as claimed in the lemma.
\end{proof}

\begin{lemma}\label{l.g}
Fix $\gamma>0$. There exists $g \in \UU(f, \eps_0)$ and
for each $i=1,2,\ldots,\ell$ there exist
a $g$-castle $U^i$ and a $g$-sub-castle $K^i$
such that:
\begin{itemize}
\item[(i)] the $U^i$ are open, pairwise disjoint, and contained in $\Gamma_p(f,m)$;

\item[(ii)] $\mu(U^i \minus Q^i)< 2 \gamma \mu(Z^i)$ and  $\mu(Q^i \minus U^i) < 2 \gamma \mu(Z^i)$;

\item[(iii)] $\mu (U^i \minus K^i) < 2 \kappa \mu(Z^i)$;

\item[(iv)] $g(U^i) = f(U^i)$ and $g$ equals $f$ outside $\bigsqcup_{i=1}^{\ell} U^i$;

\item[(v)] if $y$ is in base of $K^i$ and $n(y)$ is the height of the
tower of $K^i$ that contains $x$ then
$$
\frac{1}{n(y)} \log \| \wp (Dg^{n(y)}_y) \| < i \delta + 2\delta .
$$
\end{itemize}
\end{lemma}

\begin{proof}
By the regularity of the measure $\mu$, one can find
a compact sub-castle $J^i\subset Q^i$ such that
\begin{equation}\label{e.Q minus J}
\mu (Q^i \minus J^i) < \gamma \mu(\hQ^i).
\end{equation}
Since the $J^i$ are compact and disjoint we can find
\emph{open} pairwise disjoint castles $V^i$ such that
each $V^i$ contains $J^i$ as a sub-castle,
is contained in the open and invariant set $\Gamma_p(f,m)$,
and
\begin{equation}\label{e.V minus J}
\mu (V^i \minus J^i) < \gamma \mu(\hQ^i).
\end{equation}

For each $x \in J^i_{\mathrm{b}}$, let $n(x)$ be the height of the
tower that contains $x$. $J^i_{\mathrm{b}}$ is contained in
$Z^i_H$, so $N(x) \leq H \leq n(x)$. Let $r(x)>0$ be the radius
given by lemma~\ref{l.ufa}, with $n=n(x)$. This is defined for
almost every $x\in J^i_{\mathrm{b}}$. Reducing $r(x)$ if needed, we
suppose that the ball $\overline{B}_{r(x)}(x)$ is contained in the
base of a tower in $V^i$ (with the same height).

Using Vitali's covering lemma\footnote{First, cover the basis $J^i_{\mathrm{b}}$ of
the castle by chart domains.},
we can find a finite collection of disjoint
balls $U^i_k=B_{r_{k,i}}(x_{k,i})$,
with $x_{k,i}\in J^i_{\mathrm{b}}$ and $0<r_{k,i}<r(x_{k,i})$,
such that
\begin{equation} \label{e.ast}
\mu \Big( J^i_{\mathrm{b}} \minus \bigsqcup_k \overline{U^i_k} \Big)
< \gamma \mu(J^i_{\mathrm{b}}).
\end{equation}
Let $n_{k,i} = n(x_{k,i})$.
Notice that $n(x)=n_{k,i}$ for all $x^i\in U^i_k$.

Now we apply lemma~\ref{l.ufa} to each ball $U^i_k$.
We get, for each $k$, a measurable set
$K^i_k \subset U^i_k$ and a diffeomorphism $h_{k,i} \in \UU(f,\eps_0)$ such that
(in 3 we use that $x_{k,i} \in Z^i$)
\begin{enumerate}
\item $h_{k,i}$ equals $f$ outside the set  $\bigsqcup_{j=0}^{n^i_k-1}f^j(U^i_k)$;
\item $\mu(K^i_k) > (1-\kappa) \, \mu(U^i_k)$;
\item if $y\in K^i_k$ then
$\frac{1}{n_{k,i}} \log\| \wp(Dh_{k,i}^{n_{k,i}})_y \| < \phi(x_{k,i})+2\delta <
i\delta+2\delta$.
\end{enumerate}

Let $g$ be equal to $h_{k,i}$ in the set
$\bigsqcup_{j=0}^{n_k^i-1}f^j(U_i^k)$,
for each $i$ and $k$, and be equal to $f$ outside.
Since those sets are disjoint, $g\in \Diff$ is a well-defined diffeomorphism.
Each $h_{k,i}$ belongs to $\UU(f,\eps_0)$ and so $g$ also does.

Since each $U_k^i$ is contained in the base of a tower in the castle $V^i$,
$V^i$ is also a castle for $g$.
Let $U^i$ be the $g$-sub-castle of $V^i$ with base $\sqcup_k U_k^i$.
Analogously,
let $K^i$ be the $g$-sub-castle of $U^i$ with base $\sqcup_k K_k^i$.

It remains to prove claims~(ii) and~(iii) in the lemma.
Making use of the castle structures,
relation~\eqref{e.ast} and item~2 above imply, respectively,
\begin{equation}\label{e.begin}
\mu (J^i \minus U^i) < \gamma \mu (J^i)
\quad \text{and} \quad
\mu(U^i \minus K^i) < \kappa \mu(U^i).
\end{equation}
By~\eqref{e.V minus J} and $\hQ^i\subset Z^i$,
\begin{equation}\label{e.U minus Q}
\mu(U^i \minus Q^i) 
< \mu(V^i \minus J^i) < \gamma \mu(\hQ^i) \le\gamma\mu(Z^i).
\end{equation}
This implies the first part of item~(ii).
Combining the first part of~\eqref{e.begin} with~\eqref{e.Q minus J},
$$
\mu(Q^i \minus U^i) < \mu(Q^i \minus J^i) + \mu(J^i \minus U^i)
                    < 2\gamma\mu(\hQ^i)\le 2\gamma\mu(Z^i).
$$
This proves the second part of item~(ii). Finally, second inequality
in~\eqref{e.begin} and
$$
\mu(U^i) < \mu(Q^i) + \mu(U^i \minus Q^i) < (1+\gamma) \mu(\hQ^i) < 2\mu(Z^i).
$$
imply item (iii). The lemma is proved.
\end{proof}

\paragraph{Step 3:} Conclusion of the proof of proposition~\ref{p.global}.
\medskip

Let $U = \bigsqcup_{i=1}^\ell U^i$ and
$Q = \bigsqcup_{i=1}^\ell Q^i$ and $\hQ = \bigsqcup_{i=1}^\ell \hQ^i$.
Set $N = H\delta^{-1}$.
(Of course, we can assume that $\delta^{-1}\in \N$.)
Let
$$
G = \bigsqcup_{i=1}^\ell G^i
\quad\text{where}\quad
G^i = Z^i \cap \bigcap_{j=0}^{N-1} g^{-j}(K^i)
$$
for each $i=1,2,\ldots,\ell$. The next lemma means that
on $G$ we managed to reduce some time $N$ exponent:

\begin{lemma}\label{l.small in G}
If $x \in G^i$ then
$$
\frac{1}{N} \log \| \wp (Dg^N_x) \| < i\delta + (6C + 2)\delta .
$$
\end{lemma}

\begin{proof}
For $y\in K^i_{\mathrm{b}}$, let $n(y)$ be the height of the $g$-tower containing $y$.
Take $x\in G$, say $x\in G^i$.
Since the heights of towers of $K^i$ are less than $3H$, we
can write
$$
N = k_1 + n_1 + n_2 + \cdots + n_j + k_2
$$
so that $0\leq k_1,k_2<3H$, $1\leq n_1,\ldots ,n_j<3H$, and the
points
$$
x_1     = g^{k_1}(x), \
x_2     = g^{k_1+n_1}(x),\ \ldots, \
x_{j+1} = g^{k_1+n_1+\cdots +n_j}(x)
$$
are exactly the points of the orbit segment
$x, g(x),\ldots,g^{N-1}(x)$ which belong to $K^i_{\mathrm{b}}$.
Write
$$
\| \wp (Dg_x^N)\| \leq \| \wp (Dg_x^{k_1})\| \
\| \wp (Dg_{x_1}^{n_1})\| \  \cdots \  \| \wp (Dg_{x_j}^{n_j})\|  \
\| \wp (Dg_{x_{j+1}}^{k_2})\|.
$$
Using item~(v) of lemma~\ref{l.g}, and our choice of $N=H\delta^{-1}$,
we get
\begin{multline*}
\log \| \wp (Dg_x^N)\|
< k_1 C + (n_1 + \cdots + n_j)(i\delta + 2\delta) + k_2 C \\
< 6HC + N(i\delta + 2\delta)
< N(6C\delta + i\delta + 2\delta) ,
\end{multline*}
as claimed.
\end{proof}

We also use that $G$ covers most of $U \cup \Gamma_p^*(f,m)$, as
asserted by the next lemma.
Let us postpone the proof of this lemma for a while.

\begin{lemma} \label{l.big G}
Let $\gamma = \delta^2 / (\ell H)$ in lemma~\ref{l.g}. Then
$\mu \big( U \cup\Gamma_p^*(f,m) \minus G \big) < 12 \delta$.
\end{lemma}

Continuing with the proof of proposition~\ref{p.global}, write
$\psi(x) = \frac 1N \log \| \wp (Dg^N_x) \|$.
Since $g$ leaves invariant the set $\Gamma_p(f,m)$,
proposition~\ref{p.formula} gives
$$
\int_{\Gamma_p(f,m)} \Lambda_p(g)
\leq  \int_{\Gamma_p(f,m)} \psi .
$$
We split the integral on the right hand side as
\begin{align*}
\int_{\Gamma_p(f,m)} \psi & =
\int_{\Gamma_p(f,m) \minus (U \cup \Gamma_p^\sharp(f,m))} \psi +
\int_{(U \cup \Gamma_p^\sharp(f,m)) \minus G} \psi +
\int_{G} \psi
\\ & = (\mathrm{I})+(\mathrm{II})+(\mathrm{III}).
\end{align*}
Outside $U$, $g$ equals $f$ and so
$\psi$ equals $\frac 1N \log \| \wp(Df^N) \|$.
Thus
$$
(\mathrm{I}) \leq
\int_{\Gamma_p(f,m) \minus \Gamma_p^\sharp(f,m)} \frac 1N \log \| \wp(Df^N)\| <
\delta + \int_{\Gamma_p(f,m) \minus \Gamma_p^\sharp(f,m)} \phi,
$$
by~\eqref{e.banal}.
By lemma~\ref{l.big G} and \eqref{e.aperiodic},
$\mu \big( (U \cup \Gamma_p^\sharp(f,m)) \minus G \big) < 13\delta$.
Since $\psi<C$, we have $(\mathrm{II}) \leq 13C \delta$.
Using lemma~\ref{l.small in G},
$$
(\mathrm{III})
=    \sum_{i=1}^{\ell} \int_{G^i}  \psi
\leq \sum_{i=1}^{\ell} (i\delta + (6C+2)\delta) \mu(G^i)
<    (6C + 3)\delta + \sum_{i=1}^{\ell} (i-1)\delta \mu(G^i).
$$
Since $\phi \geq (i-1)\delta$ inside $Z^i \supset G^i$, we have
$$
(\mathrm{III})
< (6C + 3) \delta + \int_{\Gamma_p^*(f,m)} \phi.
$$
Summing the three terms, we get the conclusion of proposition~\ref{p.global}
(replace $\delta$ with $\delta /(18C+4)$ throughout the arguments):
$$
\int_{\Gamma_p(f,m)} \Lambda_p(g) <
(18C+4)\delta + \int_{\Gamma_p(f,m)} \phi.
$$
This completes the proof of the proposition, modulo proving lemma~\ref{l.big G}.
\end{proof}

\paragraph{Step 4:} Proof of lemma~\ref{l.big G}.
\smallskip

The following observations will be useful in the proof:
If $X\subset M$ is a measurable set and $N\in\N$, then
\begin{equation}\label{eq.union}
\mu \Big( \bigcup_{j=0}^{N-1}g^{-j}(X)\Big) \leq
\mu (X)+ (N-1)\mu \big(g^{-1}(X) \minus X\big) .
\end{equation}
Moreover,
$\mu \big( g^{-1}(X) \minus X\big) =\mu \big( X \minus g^{-1}(X) \big)$.

\begin{proof} [Proof of lemma~\ref{l.big G}]
We shall prove first that
\begin{equation} \label{e.hQ minus G}
\mu (\hQ^i \minus G^i) < 10 \delta \mu(Z^i).
\end{equation}
Since $\hQ^i \subset Z^i$, we have
$\hQ^i \minus G^i \subset \hQ^i \cap \bigcup_{j=0}^{N-1} g^{-j}(M \minus K^i)$.
Substituting
$$M \minus K^i \subset (U^i \minus K^i) \cup (Q^i \minus U^i) \cup
(\hQ^i \minus Q^i) \cup (M \minus \hQ^i),$$
we obtain
\begin{multline*}
\hQ^i \minus G^i \subset
\bigcup_{j=0}^{N-1} g^{-j}(U^i \minus K^i)   \ \cup \
\bigcup_{j=0}^{N-1} g^{-j}(Q^i \minus U^i)   \ \cup \
\bigcup_{j=0}^{N-1} g^{-j}(\hQ^i \minus Q^i) \\ \cup \
\Bigg[\hQ^i \cap
\bigcup_{j=1}^{N-1} g^{-j}(M \minus \hQ^i)
\Bigg] =
(\mathrm{I}) \cup (\mathrm{II}) \cup (\mathrm{III}) \cup (\mathrm{IV}).
\end{multline*}
Let us bound the measure of each of these sets.
The second one is easy: by lemma~\ref{l.g}(ii) and our choices
$\gamma=\delta^2/\ell H$ and $N=H/\delta$,
$$
\mu(\mathrm{II}) \leq N \mu(Q^i \minus U^i)
                    < 2 N \gamma \mu(Z^i)
                    < \delta \mu(Z^i).
$$
The other terms are more delicate.

\medskip

The set $X_1 = U^i \minus K^i$ is a $g$-castle whose towers have heights
at least $H$.
Hence its base, which contains $X_1 \minus g(X_1)$,
measures at most $\frac 1H \mu(X_1)$.
By~\eqref{eq.union}, we get
$$
\mu(\mathrm{I}) < \Big( 1+\frac{N}{H} \Big) \mu(X_1) < 2\delta^{-1} \mu(X_1).
$$
By lemma~\ref{l.g}(iii), we have $\mu(X_1)< 2 \kappa \mu(Z^i) = 2\delta^2 \mu(Z^i)$.
So, $\mu(\mathrm{I}) < 4 \delta \mu(Z^i)$.

\medskip

Let $X_3 = \hQ^i \minus Q^i$.
By lemma~\ref{l.hQ minus Q} and~\eqref{e.Z minus ZH}, we have
$\mu(X_3) < \delta^2 \mu(Z_i)$.
Since $f$ and $g$ differ only in $U$, we have
$$
g(X_3) \minus X_3 \subset
[f(X_3) \minus X_3] \cup g(X_3 \cap U)
= (\mathrm{V}) \cup (\mathrm{VI}).
$$
Since $X_3$ is an $f$-castle whose towers have heights of at
least $3H$,
$$
\mu(\mathrm{V}) = \mu(X_3 \minus f(X_3))\leq
\frac{1}{3H} \mu(X_3).
$$
Since $X_3 \cap U \subset \bigsqcup_k (U^k \minus Q^k)$,
lemma~\ref{l.g}(ii) gives $\mu(\mathrm{VI}) \leq 2 \ell \gamma \mu(Z^i)$.
Combining the estimates of $\mu(\mathrm{V})$, $\mu(\mathrm{VI})$,
$\mu(X_3)$ with~\eqref{eq.union} and the definitions of
$N$ and $\gamma$,
\begin{multline*}
\mu(\mathrm{III})
< \mu(X_3) + N \Big( \frac{1}{3H}\mu(X_3) + 2 \ell \gamma \mu(Z^i) \Big)
\\
< \Big(1+ \frac{1}{3\delta}\Big) \mu(X_3) + 2 \delta \mu(Z^i)
< 3 \delta \mu(Z^i).
\end{multline*}

\medskip

We also have
$$
(\mathrm{IV}) = \hQ^i \minus \bigcap_{j=1}^{N-1} g^{-j} (\hQ^i)
\subset \bigcup_{j=1}^{N-1} \big( g^{j-1}(\hQ^i) \minus g^{-j}(\hQ^i) \big)
$$
In particular, $\mu(\mathrm{IV}) \leq (N-1) \mu(\hQ^i \minus g^{-1}(\hQ^i))$.
Notice that $\hQ^i \minus g^{-1}(\hQ^i) \subset \sqcup_k U^k$
(since $\hQ^i$ is $f$-invariant).
Therefore
$$
\hQ^i \minus g^{-1}(\hQ^i) \subset
[\hQ^i \cap \sqcup_{k\neq i} U^k] \cup [U^i \minus g^{-1}(\hQ^i)]
= (\mathrm{VII}) \cup (\mathrm{VIII}).
$$
Combining
$$
(\mathrm{VII}) \subset \bigsqcup_{k\neq i} (U^k \minus \hQ^k) \subset
 \bigsqcup_{k\neq i} (U^k \minus Q^k)
$$
with lemma~\ref{l.g}(ii) we obtain
$\mu(\mathrm{VII}) \leq 2 (\ell-1) \gamma \mu(Z^i)$.
Using that $g(U^i) = f(U^i)$ and $\hQ^i = f(\hQ^i)$, we also get
$$
\mu(\mathrm{VIII}) =
\mu(g(U^i) \minus \hQ^i) =
\mu(U^i \minus \hQ^i) \leq \mu(U^i \minus Q^i) < 2 \gamma \mu(Z^i).
$$
Altogether, $\mu(\hQ^i \minus g^{-1}(\hQ^i))< 2 \ell
\gamma \mu(Z^i)$ and $\mu(\mathrm{IV}) \leq 2 N \ell \gamma
\mu(\hQ_i) < 2\delta\mu(Z_i).$

\medskip

Summing the four parts, we obtain~\eqref{e.hQ minus G}. Now
\begin{align*}
\mu\big(U\cup \Gamma_p^*(f,m) \minus G \big) & \leq \mu (
\Gamma_p^*(f,m) \minus \hQ) + \mu ( U \minus \hQ) + \mu ( \hQ
\minus G) \\ & = \mu(\mathrm{IX}) + \mu(\mathrm{X}) +
\mu(\mathrm{XI}).
\end{align*}
Using \eqref{e.Z minus hQ}, lemma~\ref{l.g}, and \eqref{e.hQ minus
G}, respectively, we get
\begin{align*}
\mu(\mathrm{IX}) &\leq \sum_i \mu(Z^i \minus \hQ^i) < \delta^2 < \delta, \\
\mu(\mathrm{X})  &\leq \mu ( U \minus Q) \leq \sum_i \mu(U^i \minus \hQ^i) <2\gamma <\delta, \\
\mu(\mathrm{XI}) &\leq \sum_i \mu(\hQ^i \minus G) < 10\delta.
\end{align*}
Summing the three parts, we conclude the proof of lemma~\ref{l.big G}.
\end{proof}

\subsection{End of the proof of theorems~\ref{t.vol} and~\ref{t.vol.continuity}} \label{ss.end t.vol}

We give an explicit lower bound for the discontinuity ``jump'' of the
semi-continuous function $\LE_p(\cdot)$.
Denote, for each $p=1, \ldots, d$,
\begin{equation*}\begin{aligned}
J_p(f) & = \int_{\Gamma_p(f,\infty)} \frac{\lambda_p(f,x) - \lambda_{p+1}(f,x)}{2} \ d\mu(x)
\end{aligned}\end{equation*}

\begin{proposition}\label{p.jump}
Given $f \in \Diff$ and $p\in\{1,\ldots,d-1\}$, and given any $\eps_0>0$ and $\delta>0$,
there exists a diffeomorphism $g\in\UU(f,\eps_0)$ such that
$$
\int_M \Lambda_p(g,x) \, d\mu(x) <
\int_M \Lambda_p(f,x) \, d\mu(x) - J_p(f) + \delta.
$$
\end{proposition}

\begin{proof}
Let $f$, $p$, $\eps_0$ and $\delta$ be as in the statement.
Using proposition~\ref{p.global}, we find $m \in \N$ and $g \in \UU(f,\eps_0)$
such that $g=f$ outside $\Gamma_p(f,m)$ and
$$
\int_{\Gamma_p(f,m)} \Lambda_p(g) < \delta +
\int_{\Gamma_p(f,m)} \frac{\Lambda_{p-1}(f) + \Lambda_{p+1}(f)}{2} .
$$
Then
\begin{equation*}\begin{aligned}
\int_M \Lambda_p(g) &=
\int_{\Gamma_p(f,m)} \Lambda_p(g) + \int_{M \minus \Gamma_p(f,m)} \Lambda_p(g) \\
&<
\delta +
\int_{\Gamma_p(f,m)} \frac{\Lambda_{p-1}(f) + \Lambda_{p+1}(f)}{2}
 + \int_{M \minus \Gamma_p(f,m)} \Lambda_p(f).
\end{aligned}\end{equation*}
Since $\Gamma_p(f,\infty)\subset\Gamma_p(f,m)$, and the integrand is non-negative,
\begin{multline*}
\int_{\Gamma_p(f,m)} \left( \Lambda_p(f) - \frac{\Lambda_{p-1}(f) + \Lambda_{p+1}(f)}{2} \right) \ge
\\ \ge
\int_{\Gamma_p(f,\infty)} \left( \Lambda_p(f) - \frac{\Lambda_{p-1}(f) + \Lambda_{p+1}(f)}{2} \right)
=J_p(f).
\end{multline*}
Therefore, the previous inequality implies
\begin{align*}
\int_M \Lambda_p(g)
& <
\delta - J_p(f) + \int_M \Lambda_p(f),
\end{align*}
as we wanted to prove.
\end{proof}

Theorem~\ref{t.vol.continuity} follows easily from
proposition~\ref{p.jump}:

\begin{proof}[Proof of theorem~\ref{t.vol.continuity}]
Let $f \in \Diff$ be a point of continuity of
$\LE_p(\cdot)$ for all $p=1, \ldots, d-1$.
Then $J_p(f) = 0 $ for every $p$.
This means that $\lambda_p(f,x) = \lambda_{p+1}(f,x)$ for
almost every $x\in \Gamma_p(f, \infty)$. Let $x \in M$ be an
Oseledets regular point. If all Lyapunov exponents of $f$ at $x$
vanish, there is nothing to prove. Otherwise, for any $p$ such
that $\lambda_p(f,x) > \lambda_{p+1}(f,x)$, the point
$x \notin \Gamma_p(f, \infty)$ (except for a zero measure set
of $x$) . This means that $x\in\DD_p(f,m)$ for some $m$:
there is a dominated splitting of index $p$, $T_{f^n x}
M = E_n \oplus F_n$, $n\in\Z$ along the orbit of $x$. Clearly,
domination implies that $E_n$ is the sum of the Oseledets
subspaces of $f$, at the point $f^n x$, associated to the Lyapunov
exponents $\lambda_1(f,x)$, \ldots, $\lambda_p(f,x)$, and $F_n$ is
the sum of the spaces associated to the other exponents.
Since this holds whenever $\lambda_p(f,x)$ is bigger than
$\lambda_{p+1}(f,x)$, it proves that the Oseledets splitting
is dominated at $x$.
\end{proof}

Theorem~\ref{t.vol} is an immediate consequence:

\begin{proof}[Proof of theorem~\ref{t.vol}]
The function $f\mapsto \LE_j(f)$
is semi-continuous for every $j=1, \ldots, d-1$, see section~\ref{sss.semicont}.
Hence, there exists a residual subset $\RR$ of $\Diff$ such that
any $f\in\RR$ is a point of continuity for
$f\mapsto (\LE_1(f), \ldots, \LE_{d-1}(f))$.
By theorem~\ref{t.vol.continuity}, every point of continuity
satisfies the conclusion of theorem~\ref{t.vol}.
\end{proof}

\section{Consequences of non-dominance for symplectic maps} \label{s.p.geom sympl}

Here we prove a symplectic analogue of proposition~\ref{p.geom}:

\begin{proposition} \label{p.geom sympl}
Given $f \in \Sympl$, $\eps_0>0$ and $0< \kappa <1$,
if $m \in \N$ is large enough then the following holds:

Let $y \in M$ be a non-periodic point and suppose we
are given a non-trivial splitting $T_y M = E \oplus F$
into two \emph{Lagrangian} spaces
such that
\begin{equation} \label{e.not dominated sympl}
\frac{\| Df^m_y|_F \|}{\mm(Df^m_y|_E)} \geq \frac 12.
\end{equation}
Then there exists a $(\eps_0, \kappa)$-realizable
sequence $\{L_{0}, \ldots, L_{m-1}\}$ at $y$ of length $m$
and there are non-zero vector $v \in E$, $w \in Df^m_y(F)$ such that
$L_{m-1} \cdots L_{0} (v) = w$.
\end{proposition}

\begin{remark}\label{r.solereason}
The hypothesis that $E$ and $F$ are Lagrangian subspaces in
proposition~\ref{p.geom sympl} is the sole reason why
theorem~\ref{t.sympl} is weaker than what is stated in~\cite{IHP}.
\end{remark}

In subsections~\ref{ss.sl1} and~\ref{ss.sn} we
prove two results, namely, lemmas~\ref{l.length 1 sympl} and~\ref{l.nested sympl},
that are used in subsection~\ref{ss.sp} to prove proposition~\ref{p.geom sympl}.
In section~\ref{s.t.sympl} we prove theorems~\ref{t.sympl.continuity}
and~\ref{t.sympl} using proposition~\ref{p.geom sympl}.

\subsection{Symplectic realizable sequences of length $1$} \label{ss.sl1}

First, we recall some elementary facts:
Let $(\cdot,\cdot)$ denote the usual hermitian inner product in $\C^q$.
Up to identification of $\C^q$ with $\R^{2q}$,
the standard inner product in $\R^{2q}$ is $\re (\cdot,\cdot)$ and
the standard symplectic form in $\R^{2q}$ is $\im (\cdot,\cdot)$.
The unitary group $\uni$ is subgroup of $\mathrm{GL}(q,\C)$ formed by the linear maps
that preserve the hermitian product.
Viewing $R\in \uni$ as a map $R: \R^{2q} \to \R^{2q}$, then
$R$ is both symplectic and orthogonal.

If $R: T_x M \to T_x M$ is a $\omega$-preserving linear map,
we shall call $R$ \emph{unitary} if it preserves the inner
product in $T_x M$ induced from the Euclidean inner product
in $\R^{2q}$ by the chart $\varphi_{i(x)}$
(recall subsection~\ref{ss.basic}).

The next lemma constructs realizable sequences of length $1$:

\begin{lemma} \label{l.length 1 sympl}
Given $f\in \Sympl$, $\eps_0>0$, $\kappa>0$, there exists
$\eps>0$ with the following properties: Suppose we are given
a non-periodic point $x \in M$ and an unitary map
$R: T_x M \to T_x M$ with $\| R - I \| < \eps$.
Then $\{Df_x R\}$ is an $(\eps_0,\kappa)$-realizable
sequence of length $1$ at the point~$x$ and
$\{R \, Df_{f^{-1}(x)} \}$ is an $(\eps_0,\kappa)$-realizable
sequence of length $1$ at the point~$f^{-1}x$.
\end{lemma}

We need the following elementary lemma, whose proof is left to the reader:

\begin{lemma} \label{l.Hamiltonian perturb}
Let $H: \R^{2q} \to \R$ be a smooth function such that the corresponding
Hamiltonian flow $\varphi^t:\R^{2q} \to \R^{2q}$ is globally defined for
every $t\in\R$.
Let $\tau: \R \to \R$ be a smooth function and let $\psi$ be a primitive
of $\tau$. Define $\tilde{H} = \psi \circ H$.
Then the Hamiltonian flow $(\tilde{\varphi}^t)$ of $\tilde{H}$ is globally
defined and it is given by $\tilde{\varphi}^t(x) = \varphi^{\tau(H(x))t}(x)$.
\end{lemma}

If $R \in \uni$ then all its eigenvalues belong to the unit circle in $\C$.
Moreover, there exists an orthonormal basis of $\C^q$ formed by eigenvectors of~$R$.
If $J \subset \R$ is an interval, we define $S_J$ as
the set of matrices $R \in \uni$ whose eigenvalues can be written as
$e^{i\theta_1}, \ldots, e^{i\theta_q}$, with all $\theta_k \in J$.
There is $C_0 > 0$, depending only on $q$, such that if
$\eps>0$ and $R \in S_{(-\eps,\eps)}$ then $\| R - I \| < C_0 \eps$.
It is convenient to consider first the case where the arguments of the
eigenvalues of $R$ have all the same sign:

\begin{lemma} \label{l.compact type}
Given $\eps_0>0$ and $0< \sigma <1$,
there exists $\eps>0$ with the following properties:
Given $R \in S_{(-\eps,0)}\cup S_{(0,\eps)}$,
there exists a bounded open set $U \subset \R^{2q}$
such that $\sigma U \subset U$, and there exists a
$C^1$ symplectomorphism $h:\R^{2q} \to \R^{2q}$ such that
\begin{itemize}
\item[(i)]  $h(z)=z$ for every $z\notin U$ and
            $h(z)=R(z)$ for every $z\in \sigma U$;
\item[(ii)]
            $\| Dh_z - I \| < \eps_0$ for all $z\in\R^{2q}$.
\end{itemize}
\end{lemma}

\begin{proof}[Proof of lemma \ref{l.compact type}]
Let $\eps_0$ and $\sigma$ be given. Let $\eps>0$ be a small number,
to be specified later. Take $R \in S_{(0,\eps)}$: the other possibility
is tackled in a similar way.
Let $\{v_1, \ldots, v_q \}$ be an orthonormal basis of eigenvectors of $R$,
with associated eigenvalues $e^{i\theta_1},\ldots,e^{i\theta_q}$,
and all $0 < \theta_k < \eps$.
Up to replacing $R$ with $SRS^{-1}$, for some $S\in \uni$, we may assume
that the basis $\{v_1, \ldots, v_q \}$ coincides with the standard basis of $\C^q$.
Therefore $R$ assumes the form
$$
R(z_1, \ldots, z_q) = (e^{i\theta_1} z_1, \ldots, e^{i\theta_q} z_q)
$$
Let $H : \C^q \to \R $ be given by $H(z) = \meio \sum_k \theta_k |z_k|^2$.
Then $R$ is the time $1$ map of the Hamiltonian flow of $H$.
Besides, there is $C_1=C_1(q)$ such that
\begin{equation}\label{eq.zDz}
\|z\| \, \|DH_z\| \leq C_1 H(z) \quad
\text{for all $z\in\C^q$.}
\end{equation}

Let $\tau:\R \to \R$ be a smooth function such that
$\tau(s)=1$ for $s\leq \sigma^2$, $\tau(t)=0$ for $t\geq 1$,
and $0 \leq -\tau'(t)\leq 2/(1-\sigma^2)$ for all $t$.
Let $\psi(s) = \int_0^s \tau(u)\, du$ and let
$\tilde{H} = \psi \circ H$.
By lemma~\ref{l.Hamiltonian perturb}, the time $1$ map $h$ of the Hamiltonian
flow of $\tilde{H}$ is
$$
h(z) = (e^{i\theta_1 \tau(H(z))} z_1 , \ldots, e^{i\theta_k \tau(H(z))} z_k)
$$
Then $h(z) = R(z)$ if $H(z) \leq \sigma^2$ and $h(z) = z$ if $H(z) \geq 1$.
Moreover, a direct calculation and \eqref{eq.zDz} give
$$
\|Dh_z - I \| \leq C_2 \eps (1-\sigma^2)^{-1} + \eps
\quad \text{for every $z\in\C^q$ with } H(z)\leq 1,
$$
where $C_2=C_2(q)$.
Take $\eps = \eps(\eps_0, \sigma)$ such that the right hand side is
less than $\eps_0$. Since $H$ is definite positive, the set
$U = \{z \in \C^q ; \; H(z) < 1\}$ is bounded.
\end{proof}

\begin{remark}\label{r.scaling trick}
We may assume that the set $U$ in lemma~\ref{l.compact type}
has arbitrarily small diameter.
Indeed, if $a>0$ then we may replace $U$ with $\widetilde{U} = aU$ and
$h$ with $\widetilde{h}(z) = a h(a^{-1} z)$.
Notice $D\widetilde{h}_z = Dh_{a^{-1}z}$, so $\widetilde{h}$
is a symplectomorphism and satisfies property (ii) of the lemma.
\end{remark}

\begin{lemma} \label{l.any type}
Given $\eps_0>0$ and $0< \sigma <1$,
there exists $\eps>0$ with the following properties:
Given $R \in S_{(-\eps,\eps)}$, there exists a bounded open set
$U \subset \R^{2q}$ such that $\sigma U \subset U$, a measurable
set $K \subset U$ with $\vol(U \minus K) < 3(1-\sigma^d) \vol(U)$,
and a $C^1$ symplectomorphism $h:\R^{2q} \to \R^{2q}$ such that
\begin{itemize}
\item[(i)]  $h(z) = z$ for every $z\notin U$ and
            $Dh_z = R$ for every $z\in K$;
\item[(ii)]
            $\| Dh_z - I \| < \eps_0$ for all $z\in\R^{2q}$.
\end{itemize}
\end{lemma}

\begin{proof}
Any $R\in S_{(-\eps, \eps)}$ can be written as a product $R = R_+ \, R_-$,
with $R_+ \in S_{(0,\eps)}$ and $R_- \in S_{(-\eps,0)}$, in fact we may
take $R_+$ and $R_-$ with the same eigenbasis as $R$.
Applying lemma~\ref{l.compact type} to $R_{\pm}$, with $\eps_0$
replaced by $\eps_0/2$, we obtain sets $U_\pm$ and
symplectomorphisms $h_\pm$. Let $U = U_+$.
Consider the family $\mathcal{F}$ of all sets of the form
$a U_- + b$, with $a>0$ and $b\in \R^{2q}$, that are contained in $U$.
This is a Vitali covering of $U$, so we may find a finite number of
disjoint sets $U_-^i = a_i U_- + b_i  \in \mathcal{F}$ that cover $U$
except for a set of volume $(1-\sigma^d) \vol(U)$.
Using lemma~\ref{l.compact type} and remark~\ref{r.scaling trick},
for each $i$ we find a symplectomorphism $h^i_-$ such that
$h^i_- = \id$ outside $U^i_-$\, and
$D(h^i_-)_z = R_-$ for $z \in K_i = a_i \sigma U_- + b_i$,
and $D(h^i_-)_z$ is uniformly close to $I$.
Let $K = (\sigma U) \cap \sqcup_i K^i$.
Define $h = h_+ \circ h^i_-$ inside each $U^i_-$, and $h = h_+$ outside.
Then $K$ and $h$ have the desired properties.
\end{proof}

\begin{proof}[Proof of lemma~\ref{l.length 1 sympl}]
Given $\eps_0$ and $\kappa$, choose $\sigma$ close to $1$ so that
$3(1-\sigma^d)<\kappa$.
Remark~\ref{r.scaling trick} also applies to lemma~\ref{l.any type}:
the set $U$ may be taken with arbitrarily small diameter.
Using lemma~\ref{l.simplifying}, we conclude that the sequences
$\{Df_x R\}$ and $\{R \, Df_{f^{-1}(x)} \}$ are $(\eps_0,\kappa)$-realizable
as stated.
\end{proof}

\subsection{Symplectic nested rotations} \label{ss.sn}

In this subsection we prove an analogue of lemma~\ref{l.nested}
for symplectic maps:

\begin{lemma} \label{l.nested sympl}
Given $f\in \Sympl$, $\eps_0>0$, $\kappa>0$, $E>1$, and $0 < \gamma
\leq \pi/2$, there exists $\beta>0$ with the following properties:
Suppose we are given a non-periodic point $x \in M$, an iterate $n \in \N$,
and a two-dimensional symplectic subspace $Y_0 \subset T_x M$ such
that:
\begin{itemize}
\item $\| Df^j |_{Y_0} \| / \mm(Df^j |_{Y_0}) \leq E^2$ for every $j=1,\ldots,n$;
\item $\ang(X_j, Y_j) \geq \gamma$ for each $j=0,\ldots,n-1$ where $X_0 = Y_0^\omega$,
      $X_j= Df^j_x (X_0)$, and $Y_j= Df^j_x (Y_0)$.
\end{itemize}
Let $\theta_0,\ldots,\theta_{n-1} \in [-\beta,\beta]$
and let $S_0,\ldots,S_{n-1} : Y_0 \to Y_0$ be the rotations
of the plane $Y_0$ by angles $\theta_0,\ldots,\theta_{n-1}$.
Let linear maps
$$
T_x M \xrightarrow{L_0} T_{fx} M \xrightarrow{L_1} \cdots
\xrightarrow{L_{n-1}}T_{f^n(x)}M
$$
be defined by
$L_j (v) =  Df_{f^j(x)} (v)$ for $v\in X_j$ and
$L_j(w) = (Df_y^{j+1}) \cdot S_j \cdot (Df_y^j)^{-1}(w)$ for $w \in Y_j$.
Then $\{L_0,\ldots ,L_{n-1}\}$ is an $(\eps_0,\kappa)$-realizable
sequence of length $n$ at the point~$x$.
\end{lemma}

\smallskip

We begin by proving a perturbation lemma that corresponds to lemma~\ref{l.rotcil}:

\begin{lemma} \label{l.rotcil sympl}
Given $\eps_0>0$ and $0<\sigma<1$, there is $\eps>0$ with the
following properties: Suppose we are given:
a splitting $\R^{2q} = X \oplus Y$ with $\dim Y =2$, $X^\omega = Y$ and
$X \perp Y$, an ellipsoid $\AA \subset X$ centered at the origin, and
a unitary map $R \in \uni$ with $R|_X = I$ and $\| R - I \| < \eps$.

Then there exists $\tau>1$ such that the following holds.
Let $\BB$ be the unit ball in $Y$.
For $a$, $b>0$ consider the cylinder $\CC=\CC_{a,b}=a\AA\oplus b\BB$.
If $a > \tau b$ and $\diam \CC < \eps_0$ then there is a $C^1$
symplectomorphism $h:\R^{2q} \to \R^{2q}$ satisfying:
\begin{itemize}
\item[(i)] $h(z)=z$ for every $z\notin \CC$ and
           $h(z)=R(z)$ for every $z\in \sigma\CC$;
\item[(ii)] $\| h(z) - z \| < \eps_0$ and
            $\| Dh_z - I \| < \eps_0$ for all $z\in\R^{2q}$.
\end{itemize}
\end{lemma}


\begin{remark} \label{r.ode}
If $H:\R^{2q} \to \R$ is a smooth function with bounded $\|DH\|$ and $\|D^2 H\|$,
then the associated Hamiltonian flow $\varphi^t :\R^{2q} \to \R^{2q}$
is defined for every time $t\in\R$, and
$$
\| \varphi^t(z) - z \|   \leq |t| \sup \|DH\|, \qquad
\| (D\varphi^t)_z - I \| \leq \exp(|t| \sup \|D^2 H\|)-1.
$$
for every $z \in \R^{2q}$ and $t\in\R$.
\end{remark}

\begin{proof}[Proof of lemma \ref{l.rotcil sympl}]
Given $\eps_0$ and $\sigma$, let
$$
K =
10 (1-\sigma)^{-2} + 20 \sigma^{-1} (1-\sigma)^{-1}
+ 30 (1-\sigma)^{-1} + 3.
$$
Fix $\bar t>0$ such that $e^{\bar t K}-1< \eps_0$\,, and let
$\eps>0$ be such that $\eps < \sqrt{2} \, \sin \bar t$.
Let $X$, $Y$, $\AA$, $\BB$, and $R$ be as in the statement.
Let $A: X \to X$ be a linear map such that $A(\AA)$ is the
unit ball in $X$. We define $\tau = \|A\|$.

Let $H: \R^{2q} \to \R$ be defined by $H(x,y) = H(y) = \meio \| y\|^2$,
where $(x,y)$ are coordinates with respect to the splitting $X \oplus Y$.
The Hamiltonian flow of $H$ is a linear flow $(R_t)_t$,
where $R_t$ is a rotation of angle $t$ in the plane $Y$, with axis $X$.
In particular, $\| R_t - I\| = \sqrt{2} \, |\sin t|$
and there exists $t_0$ with $|t_0| < \bar t$ such that $R_{t_0} = R$.

Take numbers $a$, $b>0$ such that $a/b > \tau$ and the cylinder
$\CC = a\AA\oplus b\BB$ has diameter less than $\eps_0$.
We are going to construct another Hamiltonian  $\widetilde{H}$
which is equal to $H$ inside $\sigma \CC$ and constant
outside $\CC$. The symplectomorphism $h$ will be defined as the
time $t_0$ of the Hamiltonian flow associated to $\widetilde{H}$.

For this we need a few auxiliary functions.
Let $\zeta: \R \to [0,1]$ be a smooth function such that:
\begin{itemize}
\item $\zeta(t)=1$ for $t\leq \sigma$ and $\zeta(t)=0$ for $t\geq 1$;
\item $|\zeta'(t)| \leq 10/(1-\sigma)$ and $|\zeta''(t)| \leq 10/(1-\sigma)^2$.
\end{itemize}
Let $\hat{\psi} : X \to [0,1]$ be defined by
$\hat{\psi}(x) = \zeta (a^{-1}\|x\|)$, and
    $\psi : X \to [0,1]$ be defined by
$\psi = \hat{\psi} \circ A$. It is clear that
\begin{equation}\label{e.psi}
\psi(x)=1 \text{ for } x \in \sigma a\AA
\quad\text{and}\quad
\psi(x)=0 \text{ for } x \notin a\AA.
\end{equation}
We estimate the derivatives:
\begin{align*}
D\hat{\psi}_x (v) &= a^{-1} \zeta'(a^{-1} \|x\|) \frac{\langle x,v \rangle}{\|x\|},
\\
\begin{split}
D^2\hat{\psi}_x (v,w) &=
a^{-2} \zeta''(a^{-1} \|x\|) \frac{\langle x,v \rangle \langle x,w \rangle}{\|x\|^2} +
a^{-1} \frac{\zeta'(a^{-1} \|x\|)}{\|x\|}
\left(\langle v,w \rangle - \frac{\langle x,v \rangle \langle x,w \rangle}{\|x\|^2} \right).
\end{split}
\end{align*}
Since $|{\zeta'(a^{-1} \|x\|)|}/{\|x\|} \leq 10 a (1-\sigma)^{-1}\sigma^{-1}$,
we get the bounds
\begin{align*}
\|D \hat{\psi} \| &\leq 10 (1-\sigma)^{-1} a^{-1}
\quad\text{and}\quad
\|D^2 \hat{\psi} \| &\leq \left[10 (1-\sigma)^{-2} + 20 \sigma^{-1} (1-\sigma)^{-1} \right] a^{-2}.
\end{align*}
Moreover,
$$
\|D\psi_x\| \leq \|A\| \, \|D \hat{\psi}_x \|
\quad\text{and}\quad
\|D^2\psi_x\| \leq \|A\|^2 \|D^2 \hat{\psi}_x \|
$$

Now define $\rho : \R \to \R$ by $\rho(t) = \int_0^t \zeta$ and then
let $\phi: Y \to \R$ be given by
$\phi(y) = \meio b^2 \rho(b^{-1} \|y\|)^2$.
Then
\begin{equation}\label{e.phi}
\phi(y) = H(y) \text{ for } y \in \sigma b\BB
\quad\text{and}\quad
\phi(y) = c \text{ for } y \notin b\BB,
\end{equation}
where $0< c < \meio b^2$ is a constant.
The first and second derivatives of $\phi$ are:
\begin{align*}
D\phi_y (v) &= b \rho(b^{-1} \|y\|) \rho'(b^{-1} \|y\|) \, \frac{\langle y,v \rangle}{\|y\|},
\\
\begin{split}
D^2\phi_y (v,w) &=
\left[ \rho'(b^{-1} \|y\|)^2 + \rho(b^{-1} \|y\|) \rho''(b^{-1} \|y\|) \right]
\frac{\langle y,v \rangle}{\|y\|} +
\\ &\qquad\qquad
b \rho(b^{-1} \|y\|) \, \frac{\rho'(b^{-1} \|y\|)}{\|y\|}
\left(\langle v,w \rangle - \frac{\langle y,v \rangle \langle y,w \rangle}{\|y\|^2} \right).
\end{split}
\end{align*}
Since $|\rho| \leq 1$, $|\rho'| \leq 1$, $|\rho''| \leq 10(1-\sigma)^{-1}$,
and $|\rho'(b^{-1} \|y\|)|/{\|y\|} \leq b^{-1}$, we have
$$
\|D \phi \| \leq b
\quad \text{and} \quad
\|D^2 \phi \| \leq 3 + 10(1-\sigma)^{-1}.
$$

Define $\widetilde{H} : \R^{2q} \to \R$ by
$\widetilde{H} (x,y) = c - \psi(x) (c - \phi(y))$.
Then, by~\eqref{e.psi} and~\eqref{e.phi},
\begin{equation}\label{eq.dentrofora}
\begin{aligned}
x \in \sigma a \AA \quad\text{and}\quad y \in \sigma b \BB
& \Rightarrow \widetilde{H}(x,y) = H(y),
\\
x \notin a \AA \quad\text{or}\quad y \notin b \BB
&\Rightarrow \widetilde{H}(x,y) = c.
\end{aligned}
\end{equation}
The derivatives of $\widetilde{H}$ are (write $v=v_x+v_y\in X\oplus Y$ and
analogously for $w$)
\begin{align*}
D\widetilde{H}_{(x,y)} (v) &= -(c- \phi(y)) D\psi_x (v_x) + \psi(x) D\phi_y (v_y),
\\
\begin{split}
D^2\widetilde{H}_{(x,y)} (v,w) & =
-(c- \phi(y)) D^2\psi_x (v_x,w_x) +
D\psi_x(v_x) D\phi_y(w_y) +
\\ & \qquad\qquad
D\psi_x(w_x) D\phi_y(v_y)+
\psi(x) D^2\phi_y (v_y,w_y).
\end{split}
\end{align*}
Using the previous bounds we obtain
\begin{multline*}
\| D^2\widetilde{H} \| \leq
\left[10 (1-\sigma)^{-2} + 20 \sigma^{-1} (1-\sigma)^{-1} \right] \|A\|^2  (b/a)^2
+\\
20 (1-\sigma)^{-1} \|A\| (b/a)
+ 3 + 10(1-\sigma)^{-1}.
\end{multline*}
Since $a/b > \|A\|$, we conclude that $\| D^2\widetilde{H} \| \leq K$.

Take $h : \R^{2q} \to \R^{2q}$ to be the time $t_0$ map of the Hamiltonian
flow associated to $\widetilde{H}$.
Property~(i) in the lemma follows from~\eqref{eq.dentrofora}.
Since $\diam \CC < \eps_0$, we have $\|h(z)-z\|< \eps_0$ for all $z$.
And, by remark~\ref{r.ode},
$\|Dh_z - I\| \leq e^{t_0 K}-1 \leq e^{\bar t K}-1 < \eps_0$,
proving~(ii) and the lemma.
\end{proof}

An ellipse $\BB$ contained in a $2$-dimensional symplectic subspace
$Y \subset \R^{2q}$ and centered at the origin
has \emph{eccentricity $E$} if it is the image of the unit ball under
a linear transformation $\hat B:Y \to Y$ with
$\| \hat B \| / \mm(\hat B) = E^2$.
If a map $\hR: Y \to Y$ preserves the ellipse $\BB$, then
$\hat B^{-1} \hR \hat B$ is a rotation of the plane $Y$ of some angle $\theta$.
In this case we say that \emph{$\hR$ rotates the ellipse $\BB$
through angle $\theta$}.

\smallskip

The following statement is a more flexible version of
lemma~\ref{l.rotcil sympl}.
In fact, it follows from lemma~\ref{l.rotcil sympl} just by a change of the inner product.

\begin{lemma} \label{l.rotcil 2}
Given $\eps_0>0$, $0<\sigma<1$, $\gamma>0$ and $E > 1$,
there is $\beta>0$ with the following properties:
Suppose we are given:
\begin{itemize}
\item a splitting $\R^{2q} = X \oplus Y$ with $\dim Y = 2$,
      $X^\omega = Y$ and $\ang(X,Y) \geq \gamma$;
\item an ellipsoid $\AA \subset X$ centered at the origin;
\item an ellipse $\BB \subset Y$ centered at the origin and
      with eccentricity at most~$E$;
\item a map $\hR: Y \to Y$ that rotates $\BB$ through angle $\theta$,
      with $|\theta| < \beta$.
\end{itemize}
Then there exists $\tau>1$ such that the following holds.
Let $R: \R^{2q} \to \R^{2q}$ be the linear map defined by
$R(v) = v$ if $v\in X$ and $R(w) = \hR (w)$ if $w\in Y$.
For $a, b>0$ consider the cylinder $\CC=\CC_{a,b}=a\AA\oplus b\BB$.
If $a > \tau b$ and $\diam \CC < \eps_0$ then there is a $C^1$
symplectomorphism $h:\R^{2q} \to \R^{2q}$ satisfying:
\begin{itemize}
\item[(i)] $h(z)=z$ for every $z\notin \CC$ and
           $h(z)=R(z)$ for every $z\in \sigma\CC$;
\item[(ii)] $\| h(z) - z \| < \eps_0$ and\ \
            $\| Dh_z - I \| < \eps_0$ for all $z\in\R^{2q}$.
\end{itemize}
\end{lemma}

Now lemma~\ref{l.nested sympl} is proved in the same way as we
proved lemma~\ref{l.nested}, using lemmas~\ref{l.rotcil 2}
and~\ref{l.nlinear} instead. The argument is even a bit
simpler since no truncation (like in lemma~\ref{l.trunca}) is
necessary, as we assume that the angles $\ang(X_j,Y_j)$ are
bounded from zero. The details are left to the reader.

\subsection{Proof of proposition~\ref{p.geom sympl}} \label{ss.sp}

We use the following lemma, which will also be needed in
section~\ref{s.t.cocycle}:

\begin{lemma} \label{l.transitive}
Let $G \subset \gldr$ be a closed group which acts transitively
in $\rp$.
Then for every $\eps_1>0$ there exists $\alpha>0$ such that if
$v_1$, $v_2 \in \R^d$ satisfy $\ang(v_1,v_2) < \alpha$ then there
exists $R\in G$ such that $\|R-I\|< \eps_1$ and $R(\R v_1) = \R v_2$.
\end{lemma}

\begin{proof}
For $\delta>0$, let $U_\delta = \{ R; \; R \in G,\ \|R-I\| < \delta \}$.
Given $\eps>0$, fix $\delta>0$ such that
if $R_1$, $R_2 \in U_\delta$ then $R_2 R_1^{-1} \in U_{\eps_1}$.
The hypothesis on the group implies that for any $w\in \rp$,
the map $G \to \rp$ given by $A \mapsto A(w)$ is open
(this follows from~\cite[Theorem II.3.2]{Helgason}).
Therefore, for any $\delta>0$, the set
$U_\delta(w) = \{ Rw; \; R \in U_\delta \}$
is an open neighborhood of $w$.
Cover $\rp$ by some finite union $U_\delta(w_1) \cup \cdots \cup U_\delta(w_k)$.
Now take two directions $v_1$, $v_2 \in \rp$
sufficiently close. Then both belong to some $U_\delta(w_i)$,
and so there are $R_1$, $R_2 \in U_\delta$ such that
$v_1 = R_1 w_i$ and $v_2 = R_2 w_i$.
Therefore $R = R_2 R_1^{-1}$ belongs to $U_{\eps_1}$ and
$Rv_1 = v_2$.
\end{proof}

\begin{proof}
Let $f$, $\eps_0$\,, $\kappa$ be given. Fix $0< \kappa' < \meio \kappa$.
Let $\eps>0$, depending on $f$, $\eps_0$\,, $\kappa'$,
be given by lemma~\ref{l.length 1 sympl}.
Let $\alpha>0$, depending on $\eps_1=\eps$ and $G = \uni$,
be given by lemma~\ref{l.transitive}.
Take $K$ satisfying $K \geq (\sin \alpha)^{-2}$
and $K \geq \max_x \|Df_x\|/\mm(Df_x)$.
Let $E>1$ and $\gamma >0$ be given by
$$
E^2 = 8 C_\omega^4 K (\sin \alpha)^{-4} \quad \text{and} \quad
\sin \gamma = \meio C_\omega^{-14} K^{-2} \sin^9 \alpha,
$$
where $C_\omega$ is as in~\eqref{e.def C omega}.
Let $\beta >0$ be given by lemma~\ref{l.nested sympl}.
Finally, let $m \geq 2 \pi/\beta$.
The proof is divided into three cases.

\paragraph{First case:}
Suppose that there exists $\ell \in \{0,1,\ldots,m\}$ such that
\begin{equation} \label{e.I sympl}
\ang(E_\ell, F_\ell) < \alpha.
\end{equation}
Fix $\ell$ as above and take unit vectors $\xi \in E_\ell$, $\eta \in F_\ell$
such that $\ang (\xi,\eta) < \alpha$.
By lemma~\ref{l.transitive}, there exists a unitary transformation
$R:T_{f^\ell(y)} M \to T_{f^\ell(y)} M$ such that $\|R - I \| < \eps$
and $R(\xi) = \eta$.
By lemma~\ref{l.length 1 sympl}, the sequences $\{Df_{f^\ell(x)} \, R\}$ and
$\{R \, Df_{f^{\ell-1}(x)} \}$ are $(\kappa', \eps_0)$-realizable.
Define $\{ L_0, \ldots L_{m-1} \}$ as
$$
\{ Df_y, \ldots, Df_{f^{\ell-1}(y)}, Df_{f^\ell(y)} \, R,
   Df_{f^{\ell+1}(y)}, \ldots, Df_{f^{m-1}(y)} \}
$$
if $\ell<m$ and as $\{ Df_y, \ldots, Df_{f^{m-2}(y)}, R \, Df_{f^{m-1}(y)} \}$
if $\ell=m$.
In either case, this is a $(\kappa, \eps_0)$-realizable sequence of length $m$
at $y$, whose product $L_{m-1} \cdots L_0$ sends the
direction $\R Df^{-\ell} (\xi) \subset E_0$ to the direction
$\R Df^{m-\ell} (\eta) \subset F_m$.

\paragraph{Second case:}
Assume that there exist $k, \ell \in \{0,\ldots,m\}$ with $k<\ell$ and
\begin{equation} \label{e.II sympl}
\frac{\|Df^{\ell-k}_{f^k(y)} |_{F_k} \|}{\mm (Df^{\ell-k}_{f^k(y)} |_{E_k})} > K.
\end{equation}
The proof of this case is easily adapted from
the second case in the proof of proposition~\ref{p.geom}.
We leave the details to the reader.

\paragraph{Third case:}
We suppose that we are not in the previous cases, that is,
\begin{equation} \label{e.not I sympl}
\text{for every $j \in \{0,1,\ldots,m\}$,}\quad
\ang(E_j, F_j) \geq \alpha.
\end{equation}
and
\begin{equation} \label{e.not II sympl}
\text{for every $i,j \in \{0,\ldots,m\}$ with $i<j$,} \quad
\frac{\|Df^{j-i}_{f^i(y)} |_{F_i} \|}{\mm (Df^{j-i}_{f^i(y)} |_{E_i})} \leq K.
\end{equation}
By~\eqref{e.not I sympl} and lemma~\ref{l.sympl 1}, we have,
for all $i,j \in \{0,\ldots,m\}$ with $i<j$,
\begin{equation}\label{e.opposite}
C_\omega^{-2} \sin \alpha \leq
\mm(Df^{j-i}|_{E_i}) \, \| Df^{j-i}|_{F_i} \| \leq
C_\omega^2 (\sin \alpha)^{-1}.
\end{equation}
This, together with~\eqref{e.not II sympl}, gives
\begin{align}
\mm(Df^{j-i}|_{E_i})  &\geq C_\omega^{-1} K^{-1/2} (\sin \alpha)^{1/2}, \label{e.pure1E}
\\
\| Df^{j-i}|_{F_i} \| &\leq C_\omega K^{1/2} (\sin \alpha)^{-1/2}. \label{e.pure1F}
\end{align}
Also, by~\eqref{e.opposite} and the main assumption~\eqref{e.not dominated sympl},
\begin{align}
\mm(Df^m|_{E_0})  &\leq 2^{1/2} C_\omega (\sin \alpha)^{-1/2}, \label{e.pure2E}
\\
\| Df^m|_{F_0} \| &\geq 2^{-1/2} C_\omega^{-1} (\sin \alpha)^{1/2}. \label{e.pure2F}
\end{align}

Let $v_0 \in E_0$ be such that $\|v_0\| = 1$ and
$\|Df_y^m (v_0) \| = \mm(Df_y^m|_{E_0})$.
Using lemma~\ref{l.sympl 1}.1, take $w_0 \in F_0$ with
$\|w_0\|=1$ such that $|\omega(v_0,w_0)| \geq C_\omega^{-1} \sin \alpha$.
Let $G_0 = E_0 \cap w_0^\omega$ and $H_0 = F_0 \cap v_0^\omega$.
(By $v^\omega$ we mean $(\R v)^\omega$.)
Let $X_0 = G_0 \oplus H_0$ and $Y_0 = \R v_0 \oplus \R w_0$.
Then $X_0 = Y_0^\omega$.
Let, for $j=1,\ldots,m$,
\begin{alignat*}{3}
v_j &= Df^j(v_0) / \|Df^j(v_0) \|, &\qquad G_j &= Df^j (G_0),  &\qquad X_j &= Df^j (X_0), \\
w_j &= Df^j(w_0) / \|Df^j(w_0) \|, &\qquad H_j &= Df^j (H_0),  &\qquad Y_j &= Df^j (Y_0)
\end{alignat*}
(all the derivatives are at $y$). By~\eqref{e.C omega},
$$
C_\omega^{-1} \sin \alpha \leq |\omega(v_0,w_0)| = |\omega(Df^m v_0, Df^m w_0)|
\leq C_\omega \|Df^m v_0\| \, \|Df^m w_0\|.
$$
Thus $\|Df^m w_0\| \geq C_\omega^{-2} \sin \alpha \cdot \mm(Df^m|_{E_0})^{-1}$
and, using~\eqref{e.opposite},
\begin{equation} \label{e.w0 is good}
\|Df^m w_0\| \geq C_\omega^{-4} \sin^2 \alpha \cdot \|Df^m |_{F_0} \|
\end{equation}
that is, $w_0$ is ``almost'' the most expanded vector by $Df^m$ in $F_0$.
By~\eqref{e.not dominated sympl} and~\eqref{e.w0 is good},
$$
\frac{\|Df^m w_0\|}{\|Df^m v_0\|}
\geq C_\omega^{-4} \sin^2 \alpha \frac{\|Df^m _{F_0}\|}{\mm(Df_y^m|_{E_0})}
\geq \meio C_\omega^{-4} \sin^2 \alpha.
$$
This and~\eqref{e.not II sympl} imply that, for each $j=1,\ldots,m$,
$$
K \geq
\frac{\|Df^j w_0\|}{\|Df^j v_0\|} \geq
\frac{\|Df^m w_0\|/ \|Df^{m-j}|_{F_j} \|}{\|Df^m v_0\| / \mm(Df^{m-j}|_{E_j})} \geq
\meio C_\omega^{-4} K^{-1} \sin^2 \alpha.
$$
Therefore, using~\eqref{e.not I sympl} and lemma~\ref{l.4},
\begin{equation} \label{e.eccentricity s}
\frac{\|Df^j|_{Y_0}\|}{\mm(Df^j|_{Y_0})} \leq
8 C_\omega^4 K (\sin \alpha)^{-4} = E^2.
\end{equation}

\smallskip

We now deduce some angle estimates. First, we claim that
\begin{equation}\label{e.first angle}
\sin \ang(v_0, G_0)\geq C_\omega^{-2} \sin \alpha
\quad\text{and}\quad
\sin \ang(w_0, H_0) \geq C_\omega^{-2} \sin \alpha.
\end{equation}
Indeed, write $v_0 = u + u'$ with $u'\in G_0$ and $u \perp G_0$.
Since $G_0$ is skew-orthogonal to $w_0$,
$$
C_\omega^{-1} \sin \alpha \leq |\omega(v_0,w_0)| = |\omega(u, w_0)|
\leq C_\omega \|u\|.
$$
That is, $\sin\ang(v_0,G_0) = \|u\| \geq C_\omega^{-2} \sin \alpha$.
Analogously we prove the other inequality in~\eqref{e.first angle}.
Next, we estimate $\sin \ang(v_j, G_j)$ and $\sin \ang(w_j, H_j)$
for $j=1,\ldots,m$. For this we use relation~\eqref{e.useful}
from subsection \ref{ss.angle tools}, which gives:
\begin{align}
\sin \ang(v_j, G_j) &\geq
\frac{\mm(Df^j|_{E_0})}{\|Df^j v_0 \|} \sin \ang(v_0, G_0), \label{e.angle v aux}
\\
\sin \ang(w_j, H_j) &\geq
\frac{\|Df^j w_0 \|}{\|Df^j |_{F_0}\|} \sin \ang(w_0, H_0). \label{e.angle w aux}
\end{align}
By~\eqref{e.pure1E} and~\eqref{e.pure2E}, $$
\|Df^j v_0 \|
= \frac{\|Df^m v_0 \|}{\|Df^{m-j} v_j \|}
\leq \frac{\mm(Df^m|_{E_0})}{\mm(Df^{m-j} |_{E_0})}
\leq 2^{1/2} C_\omega^2 K^{1/2} (\sin \alpha)^{-1}.
$$
for each $j=1,\ldots,m$. So, using~\eqref{e.pure1E} again,
$$
\frac{\|Df^j v_0 \|}{\mm(Df^j|_{E_0})} \leq
2^{1/2} C_\omega^3 K (\sin \alpha)^{-3/2}.
$$
This, together with~\eqref{e.first angle} and \eqref{e.angle v aux}, gives
\begin{equation} \label{e.angle v}
\sin \ang(v_j, G_j) \geq 2^{-1/2} C_\omega^{-5} K^{-1} (\sin \alpha)^{5/2}.
\end{equation}
Similarly, by~\eqref{e.w0 is good}, \eqref{e.pure1F}, and ~\eqref{e.pure2F},
$$
\|Df^j w_0 \|
= \frac{\|Df^m w_0 \|}{\|Df^{m-j} w_j \|}
\geq C_\omega^{-4} \sin^2 \alpha \, \frac{\|Df^m |_{F_0} \|}{\|Df^{m-j} |_{F_0} \|}
\geq 2^{-1/2} C_\omega^{-6} K^{-1/2} \sin^3 \alpha.
$$
Then, using~\eqref{e.pure1F} again,
$$
\frac{\|Df^j w_0 \|}{\|Df^j |_{F_0}\|} \geq
2^{-1/2} C_\omega^{-7} K^{-1} (\sin \alpha)^{7/2}.
$$
By \eqref{e.first angle} and \eqref{e.angle w aux},
\begin{equation} \label{e.angle w}
\sin \ang(w_j, H_j) \geq 2^{-1/2} C_\omega^{-9} K^{-1} (\sin \alpha)^{9/2}.
\end{equation}
Now we use lemma~\ref{l.abc} three times:
\begin{equation*}\begin{aligned}
\sin \ang(Y_j, X_j) &
\geq \sin \ang(v_j, X_j) \, \sin \ang(w_j, \R v_j \oplus X_j )\\
& \geq \sin \ang(v_j, G_j) \, \sin \ang(E_j, H_j) \, \sin \ang(w_j, \R v_j \oplus X_j )\\
& \geq \sin \ang(v_j, G_j) \, \sin \ang(w_j, H_j) \, \sin \ang(E_j, H_j)^2
\end{aligned}\end{equation*}
So, using~\eqref{e.angle v}, \eqref{e.angle w}, and $\ang(E_j, H_j)\geq \alpha$,
we obtain
\begin{equation}\label{e.last angle}
\sin \ang(X_j, Y_j) \geq \meio C_\omega^{-14} K^{-2} \sin^9 \alpha = \sin \gamma.
\end{equation}

\smallskip

Relations~\eqref{e.eccentricity s} and \eqref{e.last angle} permit us
to apply lemma~\ref{l.nested sympl}.
Since $m\beta \geq 2\pi$, it is possible to
choose numbers $\theta_0,\ldots,\theta_{m-1}$ such that $0\leq \theta_j \leq \beta$
and $\sum \theta_j = \ang(v_0 , w_0)$.
Let $S_j$ and $L_j$ be as in lemma~\ref{l.nested sympl}.
We have
$L_{m-1} \cdots L_0 |_{Y_0} = (Df^m|_{Y_0}) S_{m-1} \cdots S_0$,
so $L_{m-1} \cdots L_0 (\R v_0) = \R w_m$.
This completes the proof of proposition~\ref{p.geom sympl}.
\end{proof}

\section{Proof of theorems \ref{t.sympl.continuity} and \ref{t.sympl}} \label{s.t.sympl}

Given $f \in \Diff$ and $m \in \N$,
let $\DD(f,m)$ be the (closed) set of points $x$ such that there is
a $m$-dominated splitting of index $q=d/2$ along the orbit of $x$.
Let $\Gamma (f, m) = M \minus \DD(f,m)$ and
let $\Gamma^* (f, m)$ be the set of points $x \in \Gamma(f,m)$
which are regular, not periodic and
satisfy $\lambda_q(f,x) > 0$.
Let also $\Gamma (f, \infty) = \bigcap_{m\in \N} \Gamma(f, m)$.

We have the following symplectic analogues of propositions~\ref{p.lower},
\ref{p.global} and \ref{p.jump}:

\begin{proposition}  \label{p.lower sympl}
Let $f \in \Sympl$, $\eps_0>0$, $\delta>0$, and $0< \kappa <1$.
If $m \in \N$ is sufficiently large,
then there exists a measurable function $N: \Gamma^*(f,m) \to \N$
such that for a.e. $x\in \Gamma^*(f,m) $ and
every $n \geq N(x)$ there exists a $(\eps_0, \kappa)$-realizable sequence
$\{\hL_{0},\ldots ,\hL_{n-1}\}$
at $x$ of length $n$ such that
$$
\frac{1}{n} \log \| \mathord{\wedge}^q ( \hL_{n-1} \cdots \hL_{0} ) \|  \, \leq \,
\Lambda_{q-1}(f,x) + \delta.
$$
\end{proposition}

\begin{proposition} \label{p.global sympl}
Let $f \in \Sympl$, $\eps_0>0$  and $\delta>0$.
Then there exist $m\in\N$ and a diffeomorphism $g \in \UU(f, \eps_0)$
that equals $f$ outside the open set $\Gamma(f,m)$ and such that
$$
\int_{\Gamma(f,m)} \Lambda_q(g,x) \, d\mu(x) <
\delta +
\int_{\Gamma(f,m)} \Lambda_{q-1}(f,x)  \, d\mu(x).
$$
\end{proposition}

\begin{proposition}\label{p.jump sympl}
Given $f \in \Sympl$, let
$$
J(f) = \int_{\Gamma(f,\infty)} \lambda_q(f,x) \, d\mu(x).
$$
Then for every $\eps_0>0$ and $\delta>0$, there exists a
diffeomorphism $g\in\UU(f,\eps_0)$ such that
$$
\int_M \Lambda_q(g,x) \, d\mu(x) <
\int_M \Lambda_q(f,x) \, d\mu(x) - J(f) + \delta.
$$
\end{proposition}

The proofs of these propositions are exactly the same as
those of the corresponding results in section~\ref{s.t.vol},
in the following logical order:
\begin{center}
proposition \ref{p.geom sympl}   $\Rightarrow$
proposition \ref{p.lower sympl}  $\Rightarrow$
proposition \ref{p.global sympl} $\Rightarrow$
proposition \ref{p.jump sympl}.
\end{center}
Concerning the first implication,
notice that if $x\in \Gamma^*(f,m)$ then, by lemma \ref{l.sympl 2},
the spaces $E^+_x$ and $E^-_x$
(that correspond to positive and negative Lyapunov exponents)
are Lagrangian, so proposition~\ref{p.geom sympl} applies.

\subsection{Conclusion of the proof of theorems~\ref{t.sympl.continuity} and
\ref{t.sympl}}

\begin{proof}[Proof of theorem~\ref{t.sympl.continuity}]
Let $f \in \Sympl$ be a point of continuity of the map $\LE_q(\cdot)$.
By proposition~\ref{p.jump sympl}, $J(f) = 0$, that is,
$\lambda_q(f,x) = 0$ for a.e. $x\in \Gamma(f, \infty)$.
Let $x \in M$ be a regular point. If $\lambda_q(f,x) > 0$, we have
(if we exclude a zero measure set of $x$)
$x \notin \Gamma(f, \infty)$.
This means that there is a dominated splitting,
$T_{f^n(x)} M = E_n \oplus F_n$, $n\in\Z$ of index $q$,
along the orbit of $x$.
Then $E_n$ is the sum of
the Oseledets spaces of $f$, at the point $f^n x$, associated to the Lyapunov
exponents $\lambda_1(f,x)$, \ldots, $\lambda_q(f,x)$, and
$F_n$ is the sum of the spaces associated to the other exponents.
By part~2 of lemma~\ref{l.sympl 2}, the splitting
$T_{f^n(x)} M = E_n \oplus F_n$, $n\in\Z$ is hyperbolic.
\end{proof}

The next proposition is used to deduce
theorem~\ref{t.sympl} from theorem~\ref{t.sympl.continuity}.

\begin{proposition} \label{p.bowen}
There is a residual subset $\RR_2 \subset\Sympl$
such that if $f\in \RR_2$ then either $f$ is Anosov or
every hyperbolic set of $f$ has measure $0$.
\end{proposition}

\begin{proof}
This is a modification of an argument from~\cite{Mane2}. We use
the fact, proved in~\cite{Zehnder}, that $C^2$
diffeomorphisms are dense in the space $\Sympl$. Another key
ingredient is that the hyperbolic sets of any $C^2$ non-Anosov
diffeomorphism have zero measure. We comment on the latter near
the end.

For each open set $U\subset M$ with $\overline{U}\neq M$
and each $f\in\Diff$, consider the maximal $f$-invariant
set inside $\overline{U}$,
$$
\Lambda_f(U) = \bigcap_{n\in \Z}f^n(\overline{U}).
$$
For $\eps > 0$, let $D(\eps,U)$ be the set of diffeomorphisms
$f\in\Sympl$ such that at least one of the following properties
is satisfied:
\begin{itemize}
\item[(i)] There is a neighborhood $\UU$ of $f$ such that
$\Lambda_g(U)$ is not hyperbolic for all $g\in \UU$;
\item[(ii)] $\mu (\Lambda_f(U)) < \eps$.
\end{itemize}
Clearly, the set $D(\eps, U)$ is open. Moreover, it is dense. Indeed,
if $f$ does not satisfy (i) then there is $g$ close to $f$ such that
$\Lambda_g(U)$ is hyperbolic. Take $f_1 \in C^2$ close to $g$ in
$\Sympl$.
Then $\Lambda_{f_1}(U)$ is hyperbolic with measure zero, and so
$f_1 \in D(\eps, U)$. This proves denseness. Hence the set
$$
D(U)= \cap_{\eps>0} D(\eps,U)
\supset \{ f\in \Sympl; \; \Lambda_f(U) \text{ is hyperbolic}
\Rightarrow \mu(\Lambda_f(U)) = 0 \}
$$
is residual. Now take $\BB$ a countable basis of open sets of $M$
and let $\widehat{\BB}$ be the set of all finite unions of sets in
$\BB$. The set
$$
\RR_2 = \bigcap_{U\in \widehat{\BB},\, \overline{U}\neq M} D(U)
$$
is residual in $\Sympl$ and the hyperbolic sets for every
non-Anosov $f\in \RR$ have zero measure.

Finally, we explain why all hyperbolic sets of a $C^2$ non-Anosov
diffeomorphism have zero measure. This is well-known for
hyperbolic {\em basic\/} sets, see~\cite{Bowen}. We just outline
the arguments in the general case. Suppose $f$ has a hyperbolic
set $\Lambda$ with $\mu(\Lambda)>0$.  Using absolute continuity of
the unstable lamination, we get that
$\mu_u(W_\eps(x)\cap\Lambda)>0$ for some $x\in \Lambda$, where
$\mu_u$ denotes Lebesgue measure along unstable manifolds. By
bounded distortion and a density point argument, we find points
$x_k\in\Lambda$ such that $\mu_u(W_\eps(x_k)\setminus\Lambda)$
converges to zero. Taking an accumulation point $x_0$ we get that
$W^u_\eps(x_0)\subset\Lambda$. We may suppose that every point of
$\Lambda$ is in the support of $\mu | \Lambda$. In particular,
there are recurrent points of $\Lambda$ close to $x_0$. Applying
the shadowing lemma, we find a hyperbolic periodic point $p_0$
close to $x_0$. In particular, $W^s_\eps(p_0)$ intersects
$W^u_\eps(x_0)$ transversely. Using the $\lambda$-lemma we
conclude that the whole $W^u(p_0)$ is contained in $\Lambda$.
Define $\Lambda_0$ as the closure of the unstable manifold of the
orbit of $p_0$. This is a hyperbolic set contained in $\Lambda$,
and it consists of entire unstable manifolds. Hence,
$W^s(\Lambda_0)$ is an open neighborhood of $\Lambda_0$. Using
that $f$ preserves volume, we check that
$f(W^s_\eps(\Lambda_0))=W^s_\eps(\Lambda_0)$. This implies that
$W^s(\Lambda_0)=\Lambda_0$ and so, by connectedness, $\Lambda_0$
must be the whole $M$. Consequently, $f$ is Anosov.
\end{proof}

\begin{proof}[Proof of theorem~\ref{t.sympl}]
It suffices to take $\RR=\RR_1\cap\RR_2$ with $\RR_1$ a residual set
of continuity points of $f\mapsto\LE_q(f)$, and $\RR_2$ as in
proposition~\ref{p.bowen}.
\end{proof}

\section{Proof of theorem \ref{t.cocycle}} \label{s.t.cocycle}

Let $M$ be a compact Hausdorff space,
$\mu$ a Borel regular measure
and $f:M\to M$ a homeomorphism preserving the measure $\mu$.
Let also $G \subset \gldr$ be a closed group
which acts transitively on $\rp$.

The following result provides an analogue of proposition~\ref{p.geom}:

\begin{proposition} \label{p.geom cocycle}
Given $A\in C(M,G)$ and $\eps>0$,
if $m \in \N$ is large enough then the following holds:

Let $y \in M$ be a non-periodic point and
suppose it is given a non-trivial splitting $\R^d = E \oplus F$
such that
\begin{equation} \label{e.not dominated cocycle}
\frac{\| A^m(y)|_F \|}{\mm(A^m(y)|_E)} \geq \frac 12.
\end{equation}
Then there exists, for each $j=0,1,\ldots,m-1$, some
$L_j \in G$ with $\| L_j - A(f^j y)\|< \eps_0$,
and there are non-zero vectors $v \in E$ and $w\in A^m(y)(F)$
with $L_{m-1} \cdots L_0 (v) = w$.
\end{proposition}

\begin{proof}
Let $\eps_1 = \| A \|_\infty^{-1} \eps$, where $\| A \|_\infty = \sup_{x\in M} \|A(x)\|$.
Let $\alpha>0$, depending on $\eps_1$, be given by lemma~\ref{l.transitive}.
Let
$$
K = \max \{ 1/ \sin^2 \alpha,  \| A \|_\infty \| A^{-1} \|_\infty  \},
\quad
C = \frac{8K}{\sin^2 \alpha},
$$
and $m \geq 2C / \alpha$.
Now take $y$, $E$ and $F$ an in the statement.
For $j=0,1,\ldots, m-1$, indicate
$A_j = A(f^j y)$, $E_j = A^j(x)(E)$, $F_j = A^j(x)(F)$.
As before, we divide the rest of the proof into three cases:

\paragraph{First case:}
We assume that there exists $\ell \in \{0,\ldots,m\}$ such that
\begin{equation} \label{e.I c}
\ang (E_\ell, F_\ell) < \alpha.
\end{equation}
Fix $\ell$ as above and take $\xi \in E_\ell$, $\eta \in F_\ell$
such that $\ang (\xi, \ell) < \alpha$. Let $R \in G$ be such that
$\|R - I \| < \eps$ and $R(\R \xi) = \R \eta$. If $\ell < m$, then
we define $L_j$ as $L_\ell = A_\ell R$ and $L_j = A_j$ for $j \neq
\ell$. If $\ell = m$, then we define $L_j$ as $L_\ell = R A_\ell$
and $L_j = A_j$ for $j \neq m$. In either case, the sequence
$\{L_0, \ldots, L_{m-1} \}$ has the required properties.

\paragraph{Second case:}
Assume that there exist $k,\ell \in \{0,\ldots,m\}$ such that $k<\ell$ and
\begin{equation} \label{e.II c}
\frac{\| A_{\ell-1} \cdots A_k |_{F_k} \|}{\mm ( A_{\ell-1} \cdots A_k |_{E_k})} > K.
\end{equation}
Once more, this is similar to the second case in propositions~\ref{p.geom}
and \ref{p.geom sympl}. We leave it to the reader to spell-out the details.

\paragraph{Third case:}
We suppose that we are not in the previous cases, that is, we assume
\begin{equation} \label{e.not I c}
\text{for every $j \in \{0,1,\ldots,m\}$,}\quad
\ang(E_j, F_j) \geq \alpha.
\end{equation}
and
\begin{equation} \label{e.not II c}
\text{for every $i, j \in \{0,\ldots,m\}$ with $i<j$,} \quad
\frac{\|A_{j-1} \cdots A_i |_{F_i} \|}{\mm ( A_{j-1} \cdots A_i |_{E_i})} \leq K.
\end{equation}

Take unit vectors $\xi \in E_0$ and $\eta \in F_0$
such that
$$
\| A_{m-1} \cdots A_0 (\xi) \| = \| A_{m-1} \cdots A_0 |_{E_0} \|
\text{ and }
\| A_{m-1} \cdots A_0 (\eta) \| =\mm ( A_{m-1} \cdots A_0 |_{F_0}).
$$
Let $\xi_j = A_{j-1} \cdots A_0 (\xi)$,  $\eta_j = A_{j-1} \cdots A_0 (\eta)$
and $Y_j = \R \xi_j \oplus \R \eta_j$.
By the assumption~\eqref{e.not dominated cocycle}, we have
$\|A_{m-1} \cdots A_0 (\eta) \| / \|A_{m-1} \cdots A_0 (\xi) \| \geq 1/2$.
Also, using~\eqref{e.not II c}, we have that for each $j\in \{1,\ldots,m\}$,
$$
K \geq
\frac{\|A_{j-1}\cdots A_0 (\eta)\|}{\|A_{j-1}\cdots A_0 (\xi)\|} \geq
\frac{\|A_{m-1}\cdots A_0 (\eta)\|/\|A_{m-1}\cdots A_j\|}
     {\|A_{m-1}\cdots A_0 (\xi)\|/\mm(A_{m-1}\cdots A_j)} \geq
\frac{1}{2K}.
$$
This, together with~\eqref{e.not I c} and lemma~\ref{l.4} implies
that, for all $j\in\{1,\ldots,m\}$,
\begin{equation}\label{e.eccentricity c}
\frac{\|A_{j-1} \cdots A_0 |_{Y_0} \|}{\mm ( A_{j-1} \cdots A_0 |_{Y_0})} < C.
\end{equation}
Now assign orientations to the planes $Y_j$ such that each
$A_j |_{Y_j} : Y_j \to Y_{j+1}$ is orientation-preserving.
Let $P_j$ be the projective space of $Y_j$, with the induced orientation.
Let $v_j = \R \xi_j$ and $w_j = \R \eta_j \in P_j$
For each $z \in P_j$, let $[z] \in [0, \pi)$ be the oriented
angle between $z$ and $v_j$.
So $z \mapsto [z]$ is a bijection and $[z] \mapsto [A_j z]$ is monotonic.
If $L : Y_0 \to Y_j$ is any linear map then, by lemma~\ref{l.quase conforme},
\begin{equation} \label{e.distortion}
0 < [z_2] - [z_1] \leq \frac{\pi}{2} \ \Longrightarrow \
\frac{[Lz_2] - [Lz_1]}{[z_2] - [z_1]} \leq
\frac{2}{\pi} \cdot \frac{\|L\|}{\mm(L)}.
\end{equation}

We define directions $u_0 \in P_0, \ldots, u_m \in P_m$ by recurrence
as follows:
Let $[u_0] = 0$ and
\begin{equation} \label{e.def uj}
[u_{j+1}] = [A_j u_j] + \min \{ [w_{j+1}] - [A_j u_j], \alpha \}.
\end{equation}
Then, for each $j<m$,
$[A_j u_j] \leq [u_{j+1}] \leq  [w_{j+1}]$.
Therefore, defining $[z_j] = [(A_{j-1} \cdots A_0)^{-1} u_j]$, we have
$$
0 = [z_0] \leq [z_1] \leq \cdots \leq [z_m] \leq [w_0] < \pi.
$$
In particular, for some $i=0, \ldots, m-1$,
$[z_{i+1}] - [z_i] < \pi /m$.
Therefore, by~\eqref{e.eccentricity c} and~\eqref{e.distortion},
$$
[u_{i+1}] - [A_i u_i]  =
[A_{j-1} \cdots A_0 z_{i+1}] - [A_{j-1} \cdots A_0 z_i]
< 2C / m < \alpha.
$$
By~\eqref{e.def uj}, $[u_{i+1}] = [w_{j+1}]$.
We conclude that $[u_m] = [w_m]$.
Now for each $j$, let $R_j \in G$ be such that $\| R_j - I \| < \eps$ and
$R_j ( A_j u_j ) = u_{j+1}$.
Let also $L_j = R_j A_j$.
Then $L_{m-1} \cdots L_0 (v_0) = w_m$.
\end{proof}

Next we define sets $\Gamma_p(A,m)$, $\Gamma_p^*(A,m)$, $\Gamma_p^\sharp(A,m)$
for $p\in\{1,\ldots,d-1\}$ and $m \in\N$, in the same way as in
section~\ref{s.t.vol}, with the obvious adaptations.
Lemma~\ref{l.aperiodic} also applies in the present context.

\begin{proposition}  \label{p.lower cocycle}
Given $A \in C(M,G)$, $\eps>0$, $\delta>0$, and $p\in\{1,\ldots,d-1\}$,
if $m \in \N$ is sufficiently large then
there exists a measurable function $N: \Gamma_p^*(A,m) \to \N$
such that for a.e. $x\in \Gamma_p^*(A,m) $ and
every $n \geq N(x)$ there exist matrices
$\hL_0,\ldots ,\hL_{n-1} \in G$ such that
$\| \hL_j - A(f^j x)\|< \eps$ and
$$
\frac{1}{n} \log \| \wp ( \hL_{n-1} \cdots \hL_{0} ) \|  \, \leq \,
\frac{\Lambda_{p-1}(A,x) + \Lambda_{p+1}(A,x)}{2} + \delta.
$$
\end{proposition}
The proof is the same as proposition~\ref{p.lower}.

\begin{proposition} \label{p.global cocycle}
Let $A \in C(M,G)$, $\eps_0>0$, $p\in\{1,\ldots,d-1\}$ and $\delta>0$.
Then there exist $m\in\N$ and a cocycle $B \in C(M,G)$,
with $\| B - A \|_\infty < \eps_0$,
that equals $A$ outside the open set $\Gamma_p(A,m)$ and such that
$$
\int_{\Gamma_p(A,m)} \Lambda_p(B,x) \, d\mu(x) <
\delta +
\int_{\Gamma_p(A,m)} \frac{\Lambda_{p-1}(A,x) + \Lambda_{p+1}(A,x)}{2} \, d\mu(x).
$$
\end{proposition}

The proof of proposition~\ref{p.global cocycle} is not just an
adaptation of that of proposition~\ref{p.global}, because Vitali's
lemma may not apply to $M$. We begin by proving a weaker statement,
in lemma~\ref{l.L infinity}.
Let $L^\infty(M,G)$ denote the set of bounded measurable functions
from $M$ to $G$. Oseledets theorem also applies for cocycles in
$L^\infty(M,G)$.

\begin{lemma} \label{l.L infinity}
Let $A \in C(M,G)$, $\eps_0>0$, $p\in\{1,\ldots,d-1\}$ and $\delta>0$.
Then there exist $m\in\N$ and a cocycle $\widetilde{B} \in L^\infty(M,G)$,
with $\| \widetilde{B} - A \|_\infty < \eps_0 /2$,
that equals $A$ outside the open set $\Gamma_p(A,m)$ and such that
$$
\int_{\Gamma_p(A,m)} \Lambda_p(\widetilde{B},x) \, d\mu(x) <
\delta +
\int_{\Gamma_p(A,m)} \frac{\Lambda_{p-1}(A,x) + \Lambda_{p+1}(A,x)}{2} \, d\mu(x).
$$
\end{lemma}

\begin{proof}[Sketch of proof]
We shall explain the necessary modifications of the proof of proposition~\ref{p.global}.
The sets $Z^i$, $\hQ^i$ and $Q^i$ are defined as before.
In lemma~\ref{l.g}, the castles $U^i$ and $K^i$ become equal to $Q^i$
(as $\kappa$ and $\gamma$ were $0$).
We decompose each base $Q^i_\mathrm{b}$ into finitely many disjoint measurable sets
$U^i_k$ with small diameter.
In each tower with base $U^i_k$ we construct the perturbation $\widetilde{B}$ using
proposition~\ref{p.lower cocycle}, taking $\widetilde{B}$ constant in each floor.
The definitions of $N$ and $G^i$ are the same.
In lemma~\ref{l.big G} several bounds (those involving $\kappa$ or $\gamma$) become trivial.
Then one concludes the proof in the same way as before.
\end{proof}

\begin{proof}[Proof of proposition~\ref{p.global cocycle}]
Let $A$, $\eps_0$, $p$ and $\delta$ be as in the statement.
Let $m$ and $\widetilde{B}$ be given by lemma~\ref{l.L infinity}.
Let $N \in \N$ be such that
$$
\int_{\Gamma_p(A,m)} \frac {1}{N} \log \|\wp (\widetilde{B}^N (x))\| \, d\mu <
2\delta +
\int_{\Gamma_p(A,m)} \frac{\Lambda_{p-1}(A,x) + \Lambda_{p+1}(A,x)}{2} \, d\mu.
$$
Let $\gamma = N^{-1}\delta$. Using Lusin's theorem (see~\cite{Rudin}) and the
fact that $G$ is a manifold (see~\cite{Helgason}), one finds a continuous
$B: M \to G$ such that $B=\widetilde{B}=A$ outside the open set
$\Gamma_p(A,m)$, the norm $\| B - \widetilde{B} \|_\infty < \eps_0 /2$, and the set
$E = \{ x\in M; \; B(x) \neq \widetilde{B}(x) \}$ has measure $\mu(E) < \gamma$.
Let $G = \bigcap_{j=0}^{N-1} f^{-j}\big(\Gamma_p(A,m) \minus E \big) \subset \Gamma_p(A,m).$
Then $\mu \big( \Gamma_p(A,m) \minus G \big) \leq N \mu(E) < \delta$.
Then, letting $C$ be an upper bound for $\log \|\wp (\widetilde{B} (x))\|$, we have
\begin{align*}
\int_{\Gamma_p(A,m)} \Lambda_p(B,x) \, d\mu
&\leq
\int_{\Gamma_p(A,m)} \frac {1}{N} \log \|\wp (B^N (x))\| \, d\mu
\\
&<
C \delta + 2 \delta +
\int_{\Gamma_p(A,m)} \frac{\Lambda_{p-1}(A,x) + \Lambda_{p+1}(A,x)}{2} \, d\mu.
\end{align*}
Up to replacing $\delta$ with $\delta/(C+2)$, this completes the proof.
\end{proof}

Using proposition~\ref{p.global cocycle}, one concludes the proof of
theorem~\ref{t.cocycle} exactly as in subsection~\ref{ss.end t.vol}.
The fact that either vanishing of the exponents or dominance of the
splitting is also a sufficient condition for continuity is an easy
consequence of semi-continuity of Lyapunov exponents and robustness
of dominated splittings under small perturbations of the cocycle.



\bigskip

\begin{flushleft}
Jairo Bochi ({\tt bochi@impa.br}) \quad Marcelo Viana ({\tt
viana@impa.br})

\medskip

IMPA, Estrada D. Castorina 110, 22460-320 Rio de Janeiro, Brazil.
\end{flushleft}

\end{document}